\documentclass[12pts]{article}
\usepackage{latexsym,amssymb,amsmath} 
\usepackage[utf8]{inputenc} 
\usepackage{path}
\usepackage{empheq}
\usepackage{stmaryrd} 
\usepackage[dvipsnames]{xcolor}
\usepackage{tikz}
\usepackage{url}
\usepackage{verbatim}
\usepackage{cancel}
\usepackage{breqn}
\usepackage{latexsym,amssymb,amsmath} 
\usepackage{lmodern}
\usepackage{lipsum}
\usepackage{marvosym}
\usepackage{graphicx} 
\usepackage{float}
\usepackage{multicol}
\usepackage[font=normal, justification=centering]{caption}
\usepackage{algorithmic}
\usepackage{algorithm}
\usepackage{bm} 
\usepackage{natbib}
\usepackage{dsfont}
\usepackage{hyperref}
\usepackage{booktabs}
\usepackage{bbm}
\usepackage{enumitem} 
\usepackage{xspace}
\usepackage{xifthen}
\usepackage{authblk}   
\usepackage{booktabs,tabularx,ragged2e}
\usepackage{subcaption}
\usepackage{stackengine}
\usepackage{textcomp}
\usepackage{wasysym}

\definecolor{ocre}{RGB}{243,102,25}

\numberwithin{equation}{section}
\numberwithin{table}{section}
\numberwithin{figure}{section}

\newcommand{\emphara}[1]{\emph{#1}}

\newlength{\fwtwo} 

\newcommand{\pupif}[2]{\ifthenelse{\isempty{#1}}{\ensuremath{#2}}{\ensuremath{#2^{(#1)}}}}

\newcommand{\IR}{\mathbb{R}}

\newcommand{\E}{\mathrm{E}}

\DeclareMathOperator{\card}{card}

\newcommand{\xx}{\ensuremath{x}\xspace}
\newcommand{\uu}{\ensuremath{u}\xspace}
\newcommand{\ww}{\ensuremath{w}\xspace}
\newcommand{\junk}[1]{{}}
\newcommand{\EIsymb}{\text{EI}\xspace}
\newcommand\EI[2][]{{\pupif{#1}{\EIsymb}(#2)}\xspace}
\newcommand{\Xset}{\ensuremath{\mathcal{X}}\xspace}
\newcommand{\Xdoe}{\ensuremath{\mathbf{X}}\xspace}
\newcommand{\Uset}{\ensuremath{\mathcal{U}\xspace}}
\newcommand{\Udoe}{\ensuremath{\mathbf{U}}\xspace}
\newcommand{\Ydoe}{\ensuremath{\mathbf Y}\xspace}
\newcommand{\Lset}{\ensuremath{\mathcal{L}\xspace}}
\newcommand{\Ldoe}{\ensuremath{\mathbf L}\xspace}
\newcommand{\rlat}{\ensuremath{\ell}\xspace} 
\newcommand{\rlatvec}{\rlat\xspace} 
\newcommand{\rlatdoe}{\text{\textbf{L}}\xspace} 
\newcommand{\Lat}{\phi} 
\newcommand{\Latvec}[1]{\ensuremath{\pupif{#1}{\Lat}}\xspace} 
\newcommand{\ndis}{n_d}
\newcommand{\ncont}{n_c}
\newcommand{\nlev}[1]{m_{#1}}
\newcommand{\nlat}{q} 
\newcommand{\nlatvec}{\ensuremath{n_{\rlat}}\xspace} 
\newcommand{\ndoe}{\ensuremath{N_\text{DoE}}\xspace} 
\newcommand{\dualF}{\ensuremath{D\xspace}}

\title{	A comparison of mixed-variables Bayesian optimization approaches 
}
\author[1,3]{Jhouben Cuesta-Ramirez}
\author[3]{Rodolphe Le Riche}
\author[2]{Olivier Roustant}
\author[4]{Guillaume Perrin}
\author[5]{C\'edric Durantin}
\author[1]{Alain Gli\`ere}
\affil[1]{CEA, LETI, Univ. Grenoble Alpes, 38000 Grenoble, France}
\affil[2]{INSA Toulouse, France}
\affil[3]{LIMOS (CNRS, Mines Saint-Etienne, UCA), Saint-Etienne, France}
\affil[4]{COSYS, Université Gustave Eiffel, 77420 Champs-sur-Marne, France}
\affil[5]{CEA-DAM, France}
\begin{document}
	\DeclareGraphicsExtensions{.jpg}	
	\maketitle
	\begin{abstract}
		Most real optimization problems are defined over a mixed search space where the variables are 
		both discrete and continuous. 
		In engineering applications, the objective function is typically calculated with a numerically costly black-box simulation. 
		General mixed and costly optimization problems are therefore of a great practical interest, yet their resolution remains in a large part an open scientific question.
		
		In this article, costly mixed problems are approached through Gaussian processes where the discrete variables are relaxed into continuous latent variables. 
		The continuous space is more easily harvested by classical Bayesian optimization techniques than a mixed space would.
		Discrete variables are recovered either subsequently to the continuous optimization, or simultaneously with an additional continuous-discrete compatibility constraint that is handled with augmented Lagrangians. 
		
		Several possible implementations of such Bayesian mixed optimizers are compared. 
		In particular, the reformulation of the problem with continuous latent variables is put in competition with searches working directly in the mixed space. 
		Among the algorithms involving latent variables and an augmented Lagrangian, a particular attention is devoted to the Lagrange multipliers for which a local and a global estimation techniques are studied.
		
		The comparisons are based on the repeated optimization of three analytical functions and a beam design problem.
		
	\end{abstract}
	
	
	\section*{Introduction}
	\label{sec:intro}
	A key task in engineering design is to find an optimal configuration from a very large set of alternatives. 
	When the performance of the candidate solutions is measured through a realistic simulation, the numerical cost of the procedure becomes a bottleneck. 
	The optimization of computationally expensive simulators is a topic widely studied in the literature \cite{Thi2019}. \\
	
	In this work, we focus on Bayesian optimization (BO), which is particularly suitable for solving such problems \cite{frazier_tutorial_2018}. 
	Bayesian optimization is a sequential design strategy that requires a data-driven mathematical model or metamodel that provides predictions along with their uncertainty \cite{Bartz-Beielstein2019}. 
	The metamodel replaces some of the calls to the expensive simulation and is a key ingredient to the optimization of costly functions.
	An acquisition criterion \cite{wilson_maximizing_2018} aggregates the spatial predictions and uncertainties. 
	The metamodel is trained from a reduced set of simulation data and the acquisition criterion is maximized to propose new configurations to be simulated at the next iteration. 
	When the acquisition criterion is the expected improvement (\EIsymb), as first introduced in \cite{mockus1978application}, the BO algorithm is often called EGO (Efficient Global Optimization, \cite{Jones1998}).
	EGO is currently a state-of-the-art approach to medium size, continuous and costly optimization problems, both from an empirical \cite{leriche2021revisiting} and a theoretical point of view \cite{vazquez_convergence_2010}.

	\vskip\baselineskip
	However, in realistic settings, some of decision variables are categorical. In structural design for example, the type of material, the number of components, the choice between alternative technologies lead to discrete variables with no obvious distance between them. 
	The combination of continuous and categorical variables is called a mixed optimization problem.
	
	In non-costly cases, mixed optimization problems can be approached by Mixed-Integer NonLinear Programming \cite{belotti2013mixed} (when the discrete variables are integers), by sampling based techniques such as evolutionary optimization \cite{cao2000evolutionary,emmerich2008mixed,ocenasek2002estimation} or by alternating mixed programming \cite{audet2001pattern}.
	
	When the objective function is costly, mixed optimization problems remain challenging and a topic for research.
	It is customary to replace some of the calls to the original objective function by calls to a (meta)model of it. 
	\cite{bartz2017model} {provide} an overview of metamodels that have or can be used in optimization when the variables are continuous or discrete.
	Bayesian optimization methods, which rely on metamodels to save computations, have already been extended to mixed problems.
	It was made possible by the realization that GP kernels (covariance functions) in mixed variables can be created by composing continuous and discrete kernels.
	The acquisition function is defined over the same space as the objective function. Therefore maximizing the acquisition function is also a mixed variables problem.
	
	To the best of our knowledge, the first EGO-like algorithm for mixed variables has been proposed in \cite{hutter2011sequential}.
	In this article, the mixed kernel is a product of continuous and discrete Gaussian kernels, and random forests constitute an alternative choice of mixed metamodel. More precisely, the discrete kernel is a Gaussian of integer or hamming (also known as Gower) distance for ordinal or nominal variables, respectively. In \cite{hutter2011sequential}, the expected improvement is first optimized with a multi-start local search for both the continuous and discrete variables (thus a neighborhood for the discrete variables is defined) which is then complemented by a random search.
	This work was continued with the REMBO method in \cite{wang2016bayesian}, where a random linear embedding is introduced to tackle high-dimensional problems.
	Discrete variables were relaxed into continuous variables thanks to a mapping function. The optimization of the acquisition function was made with a combination of the DIRECT and CMA-ES continuous global optimizers.
	Both \cite{hutter2011sequential} and \cite{wang2016bayesian} have been motivated by applications to the automatic configuration of algorithms. The goal of reaching very high dimensions (millions) probably forced the authors to use isotropic kernels {as a way to keep the number of hyper-parameters low (only one length-scale for all dimensions)}.
	
	A Bayesian mixed optimizer is presented in \cite{pelamatti_efficient_2019}. The GP kernels are products of continuous and discrete kernels. Different discrete kernels are compared, namely the homo- and hetero-scedastic hypersphere decomposition and the compound symmetric kernels. The optimization of the acquisition function is performed with a genetic algorithm in mixed variables.
	A similar BO with mixed kernel is described in \cite{Zuniga_Sinoquet_mixedEGO}, but the expected improvement is optimized with the mixed version of the MADS algorithm \cite{audet2001pattern} and the neighborhood of the categorical variables is defined through a probabilistic model.
	
	Random forests can replace the kriging model in BO with mixed inputs as they natively have a measure of prediction uncertainty. 
	Such an implementation, first done in \cite{hutter2011sequential}, is part of the \texttt{mlrMBO} \texttt{R} package \cite{mlrMBO}, in conjunction with several acquisition criteria that can be optimized with a ``focus-search'' algorithm. The focus-search algorithm hierarchically samples the search space of the chosen acquisition criterion.

	\vskip\baselineskip
	Recent developments in metamodels involving mixed variables show that it is possible to map the categorical variables into quantitative non-observed \emph{latent variables} that are then considered as continuous \cite{zhang_latent_2019}.
	Whenever it is possible to write a model of the studied system, quantitative latent variables exist that describe the effects of the categorical variables. Typically, there are more latent variables than categorical variables.
	The existence of continuous latent variables can sometimes be established from the physics of the considered phenomena, 
	e.g. in material science \cite{zhang_bayesian_2020}. In structural mechanics for example, if the categorical 
	variable describes the shape and the material of an element load in flexion, its bending moment of inertia is a candidate latent variable. 
	Latent variables can emulate the properties of the original categorical variables, in particular within the metamodel, and 
	open the way to reasonings with continous quantities: the kernels of the Gaussian processes can be taken as continuous, 
	gradients and neighborhoods are naturally defined during the optimization. On the contrary, categorical variables and their 
	inherent lack of distance definition is the cause of complications in the kernel definition and in the optimization.
	
	\vskip\baselineskip
	This article presents a new Bayesian optimization algorithm for mixed variables called LV-EGO (for Latent Variable EGO). 
	Our contribution with respect to \cite{zhang_bayesian_2020} is that the continuity of the latent variables is also taken 
	advantage of during the optimization of the acquisition criterion. This implies that categorical variables must be 
	recovered from the continuous latent variables proposed by the optimizer, which creates a new ``pre-image'' problem.
	
	Section~\ref{sec:pb_statement} introduces the problem and the principles of Bayesian optimization. 
	In Section~\ref{sec:LVEGO}, several variants of LV-EGO are described. They differ in the handling of the relationship 
	between the categorical and the latent variables: the ``vanilla'' LV-EGO just recovers categorical variables after the 
	optimization while augmented Lagrangian versions account for the link during the optimization through constraints.
	Section \ref{sec:applications} presents a set of benchmarks comparing our method to other state-of-the-art techniques. 
	One of the benchmarks is a beam design problem and gives the opportunity to discuss the interpretation of the latent variables. 
	Finally, Section~\ref{sec:Concl} offers conclusions and perspectives to this work.
	
	\subsection*{Notations and abbreviations}
	\begin{multicols}{2}
		{\small (by alphabetical order)}
		\begin{itemize}[label={},itemsep=0mm]
\item {ALV~: Augmented Lagrangian latent Variable.}
			\item DoE~: Design of Experiment.
			\item $\dualF(),{\widehat{\dualF}()}$~: dual {and approximate dual} functions.
			\item MLE~: Maximum Likelihood Estimation
\item {$\epsilon$~: relaxation constant for the discreteness constraint.}
			\item E[$\cdot$]~: mathematical expectation
\item {EGO~: Efficient Global Optimization algorithm.}
			\item \EI[]{}, \EI[]{\xx,\rlatvec}, \EI[t]{\xx,\rlatvec}~: expected improvement (at point $(\xx,\rlatvec)$ and iteration $t$).
\item {ES~: Evolution Strategy optimization algorithm.}
			\item \Latvec{}~,~\Latvec{t}: vector of latent mapping functions stemming from MLE maximization (at iteration $t$), $\Latvec{t}:~\uu \in \Uset \rightarrow \Latvec{t}(\uu) \in \Lset$. 

			\item {$f(),f^{(t)}()$~: non-costly objective function to minimize,  based on the GP, typically $-\log(1+\EI[t]{})$.}
			\item $g(),{g^{(t)}()}$~: inequality ($\le 0$) or equality constraint function.

\item {GP~: Gaussian process.}
\item {$\lambda,\lambda_t$~: Lagrange multiplier.}
			\item \rlatvec~: vector of relaxed latent variables. They take a value in $\Lset$.    
\item {LV~: Latent Variable.}
\item {MK~: Mixed Kriging, a GP indexed on mixed variables.}
\item {MS~: Mixed Space formulation (as opposed to relaxed with latent variables).}
			\item $\nlev{},\nlev{j}$~: number of levels for all discrete variables or for the discrete variable $\uu_j$.
			\item $\ncont$~: number of continuous variables.
			\item $\ndis$~: number of discrete variables.
			\item $\nlatvec$~: total number of latent variables for all discrete variables, in this article $=\ndis \nlat$.
			\item \nlat~: number of latent variables per discrete variable, in this article $=2$.
\item {$\rho,\rho_t$~: penalty parameter.}
\item {RFO~: Random Forest Optimization algorithm.}
			\item $t, ~^{(t)}$~: number of calls to the expensive objective function, superscript for functions redefined at each iteration {(depending on the GP)}.
			\item \uu~: vector of discrete (ordinal or nominal) variables, $\in \Uset$.
			\item \Udoe~: set of the discrete part of already evaluated points, $\in \Uset^{t}$.
			\item \xx~: vector of continuous variables, $\in \Xset \subset \IR^{\ncont}$.
			\item \Xdoe~: set of the continuous part of already evaluated points, $\in \Xset^t \subset {\IR^{\ncont}}^t$.
\item {$y(,)$~: ``costly'' objective function to minimize, typically based on a numerical simulation.}
			\item \Ydoe~: current set of outputs of the evaluated points, $\in \IR^{t}$.
		\end{itemize}
	\end{multicols}
	
	\section{Problem statement and background}
	\label{sec:pb_statement}
	We consider the problem of minimizing a function $y(x, u)$ depending on a vector of continuous variables $x = (x_1, \dots, x_\ncont)$ and a vector of discrete variables $u = (u_1, \dots, u_\ndis)$, where each $u_i$ has $\nlev{i}$ levels encoded $1, \dots, \nlev{i}$. 
	We denote $\Xset$ the domain of definition for the continuous inputs, typically, after rescaling, the hypercubic domain $[0, 1]^\ncont$. Similarly, we denote $\Uset = \prod_{j=1}^\ndis \{1, \dots, m_j\}$ the domain of definition for the discrete inputs. 
$\Xset\times\Uset$ is the set of the mixed optimization variables. \\
	We focus on costly functions, meaning that each evaluation of $y$ is time-consuming, and we aim at minimizing $y$ with a tiny budget of evaluations.
	In this context, minimizing directly $y$ is hardly possible. An alternative is to use Bayesian optimization (BO). In BO approaches, there are two main ingredients: a Gaussian process (GP) serving as a fast proxy, often called metamodel, built from the current learning set, and a sampling criterion, often called acquisition criterion, used to update the learning set with a new data point computed with $y$. 
	A famous acquisition criterion is the expected improvement (EI).
	In that case, the BO approach is often called Efficient Global Optimization (EGO) algorithm. \\
	
	To be more precise, let $(\Xdoe,\Udoe) = \{(\xx,\uu)^{(1)}, \dots, (\xx,\uu)^{(t)} \} \in (\Xset\times\Uset)^t$ be a design of experiments (DoE), and $y_i = y(\xx^{(i)},\uu^{(i)})$ be the corresponding function evaluations ($i=1, \dots, t$).
	Let $y_\min = \min(y_1, \dots, y_t)$ be the current minimum.
	Let us now assume that $y$ is a particular realization of the GP $Y$ defined on $\Xset\times\Uset$. In that case, the EI criterion is defined by
	$$ \EI{\xx,\uu} = \E \left[ \max( y_{\min} - Y^t(\ww), 0) \right], ~(\xx,\uu)\in (\Xset\times\Uset), $$
	where $Y^t$ is the conditional GP knowing the observations:
	$$ Y^t \coloneqq Y \ \vert \ \{ Y((\xx,\uu)^{(1)}) = y_1, \dots, Y((\xx,\uu)^{(t)}) = y_t \}. $$	
	
	Notice that $\EI{\xx,\uu}$ is large when exploiting interesting area, 
that is to say when there is a good chance that $Y^t(\xx,\uu)$ is smaller than $y_{\min}$.  
This may occur when $\E[Y^t(\xx,\uu)]$ is close to $y_{\min}$, or when exploring unvisited areas, 
i.e. when the variance of $Y^t(\xx,\uu)$ is large compared to $(\E[Y^t(\xx,\uu)]-y_{\min})^2$. The idea of EGO is to evaluate $y$ at a new point maximizing the EI criterion until a stopping criterion is reached. See Algorithm \ref{alg:EGOgeneric} for a synthetic description of the EGO algorithm when the stopping criterion is a maximum number of evaluations of $y$, noted $\textit{budget}$. 
{A maximum budget is the logical stopping criterion in our context of costly optimization. Other stopping conditions are possible in the form of lower bounds on the acquisition criteria (expected improvement, knowledge gradient \cite{frazier_tutorial_2018},\ldots) i.e., minimal measures of progress below which the search should stop. }
{In line 9, the solution returned by the algorithm is the best point of the last DoE, $(\Xdoe,\Udoe)$.}
	\begin{algorithm}[H]
		\begin{algorithmic}[1]
			\STATE Generate the initial DoE of size \ndoe, $(\Xdoe,\Udoe)$, and calculate $\Ydoe = (y_1, \dots, y_\ndoe)$, $t\leftarrow \ndoe$.
			\WHILE {$t ~ \leq ~ \text{budget}$}
			\STATE Estimate the GP $Y^{t}$ from the learning set formed by $(\Xdoe,\Udoe)$ and $\Ydoe$.
			\STATE Look for the current minimum $y_{\min}$ and maximize $(\xx,\uu)\mapsto \EI{\xx,\uu}$ on $\Xset \times \Uset$: $(\xx^{t+1},\uu^{t+1}) \in \text{argmax}_{(\xx,\uu) \in \Xset\times\Uset} \EI{\xx,\uu} $.
			\STATE Evaluate $y$ at $(\xx^{t+1},\uu^{t+1})$, $y^{t+1} = y(\xx^{t+1},\uu^{t+1})$.
			\STATE Update the learning set: $(\Xdoe,\Udoe) \leftarrow (\Xdoe,\Udoe) \cup (\xx^{t+1},\uu^{t+1})$, $\Ydoe \leftarrow \Ydoe \cup \{y^{t+1}\}$.
			\STATE $t \leftarrow t+1$
			\ENDWHILE
			\STATE $(\xx^\star,\uu^\star) = \arg \min_{(\xx,\uu) \in (\Xdoe,\Udoe)} y(\xx,\uu)$, $y^\star = y(\xx^\star,\uu^\star)$
			\RETURN ($\xx^\star,\uu^\star,y^\star$)
		\end{algorithmic}
		\caption{EGO algorithm on a mixed space}
		\label{alg:EGOgeneric}
	\end{algorithm}

	This EGO algorithm has been intensively studied to minimize  nonlinear functions that are expensive to be evaluated in the case $\Uset=\emptyset$, i.e. when all input variables are continuous (see \cite{leriche2021revisiting} for numerical illustrations of its efficiency). 
The application of this algorithm in the presence of categorical variables is much less documented (see e.g. \cite{pelamatti_efficient_2019, Zuniga_Sinoquet_mixedEGO}), which can be explained by two main difficulties. 
The first one is related to the difficult estimation of covariance kernels on mixed spaces. 
Indeed, multi-dimensional covariance functions are often built by combination of one-dimensional ones. 
Therefore, covariance functions on mixed spaces can be obtained by combining covariance functions on $\Xset$ and $\Uset$:
	\begin{equation} \label{eq:kernelCombination}
	\text{Cov}(Y(\xx,\uu),Y(\xx',\uu'))=k^x_1(x_1,x_1')\ast\cdots\ast k^x_{n_c}(x_{n_c},x_{n_c}')\ast k^u_1(u_1,u_1')\ast\cdots\ast k^u_{n_d}(u_{n_d},u_{n_d}'),
	\end{equation}
where $k^x_1,\ldots,k^x_{n_c},k_1^u,\ldots,k_{n_d}^u$ are covariance functions and $\ast$ is an operation that preserves positive definiteness, such as sum or product. If we focus on the single categorical variable $u_j$ with levels $1, \dots, \nlev{j}$, we can identify the covariance function $k_j^u$ to a $(m_j\times \nlev{j})$-dimensional positive semidefinite matrix $\mathbf{T}$, such that for all $1\leq k,\ell\leq \nlev{j}$,
	
	\begin{equation}
	(\mathbf{T})_{k\ell}=k_j^u(k,\ell).
	\end{equation}
	
	This means that $\sum_{j=1}^{n_d}\nlev{j}(\nlev{j}+1)/2$ coefficients need to be estimated to determine a covariance on $\Uset$ in the general case. That number can be large when $m$ is large, which very often makes this estimation very difficult in practice. Furthermore, the optimization problem is often harder than the box-constrained one met with continuous variables. Indeed it is either constrained by the positive definiteness of $\mathbf{T}$, which is non-linear, or defined on a manifold if $\mathbf{T}$ is parameterized in spherical coordinates. We refer to \cite{roustant_group_2020} for more details and other parsimonious representations of $k_j^u$, which can reduce but not totally fix these issues.
	The second reason that can explain the few number of direct applications of EGO algorithm on mixed space is related to the difficult maximization of the expected improvement, i.e. the search of the new input points where to call the function $y$, which are solutions of:
	
	\begin{equation}
	\max_{\xx,\uu \in \Xset \times \Uset} \EI {\xx,\uu} ~.
	\label{eq-EImixed}
	\end{equation}
	
	Indeed, classical optimization algorithms on continuous spaces usually try to exploit information related to the gradient of the function to be maximized, as well as notions of proximity in the space of the inputs. 
However, these two notions are difficult to exploit when dealing with categorical inputs, i.e. without any a priori ordering between the input instances. 
To circumvent this difficulty, a naive approach of resolution would consist in no longer considering a single maximization problem on $\Xset\times\Uset$, but the resolution in parallel of $\prod_{j=1}^{n_d}\nlev{j}$ maximization problems on $\Xset$, i.e. one problem per combination of instances of the categorical inputs $u$. 
	Such an approach is not tractable when the number of optimization problems to be solved becomes large, which has motivated the definition of heuristics, such as evolutionary algorithms \cite{li2013mixed,cao2000evolutionary,lin_hybrid_2018}, which seek to concentrate the searches only on the interesting instances of $u$. 
These approaches still rely on a large number of calls to the function to be optimized, and their convergence is not always easy to quantify. \\
	
	Because mixed optimization problems are difficult, an alternative approach is proposed in the rest of this paper. 
	It is based on the possibility to relax the discrete variables into continuous latent variables, therefore benefiting from the more efficient search mechanisms that exist in continuous spaces (e.g. gradients).

	\section{EGO with latent variables}
	\label{sec:LVEGO}
	
	\subsection{Latent variables}
	\label{sec:latentVar}
	For an easier handling of categorical inputs, it was proposed in \cite{zhang_latent_2019} to replace each categorical input $u_j$ by a vector of $q_j\geq 1$ continuous inputs with values in $\IR^{q_j}$, noted $\ell_j$. 
	To give an intuition of the underlying idea in the automotive domain, a category of lubricant may be determined by physical continuous features such as boiling temperature, viscosity, etc that act as latent variables. 
	In structural mechanics, the shape of a load carrying structure, which is categorical, has underlying continuous flexural and membrane moments that drive its behavior.
	This amounts to associating to the Gaussian process (GP) $Y$ a new GP $\widetilde{Y}$, such that for each instance $u$ of the categorical inputs there exists a particular value of $\ell \coloneqq (\ell_1,\ldots,\ell_{n_{d}})\in \Lset \subset \mathbb{R}^{q_1}\times \cdots \times \mathbb{R}^{q_{n_d}}$, which is called \textit{latent variable}, allowing us to write:
	
	\begin{equation}
	Y(x,u)  \stackrel{\tiny{\text{in law}}}{=} \widetilde{Y}(x,\ell), \ x\in\Xset.
	\end{equation}
	
	An important point is that the values of $\ell$ are unobserved and therefore $\widetilde{Y}$ is unknown. Nevertheless, in order to replace the EI maximization problem on $\Xset \times \Uset$ by a new optimization problem on $\Xset \times \Lset$, a precise knowledge of $\widetilde{Y}$ is not necessary. 
	Indeed, assuming that kernels for mixed inputs are built by combining $1$-dimensional ones as in (\ref{eq:kernelCombination}), it is sufficient to identify the mappings $\phi_j$ from $\{ 1,\ldots,\nlev{j}\}$ to $\IR^{q_j}$ to each variable $u_j$ such that
	
	\begin{equation}  \label{eq:kerCatLV}
	k_{ j}^{u}(u_j, u_j') \approx k_j(\Latvec{}_j(u_j), \Latvec{}_j(u_j')),
	\end{equation}

	\noindent{}where $k_j$ is a continuous kernel on $\IR^{q_j} \times \IR^{q_j}$. Thus, it is not so much the values of $\Latvec{}_j(u_j)$ that are important, but their relative positions in $\IR^{q_j}$ in order to allow a reasonable reconstruction of the dependency structure between $Y(x,u)$ and $Y(x',u')$. \\
	
	According to the works achieved in \cite{zhang_latent_2019}, it appears that interesting mappings can be obtained by likelihood maximization and that relatively small values of $q_j$ can give a satisfying reconstruction. Following their recommendations, $q_j$ can be chosen equal to $1$ if $\nlev{j}\leq 3$ and to $2$ otherwise, which will be the values chosen in the rest of this paper. We denote by $n_{\ell}=\sum_{j=1}^{n_d}q_j$ the total number of latent variables. 
	Following \cite{roustant_group_2020}, the continuous kernel $k_j$ associated to the latent variables was chosen as the dot product kernel $k_j(t, t') = \langle t,  t' \rangle$.
	The corresponding covariance matrix is then low-rank, and provided better performances than the Gaussian kernel in the examples considered in the latter reference.
	
	This new parametrization leads us to the following adaptation of the EI maximization problem defined by Eq. (\ref{eq-EImixed}), which we name \textit{acquisition problem} as it allows to acquire a new point to evaluate:
	
	\begin{equation}
	\begin{split}
	& \max_{\xx,\rlatvec \in \Xset \times \Lset \subset \IR^{\ncont+\nlatvec}} \EI[t]{\xx,\rlatvec} \\
	& \text{such that } \exists \uu \in \Uset \text{ with } \rlatvec = \Latvec{t}(\uu).
	\end{split}
	\label{eq:EIx_lat} 
	\end{equation}
	
	Here, $\EI[t]{\xx,\rlatvec}$ is the expected improvement associated with GP $\widetilde{Y}$ at iteration $t$, $\Latvec{t}=(\Latvec{t}_1,\ldots,\Latvec{t}_{n_d})$ is the vector-valued mapping from $\prod_{j=1}^{n_d}\{1,\ldots,m_j\}$ to $\IR^{q_1}\times \cdots \times \IR^{q_{n_d}}$ at iteration $t$, and the constraint on the values of $\ell$ is driven by the fact that 	
	the values of the latent variables at the new point have to remain compatible with the current mapping functions. \\ 
	
	We follow two paths to solve this acquisition problem. 
	In the vanilla LV-EGO approach, which will be described soon, the \EIsymb maximization and the latent-discrete compatibility constraint are addressed one after each other. 
	Alternatively, with the augmented Lagrangian approaches, which will be described in Section \ref{augmentedEGOLag}, the full constrained optimization problem is treated.
	
	\subsection{The vanilla LV-EGO algorithm}		
	\label{vanillaEGO}
	
	At each iteration, the vanilla LV-EGO algorithm first maximizes \EIsymb in a relaxed, fully continuous, formulation where the discrete variables are replaced by relaxed continuous latent variables. 
	Then, a pre-image problem is solved where \EIsymb is maximized over the discrete variables only, the continuous variables being fixed at their value of the relaxed problem.
	The LV-EGO methodology is summarized in Algorithm \ref{alg:vanilla}. 
	\begin{algorithm}[H]
		\begin{algorithmic}[1]
			\STATE \textbf{Generate} the initial DoE of size \ndoe : \Xdoe, \Udoe
			\STATE \textbf{Costly function evaluations} $y(\xx^i,\uu^i) ~,~ i=1,\ldots,\ndoe$, $t\leftarrow \ndoe$
			\WHILE {$t \le \text{budget}$}
			\STATE \textbf{Estimate} the latent variable mappings \Latvec{t} and the parameters of the continuous GP $\widetilde{Y}$.
			\STATE \textbf{Perform} one EGO iteration in the \emph{relaxed continuous} space :\\ 
			$(\xx^{t+1},\rlatvec^{t+1}) = \arg \max_{\xx,\rlatvec \in \Xset \times \Lset \subset \IR^{\ncont+\nlatvec}} \EI[t]{\xx,\rlatvec}$. 
			\STATE \textbf{Recover} the \emph{discrete pre-image} component $u^{t+1}$ as: 
			$u^{t+1} = \arg \max_{\uu \in \Uset} \EI[t]{\xx^{t+1},\Latvec{t}(\uu)}$. 
			\STATE \textbf{Update} the DoE with $(\xx^{t+1},\uu^{t+1})$ with output value $y(\xx^{t+1},\uu^{t+1})$. 
			\STATE $t \leftarrow t+1$
			\ENDWHILE
			\STATE \textbf{Return} $(\xx^\star, \uu^\star) = \arg \min_{\xx^t,\uu^t \in (\Xdoe,\Udoe)} y(\xx^t,\uu^t)$ 	
		\end{algorithmic}
		\caption{Vanilla LV-EGO with mixed inputs}
		\label{alg:vanilla}
	\end{algorithm}
	
	The main difference with the generic Bayesian algorithm \ref{alg:EGOgeneric} is the new discrete pre-image problem in line 6.
Notice that the pre-image is formulated in terms of the \EIsymb objective, as opposed to a more arbitrary distance like $\lVert \rlatvec^{t+1} - \Latvec{t}(\uu)\rVert$.
{Solving the pre-image in terms of the iterative figure of merit, the expected improvement, is meant to provide a gain in efficiency with respect to a pre-image minimizing an Euclidean distance between the map of a discrete level and the latent variables. In the particular situation where the latent variable coincides with the image of a discrete level, $\rlatvec^{t+1} = \Latvec{t}(\uu^{t+1})$, both approaches yield the same result since $\rlatvec^{t+1}$ is a maximizer of \EIsymb (see line 5 of Algorithm \ref{alg:vanilla}).
}
	\vskip\baselineskip
	In terms of implementation, the \EIsymb maximization (line 5) is done with the COBYLA algorithm, a gradient free non-linear optimization technique \cite{Powell1994}.
	Since COBYLA is a local optimizer and the \EIsymb is a multimodal function, the maximization is repeated (10 times, {which is more than the maximum dimension of the test cases studied in this article and more than the default -- 3 -- of the \texttt{kergp} package}) from randomly chosen initial points and the best result is kept.
	An exhaustive search is carried out for the \EIsymb maximization of the pre-image problem (line 6).
	
	A comparison of the numerical complexities of the vanilla LV-EGO (Algorithm~\ref{alg:vanilla}) and the generic EGO (Algorithm~\ref{alg:EGOgeneric}) shows that the cost of the latent variables is limited. 
	Let us consider that the discrete space can be searched essentially by enumeration in $\mathcal O(\card \Uset) = \mathcal O(\prod_{i=1}^\ndis \nlev{i})$ operations (where $\nlev{i}$ is the number of levels per discrete variable) while a continuous space can be searched more efficiently in linear time.
	At each iteration, the Bayesian algorithms of this paper have three steps: first a GP is learned, then an acquisition criterion (\EIsymb for now and an augmented Lagrangian later) is maximized and finally a pre-image problem is solved. In the vanilla LV-EGO algorithm, these steps take place at lines 4, 5 and 6 of Algorithm~\ref{alg:vanilla}, respectively. 
	Table~\ref{tab:numComplexities} summarizes the number of operations per step. 
	The number of operations for learning the GPs is proportional to the cube of the number of points evaluated ($t$) because of the inversions of the covariance matrices, times the number of (continuous) parameters of the GP for the likelihood maximization. 
	
	The two other steps, the acquisition and the pre-image, imply predictions by the GP in $t^2$ operations times a number of operations that depends on the specific algorithm. 
	Comparing in Table~\ref{tab:numComplexities} the column of the generic EGO with that of the vanilla LV-EGO, and assuming that for all $i$ $\nlev{i}=\nlev{}$ to keep the discussion simple, it can be seen that the latent variables induce a slight extra cost to be learnt. 
	When $\nlat=2$, which is our default here, this extra cost is $\ndis \times \nlev{i} \times t^3$ operations. {Setting} $q=1$ would not add any cost to the learning.
	An advantage, which comes from the sequential resolution of the mixed problem, occurs in the maximization of the acquisition criterion when $\ncont + \nlat \times \ndis \times \nlev{} < \nlev{}^\ndis \times \ncont$, at the cost of an additional pre-image problem to solve. 
	Thus, LV-EGO will be faster than a mixed EGO once the latent variables are estimated if $\nlev{}^\ndis + \ncont + \nlat \times \nlev{} \times \ndis < \nlev{}^\ndis \times \ncont$, which happens frequently (take for example $\ncont=4,\ndis=2, \nlev{}=10, \nlat=2$). 
	
	\begin{table}
		
		\begin{tabular}{|p{1.5cm}|c|p{3.2cm}|p{3.2cm}|p{3.2cm}|}
			\hline
			~ & Mixed space search & Vanilla LV-EGO & ALV-EGO-g & ALV-EGO-l \\
			~ & (Alg.~\ref{alg:EGOgeneric}) & (Alg.~\ref{alg:vanilla}) & (Alg.~\ref{alg:ALV-EGO}+\ref{alg:globDual}) & (Alg.~\ref{alg:ALV-EGO}+\ref{alg:locDual}) \\
			\hline
			GP learning & $(\ncont+\sum_{i=1}^\ndis \nlev{i})\times t^3$ & $(\ncont+\nlat \times  \sum_{i=1}^{\ndis} \nlev{i})\times t^3$ & $(\ncont+\nlat \times \sum_{i=1}^{\ndis} \nlev{i})\times t^3$ & $(\ncont+\nlat \times \sum_{i=1}^{\ndis} \nlev{i})\times t^3$ \\
			\hline
			max acquisition & $ (\prod_{i=1}^\ndis \nlev{i}) \times \ncont \times t^2$ & $(\ncont + \nlat \times \sum_{i=1}^{\ndis} \nlev{i}) \times t^2$ & $(N_{\text{DoE}}' + \ncont + \nlat \times \sum_{i=1}^{\ndis} \nlev{i}) \times t^2$ & $(\ncont + \nlat \times \sum_{i=1}^{\ndis} \nlev{i}) \times t^2$ \\
			\hline
			pre-image & 0 & $(\prod_{i=1}^{\ndis} \nlev{i}) \times t^2$ & $(\prod_{i=1}^\ndis \nlev{i}) \times t^2$ & $ (\prod_{i=1}^\ndis \nlev{i}) \times t^2$ \\
			\hline
		\end{tabular}
		\caption{Numerical complexities of the algorithms compared at each iteration (for a given $t$).}
		\label{tab:numComplexities}
	\end{table}
	
	\subsection{LV-EGO algorithms with Augmented Lagrangian}
	\label{augmentedEGOLag}
	
	A possible pitfall of the vanilla LV-EGO detailed in Algorithm \ref{vanillaEGO} is that the link between the discrete variables \uu and their relaxed continuous counterparts \rlatvec is lost when maximizing $\EI[t]{\xx,\rlatvec}$ in line 5.
	Recovering it during the discrete pre-image problem where $\xx$ is fixed to a value optimal in the relaxed formulation but possibly non-optimal with respect to the mixed problem (\ref{eq-EImixed})  may yield a sub-optimal solution.
	For this reason, we now propose LV-EGO algorithms that account for the discreteness constraint during the optimization using augmented Lagrangians.\\

	In that prospect, notice that problem (\ref{eq:EIx_lat}) can be approximated as an optimization problem with an inequality constraint:
	\begin{equation}
	\begin{split}
	& \min_{\xx,\rlatvec \in \Xset \times \Lset \subset \IR^{\ncont+\nlatvec}} f^{(t)}(\xx,\rlatvec) \coloneqq -\log(1+\EI[t]{\xx,\rlatvec}) \\
	& \text{such that }  g^{(t)}(\rlatvec) \coloneqq \min_{\uu \in \Uset} \lVert \rlatvec - \Latvec{t}(\uu) \rVert - \epsilon ~\le~ 0
	\end{split}
	\label{eq:OP}
	\end{equation}
	where $\epsilon$ is a small positive relaxation constant and $\lVert \cdot \rVert$ the Euclidean norm.
	In this reformulation, called \textit{relaxed acquisition problem}, notice the $\log$ scaling of the \EIsymb which does not change the solution but improves the conditioning of the problem.
	Two values of $\epsilon$ will be discussed in the sequel, $\epsilon=0$ in which case the constraint becomes an equality constraint, 
	$\min_{\uu \in \Uset} \lVert \rlatvec - \Latvec{t}(\uu) \rVert = 0 $, 
	and $\epsilon>0$ but small which corresponds to a relaxation of the equality.
	In the sequel, $\epsilon$ is normalized with respect to the size of the vector of latent variables and set to $\epsilon=0.01$.

	The constrained optimization problem (\ref{eq:OP}) is solved through an
	augmented Lagrangian approach \cite{minoux1986mathematical,nocedal_numerical_2006}. 
	The augmented Lagrangian is that of Rockafellar \cite{rockafellar_lagrange_1993} which, specified for Problem~(\ref{eq:OP}), is,
	\begin{equation}
	L_A^{(t)}(\xx,\rlatvec;\lambda,\rho) =
	\begin{cases}
	f^{(t)}(\xx,\rlatvec) - \frac{\lambda^2}{2 \rho} &\text{ , if } g^{(t)}(\rlatvec) \le \frac{-\lambda}{\rho}{~,}\\
	f^{(t)}(\xx,\rlatvec)+ \lambda g^{(t)}(\rlatvec) + \frac{\rho}{2} g^{(t)}(\rlatvec)^2 &\text{ , otherwise{.}}
	\end{cases} 
	\label{eq-AL_Rockafellar}
	\end{equation}
	
	When $\epsilon=0$, the constraint $g^{(t)}(\rlatvec) \leq 0$ becomes an equality constraint, $g^{(t)}(\rlatvec)=0$. 
	In this case, the augmented Lagrangian connected to that of Rockaffelar is that of Hestenes \cite{hestenes1969multiplier} and takes the form 
	\begin{equation}
	L_A^{(t)}(\xx,\rlatvec;\lambda,\rho) = f^{(t)}(\xx,\rlatvec)+ \lambda g^{(t)}(\rlatvec) + \frac{\rho}{2} g^{(t)}(\rlatvec)^2
	\label{eq-AL_equality}
	\end{equation}
	
	Complementary explanations about the augmented Lagrangians are given in Appendix~\ref{sec-AL}.

	Augmented Lagrangians require to specify the values of the Lagrange multiplier, $\lambda$, and of the penalty parameter, $\rho$. 
	The general principle to fix them is to calculate the generalized Lagrange multiplier with a dual formulation \cite{minoux1986mathematical}:
	the dual function $\dualF^{(t)}$ is maximized with respect to the multiplier $\lambda$ while the penalty parameter $\rho$ should take the smallest value that allows {one} to find feasible solutions,
	\begin{equation}
	\begin{split}
	& \rho_t = \arg \min_{\rho \ge 0} \rho ~\text{ such that }~ g(\rlatvec^{t}) \le 0 \\
	& \text{where } \lambda_t = \arg \max_{\lambda \ge 0} \dualF^{(t)}(\lambda,\rho)~, \\
	& \quad \dualF^{(t)}(\lambda,\rho) = \min_{\xx,\rlatvec \in \Xset \times \Lset \subset \IR^{\ncont+\nlatvec}} L_A^{(t)}(\xx,\rlatvec;\lambda,\rho)~,\\
	& \quad \text{and } (\xx^t,\rlatvec^t) \in \arg \min_{\xx,\rlatvec \in \Xset \times \Lset \subset \IR^{\ncont+\nlatvec}} L_A^{(t)}(\xx,\rlatvec;\lambda,\rho)~.
	\end{split}
	\label{eq:dualFormulation}
	\end{equation}
	There are two logics to solve Problem~(\ref{eq:dualFormulation}), both of which have been investigated in this study.
	Following an idea presented in \cite{leriche_dual_evo_2002} for classical Lagrangians, we first propose to approximate the dual function $\dualF()$ as the lower front of the augmented Lagrangians of a finite set of calculated points.
	The approximated dual is 
	\begin{equation}
	\widehat{\dualF}(\lambda,\rho)  =  \min_{(\xx,\rlatvec)\in (\Xdoe',\rlatdoe')} L_A^{(t)}(\xx,\rlatvec;\lambda,\rho) 
	\label{eq-dualFHat}
	\end{equation}
	where $(\Xdoe',\rlatdoe')$ is a DoE that should not be mistaken for $(\Xdoe,\Udoe)$, the DoE of the original expensive problem. 
	$(\lambda_t,\rho_t,\xx^t,\rlatvec^t)$ comes from solving Problem~(\ref{eq:dualFormulation}) with minimizations over the finite set $(\Xdoe',\rlatdoe')$ instead of the initial $\Xset \times \Lset$.
	The functions in Problem~(\ref{eq:OP}) are not costly, $(\Xdoe',\rlatdoe')$ can be quite large.
	This approach is called \textit{global dual} as a global approximation to the dual function is built and maximized.
	It applies to very general functions, e.g., non differentiable functions.
	Another advantage of this approach is to allow large changes in the dual space. 
	Figure \ref{fig:ALg} provides an illustration of the approximated dual function and the effect of $\rho$ on the dual problem.
	The sketch is done for an inequality constraint, yet it also stands with marginal changes for an equality (cf. Appendix~\ref{sec-AL} and the caption to the Figure).
	Under the non-restrictive hypothesis that there is a $\rho$ beyond which the solution to the primal problem (\ref{eq:OP}) maximizes the dual function, maximizing the dual function preserves the global aspect of the search. However, the accuracy of the obtained $(\lambda_t,\rho_t)$'s will depend on the DoE. 
	Because there is only one constraint in the current problem and evaluating it does not require calling the costly function, the maximization on $\lambda$ and $\rho$ is done by enumeration on a $100 \times 20$ grid and $(\Xdoe',\rlatdoe')$ is a 100 LHS sample.

	\vskip\baselineskip
	The other path to updating the multiplier is to progressively change them based on the minimizers of the augmented Lagrangian at the 
	current step. This updating can be seen as a step in the dual space which makes it general, although it is usually proved by analogy 
	with the Karush Kuhn and Tucker optimality conditions \cite{nocedal_numerical_2006} which add unnecessary conditions (like differentiability), cf. Appendix~\ref{sec-AL}. Let $(\xx^t,\rlatvec^t)$ be a solution to 
	\begin{equation}
	\min_{\xx,\rlatvec \in \Xset \times \Lset \subset \IR^{\ncont+\nlatvec}} L_A^{(t)}(\xx,\rlatvec;\lambda_t,\rho_t)
	\label{eq-minLag}
	\end{equation}
	The update formula reads
	\begin{equation}
	\lambda_{t+1} ~=~\lambda_{t} + \rho_{t}\left( g^{(t)}(\rlatvec^t) + \max(0,\frac{-\lambda_t}{\rho_t}-g^{(t)}(\rlatvec^t))\right)
	\label{eq-lambdaUpdateKKT}
	\end{equation}
	As in \cite{picheny_bayesian_2016}, the penalty parameter $\rho$ is simply increased if the constraint is not satisfied,
	\begin{equation}
	\rho_{t+1} ~=~
	\begin{cases}
	\rho_t &\text{ if }\quad g^{(t)}(\rlatvec^t) \le 0 \\
	2 \rho_t &\text{ otherwise}
	\end{cases}
	\label{eq-rhoUpdateKKT}
	\end{equation}
	The update scheme based on equations~(\ref{eq-lambdaUpdateKKT}) and (\ref{eq-rhoUpdateKKT}) is called \textit{local dual} as a local step in the dual $(\lambda,\rho)$ space is taken.
	
	\begin{algorithm}[H]
		\begin{algorithmic}[1]
			\STATE \textbf{generate} the initial DoE of size \ndoe for $(\Xdoe,\Udoe)$
			\STATE \textbf{costly function evaluations} $y(\xx^i,\uu^i) ~,~ i=1,\ldots,\ndoe$, $t\leftarrow \ndoe$
			\STATE \textbf{initialize} budget, $\epsilon$
			\WHILE {$t \le \text{budget}$}
			\STATE \textbf{estimate} the latent variables \Latvec{t} and the GP parameters from current DoE.
			\STATE\COMMENT{approximately solve the relaxed acquisition problem~(\ref{eq:OP}) with $f^{(t)}(\cdot) = -\log(1+\EI[t]{\cdot})$} \\
			$(\xx^{t+1},\rlatvec^{t+1}) = \arg \min_{\xx,\rlatvec} f^{(t)}(\xx,\rlatvec)$ s.t. $g^{(t)}(\rlatvec^{t+1}) =  \min_{\uu \in \Uset} \lVert \rlatvec - \Latvec{t}(\uu) \rVert - \epsilon \le 0$, \\
			ALV-EGO-g variant: with the global dual scheme, cf. Algorithm~\ref{alg:globDual}  \\
			ALV-EGO-l variant: with the local dual scheme, cf. Algorithm~\ref{alg:locDual}
			\STATE \textbf{recover} the \emphara{discrete pre-image} component $\uu^{t+1}$ as: 
			$\uu^{t+1} = \arg \max_{\uu \in \Uset} \EI[t]{\xx^{t+1},\Latvec{t}(\uu)}$ 
			\STATE \textbf{update DoE}: add $(\xx^{t+1},\uu^{t+1})$ and its costly evaluation $y(\xx^{t+1},\uu^{t+1})$ to the DoE $(\Xdoe,\Udoe)$.
			\STATE $t \leftarrow t+1$
			\ENDWHILE
			\RETURN $(\xx^\star, \uu^\star) ~=~ \arg \min_{(\Xdoe,\Udoe)} y(\xx,\uu)$ 
		\end{algorithmic}
		\caption{Augmented Lagrangian Latent Variables EGO with global or local dual scheme (ALV-EGO-g or ALV-EGO-l)
			\label{alg:ALV-EGO}
		}
	\end{algorithm}
	
	Algorithm~\ref{alg:ALV-EGO} gathers all these changes and is called ALV-EGO.
	The essential difference between this ALV-EGO algorithm and the vanilla counterpart (Algorithm~\ref{alg:vanilla}) is that the \EIsymb maximization step is constrained so that the link between the discrete variables and the relaxed latent variables (hence the continuous $\xx$) is not lost and left to the pre-image step. The coupling between the continuous and the discrete variables is better accounted for. However, a pre-image step (line 7) is still necessary to fully recover a discrete solution in cases when the constraint is relaxed ($\epsilon>0$).
	In ALV-EGO like in the vanilla LV-EGO, there are $\nlat=2$ continuous latent variable per discrete variable. \\
	The global and local dual schemes are further detailed in Algorithms~\ref{alg:globDual} and \ref{alg:locDual}.
	The continuous minimizations of the Augmented Lagrangians once the Lagrange multipliers are set are always done with 10 random 
	restarts of the COBYLA algorithm \cite{Powell1994}. They occur in Algorithm~\ref{alg:globDual}, line 4 and Algorithm~\ref{alg:locDual} line 5. 
	To allow comparisons, this implementation is identical to the \EIsymb maximization of the vanilla LV-EGO (step 5 of Algorithm~\ref{alg:vanilla}).
	
	\begin{algorithm}[H]
		\begin{algorithmic}[1]
			\ENSURE An estimation of the solution to the relaxed acquisition problem~(\ref{eq:OP})
			\REQUIRE $f^{(t)}()$, an objective function, $g^{(t)}()$, a constraint \\
			$N'_{\text{DoE}}$, $N_\lambda$, $N_\rho~>0$
			\STATE Calculate a DoE $(\Xdoe',\Ldoe') \in (\Xset,\Lset)^{N'_{\text{DoE}}}$.\\
			Half of the points are feasible by \textit{i)} sampling a $\uu \in \Uset$ and \textit{ii)} setting $\rlatvec' = \Latvec{t}(\uu)$
			\STATE Create a grid of Lagrange multipliers and penalty parameters, $(\boldsymbol\lambda,\boldsymbol\rho)= \{\lambda_1,\ldots,\lambda_{N_\lambda}\}\times \{\rho_1,\ldots,\rho_{N_\rho}\}$, with $\lambda_i\ge 0$ and $\rho_j \ge 0$ for all $i,j$
			\STATE Approximately solve the dual problem by enumeration: \\ 
			$\rho_t$ smallest $\rho \in \boldsymbol\rho$ that yields a feasible solution, $g(\rlatvec^t) \le 0$ where \\
			\quad $(\lambda_t,\xx',\rlatvec') = \arg \max_{\lambda \in \boldsymbol\lambda} \min_{(\xx,\rlatvec) \in (\Xdoe',\Ldoe')} L_A^{(t)}(\xx,\rlatvec;\lambda,\rho)$
			\STATE Fine tune the next candidate: $(\xx^{t+1},\rlatvec^{t+1}) = \arg \min_{(\xx,\rlatvec) \in (\Xset,\Lset)} L_A^{(t)}(\xx,\rlatvec;\lambda_t,\rho_t) $
			\RETURN $\xx^{t+1}, \rlatvec^{t+1}$
		\end{algorithmic}
		\caption{Global dual scheme (makes ALV-EGO-g when used in Algorithm~\ref{alg:ALV-EGO})
			\label{alg:globDual}
		}
	\end{algorithm}
	
	\begin{algorithm}[H]
		\begin{algorithmic}[1]
			\ENSURE An estimation of the solution to the relaxed acquisition problem~(\ref{eq:OP})
			\REQUIRE $f^{(t)}()$, an objective function, $g^{(t)}()$, a constraint \\
			initial values of the Lagrange multiplier and penalty, $\lambda_\ndoe = 0$ and $\rho_\ndoe = 1$, $t$
			\IF{$t>\ndoe$}
			\STATE\COMMENT{when $t=\ndoe$ the initial $\lambda_\ndoe, \rho_\ndoe$ are used} \\
			Update $\lambda$ according to Eq.~(\ref{eq-lambdaUpdateKKT}) \\
			$\lambda_t = \lambda_{t-1} + \rho_{t-1}\left( g^{(t-1)}(\rlatvec^{t}) + \max(0,\frac{-\lambda_{t-1}}{\rho_{t-1}}-g^{(t-1)}(\rlatvec^{t}))\right)$ 
			\STATE Update $\rho$ according to Eq.~(\ref{eq-rhoUpdateKKT}) \\
			$\rho_{t} ~=~ \rho_{t-1}$ if $g^{(t-1)}(\rlatvec^{t}) \le 0$, $2 \rho_{t-1}$ otherwise
			\ENDIF
			\STATE $(\xx^{t+1},\rlatvec^{t+1}) = \arg \min_{(\xx,\rlatvec) \in (\Xset,\Lset)} L_A^{(t)}(\xx,\rlatvec;\lambda_t,\rho_t)$
			\RETURN $\xx^{t+1}, \rlatvec^{t+1}$
		\end{algorithmic}
		\caption{Local dual scheme (makes in ALV-EGO-l when used in Algorithm~\ref{alg:ALV-EGO})
			\label{alg:locDual}
		}
	\end{algorithm}
	While the local update of $\lambda$ and $\rho$ might seem less robust, it is the most common implementation and it might be sufficient for the constrained \EIsymb maximization. Indeed, between two iterations, the \EIsymb changes only locally around the current iterate. Providing the latent mapping functions do not change too much, a local update of $\lambda$ and $\rho$ seems appropriate.
	The numerical complexity of the ALV-EGO-g and -l algorithms is essentially the same as that of the vanilla LV-EGO, cf. Table~\ref{tab:numComplexities}. The global dual scheme has a slight extra-cost because of the search for the Lagrange multiplier and penalty parameter that require $N_{\text{DoE}}'$ extra GP predictions.
	
	Eventually, four variants of ALV-EGO are considered, ALV-EGO-ge or -gi or -le or -li where g stands for global, l for local, e for equality ($\epsilon=0$) and i for inequality ($\epsilon>0$).
	
	\section{Description of the numerical experiments}
	\label{sec:applications}
	
	\subsection{Algorithms tested}
	\label{sec:algoTested}
	
	The various algorithms tested are summarized in the Table \ref{tab:algos} which provides their names, the type of formulation for the mixed variables, the type of metamodel, the acquisition criterion and the technique to optimize the acquisition criterion.
	The two possible formulations for the mixed variables are either by searching in a mixed space (MS) or by a formulation in latent variables (LV).
	All Gaussian processes (GPs) are built with the \texttt{kerpg} package \cite{kergp}.
	The meaning of the acronyms is: LV-EGO, Latent Variables EGO; LV-RFO, Latent Variables Random Forest Optimization ; ALV-EGO-ge/-gi/-le/-li, Augmented Lagrangian Latent Variables global/local dual scheme with equality/inequality pre-image constraints; MS-RFO, Mixed Space search with Random Forest Optimization; MS-ES, Mixed Space search with Evolution Strategy; MS-MKES, Mixed Space search with Mixed Kriging metamodel and Evolution Strategy.
	\begin{table}
		\begin{tabular}{p{2cm}|p{1.8cm}|p{3cm}|p{1.8cm}|p{4cm} }
			\hline
			name & formulation & metamodel & acq. crit. & optimizer of the acq. crit. \\
			\hline
			LV-EGO & LV & GP & \EIsymb & restarted COBYLA \\
			\hline
			LV-RFO & LV & \texttt{randomForest} toolbox & \EIsymb & focus-search (from \texttt{mlrMBO}) \\
			\hline
			ALV-EGO-ge or -gi & LV & GP & \EIsymb & DoE (for $\lambda_t$ and $\rho_t$) and restarted COBYLA \\
			\hline
			ALV-EGO-le or -li & LV & GP & \EIsymb & restarted COBYLA \\
			\hline
			MS-RFO &  MS &  \texttt{randomForest} toolbox &  \EIsymb & focus-search (from \texttt{mlrMBO}) \\
			\hline
			MS-ES & MS & none & $-y(\xx,\uu)$ &  evolution strategy (from \cite{li2013mixed} in \texttt{CEGO} implementation \cite{cego}) \\
			\hline
			MS-MKES & MS & GP (sym. compound disc. kernel)& EI & evolution strategy (from \cite{li2013mixed} in \texttt{CEGO} implementation \cite{cego}) \\
			\hline
		\end{tabular} 
		\caption{
			Summary of the 9 algorithms tested: name, space over which it is defined (mixed versus continuous with latent variables), metamodel used, acquisition criterion, optimizer of the acquisition criterion.
			\label{tab:algos}
		}
	\end{table}
	The different algorithms will be tested on the suite of test problems described hereafter. 
	
	\subsection{Test cases}
	\label{sec:testCases}
	There are 3 analytical test cases and a beam bending problem. 
	The analytical test cases have all been designed by discretizing some of the variables of classical multimodal continuous test functions.
	The following notation is introduced to describe the discretization: if the continuous variable $\xx_i$ is discretized with $\uu_j$ that 
	takes values in $\{1,\ldots,\nlev{j}\}$, then $\uu_j(k)=\beta$ means $\xx_i=\beta$ when $\uu_j = k$, $\beta$ a scalar, $1\le k \le \nlev{j}$.
	
	\begin{figure}[H] 
		\begin{subfigure}[b]{0.5\linewidth}
			\centering
			\includegraphics[width=1\textwidth]{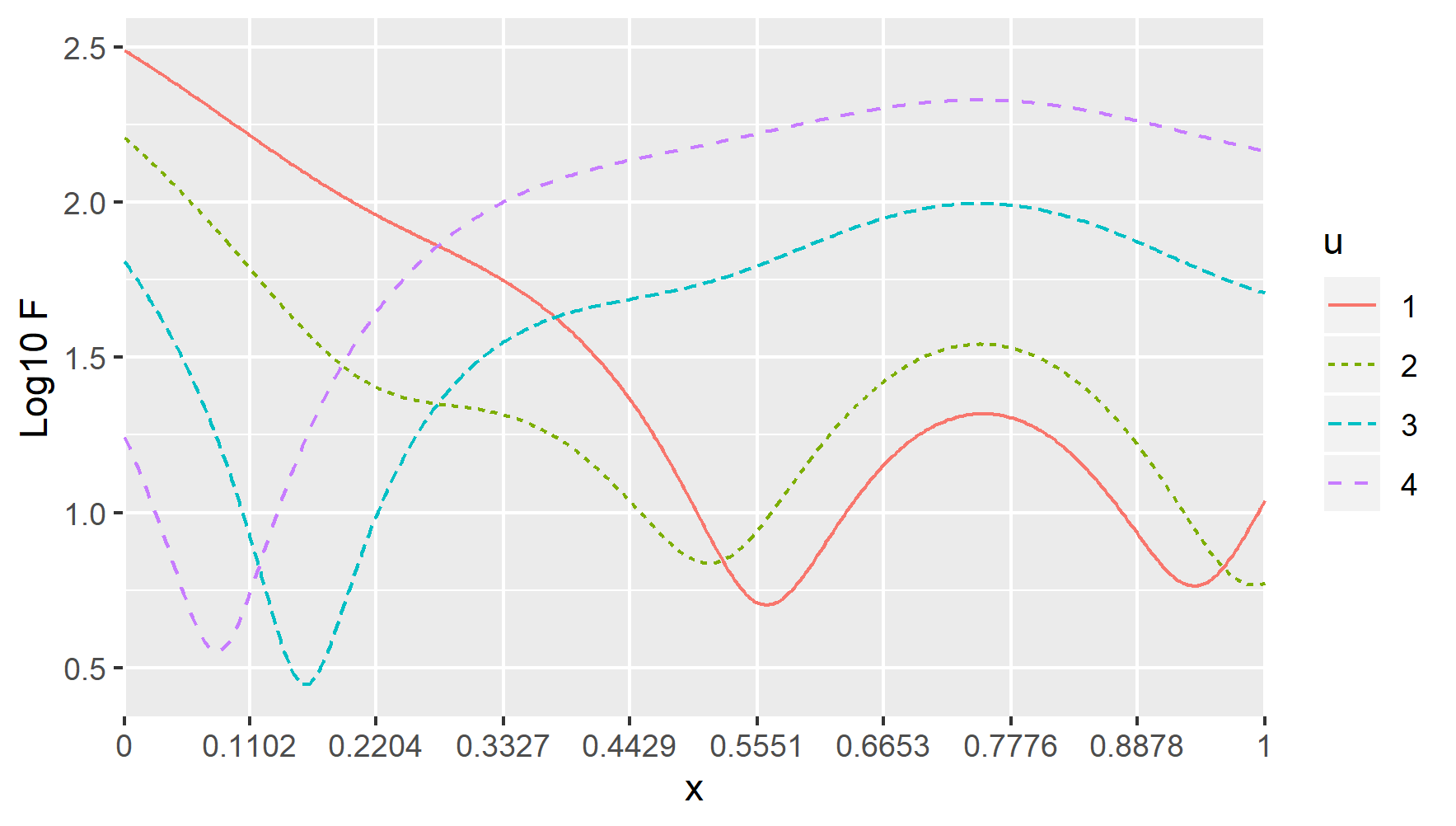}
			\caption{Discretized Branin-Hoo function}				
			\label{3dfig:bralev} 
		\end{subfigure}
		\begin{subfigure}[b]{0.5\linewidth}
			\centering
			\includegraphics[width=1\textwidth]{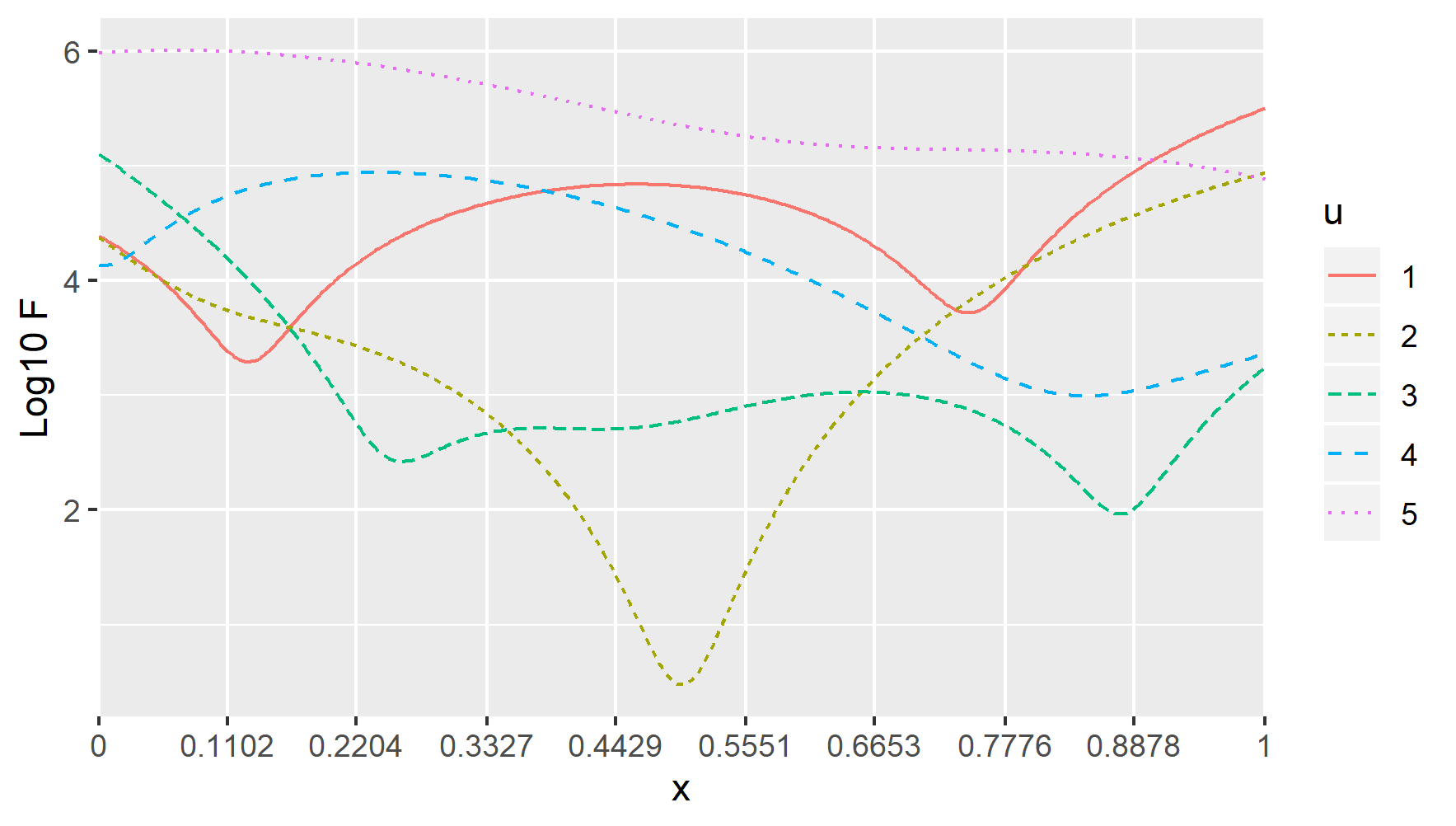}
			\caption{Discretized Goldstein-Price function}
			\label{3dfig:gslev} 
		\end{subfigure} 
		\caption{Two of the test functions with 1 discrete variable.
{Each curve is a 1-dimensional cross-section of the initial 2-dimension function, where the second variable is fixed at a given discrete value \uu.}}
		\label{1dfigs} 
	\end{figure}		
	
	\paragraph{Test case 1: discretized Branin function.}
	\label{sec:Branin}

	We modified the $2$ dimensional \textit{Branin-Hoo} function whose expression is
	\begin{align*}
	y(x_1,x_2) & = ({x'}_2 - b {x'}_1^2 + c {x'}_1 - r)^2 + s(1 - t) \cos{({x'}_1)} + s, \\
	x' & = {x'}^{\text{min}} + ({x'}^{\text{max}}-{x'}^{\text{min}}) \times x 
	\end{align*}
	where $b = 5/(4\pi^2), c = 5/\pi , r = 6, s = 10, t = 1/(8\pi)$, ${x'}^{\text{min}}=[-5;0] , {x'}^{\text{max}}=[10;15]$ 
	by keeping $x_1$ continuous in $[0;1]$ and making 
	$x_2$ discrete with $4$ levels $\{\uu(1) = 0; \uu(2) = 0.333; \uu(3)= 0.666; \uu(4) = 1\}$. 
	The discretized Branin, which was already used in \cite{zhang_bayesian_2020}, has several local minima as shown in Figure~\ref{3dfig:bralev}.
	The global optimum is located at $(\xx_1^\star,\uu^\star) = (0.182;\uu(3))$ with $y(\xx_1^\star,\uu^\star)=2.791$.
	\\
	
	\paragraph{Test case 2: discretized Goldstein function.}
	\label{sec:Goldstein}
	As a second test case, the continuous Goldstein function
	\begin{align*}
	y(x_1,x_2) = &[1+({x'}_1 + {x'}_2 + 1)^2  (19 - 14{x'}_1 + 3 {x'}_1^2  - 14{x'}_2 + 6{x'}_1 {x'}_2 + 3{x'}_2^2) ]  \\
	\times & [30 + (2{x'}_1 -3{x'}_2)^2  (18 - 32{x'}_1 + 12{x'}_1^2 + 48{x'}_2 - 36{x'}_1 {x'}_2 + 27 {x'}_2^2)] ~, \\
	x'  =& {x'}^{\text{min}} + ({x'}^{\text{max}}-{x'}^{\text{min}}) \times x 
	\quad,\quad  {x'}^{\text{min}} = [-2,-2]~,~{x'}^{\text{max}} = [2,2]
	\end{align*} 
	is partly discretized by replacing $\xx_2$
	by $\uu$ with $5$ levels $\{\uu(1) = 0; \uu(2) = 1/2; \uu(3)= 1/2; \uu(4) = 3/4 ; \uu(5) = 1\}$.
	The discretized Goldstein, which has also been studied in \cite{zhang_bayesian_2020}, is drawn in Figure~\ref{3dfig:gslev}.
	It has several local optima. 
	The global optimum is located at $(\xx_1^\star,\uu^\star) = (0.5; \uu(2))$ with $y(\xx_1^\star,\uu^\star)=3$.
	\\
	
	\paragraph{Test case 3: discretized Hartman function.}
	\label{sec:Hart}
	Two variables are discretized in the 6 dimensional Hartman function, 
	\begin{align*}
	y(x) &= -\sum_{i=1}^{4}\alpha_i \exp \left( -\sum_{j=1}^{d} A_{ij}(x_j-P_{ij})\right),
	\end{align*}
	where $x \in [0,1]^d$, $d = 6$, $\alpha= [1,1.2,3,3.2]^\top$  and
	\begin{align*}
	\hspace*{-.5cm}
	A = \begin{pmatrix}
	10& 3& 17& 3.5& 1.7& 8\\
	0.05& 10& 17& 0.1& 8& 14&\\
	3& 3.5& 1.7& 10& 17& 8 \\
	17& 8& 0.05& 10& 0.1& 14
	\end{pmatrix}, \,
	P= 10^{-4}\begin{pmatrix}
	1312& 1696& 5569& 124& 8283& 5886\\
	2329& 4135& 8307& 3736& 1004& 9991\\
	2348& 1451& 3522& 2883& 3047& 6650\\
	4047& 8828& 8732& 5743& 1091& 381
	\end{pmatrix}~.
	\end{align*}
{The variables} $x_5$ and $x_6$ are discretized with $5$ and $4$ levels respectively such that 
	$\{\uu_1(1) = 0.350; \uu_1(2) = 0.257; \uu_1(3)= 0.477; \uu_1(4) = 0.312; \uu_1(5) = 0.657\}$ and 
	$\{\uu_2(1) = 0.150; \uu_2(2) = 0.657; \uu_2(3)= 0.512; \uu_2(4) = 0.741\}$.
	Again, there are multiple local minima and the global optimum is 
	located at $(\xx^\star,\uu^\star) = (0.202;0.150;0.477;$ $0.275;\uu_1(4),\uu_2(2))$ with $y(\xx^\star,\uu^\star)=-3.322$.
	\\
	
	\paragraph{Euler-Bernoulli beam bending problem.}
	\label{sec:Bending}
	This test case corresponds to an horizontal beam that is clamped at one end and subject to a vertical force at the other end. 
	If the length of the beam is sufficiently long compared to the dimensions of its cross section, and if it is operating within its linear elastic range, the final beam deflection $y$ (to be minimized) is expressed as
	\begin{align}
	D(L,S,\tilde{I}) &= {\frac{P L^3 }{3 \, E\, S^2 \,\tilde{I}}}
	\label{eq:bending}
	\end{align}
	where {$P = 600 N$ is the vertical load, $E= 600 GP\!a$ is the Young's modulus,}
{$L \in [10,20]$} is the horizontal length of the beam, {$S \in [1,2]$} is the cross-section area and $\tilde{I} = I/S^2, \in \{\tilde{I}(1), \tilde{I}(2),\dots, \tilde{I}(12)\}$ is the normalized moment of inertia that can explicitly be derived for a given catalog of beam profiles. 
	The 12 levels of the normalized moment of inertia are 
	\begin{equation}
	\tilde{I} = \{0.083; 0.139; 0.380; 0.080; 0.133; 0.363; 0.086; 0.136; 0.360; 0.092; 0.138; 0.369 \}~.
	\label{eq:Iset}
	\end{equation} 
	We are interested in finding the best compromise between a minimization of the vertical deflection and the total weight, as expressed in the objective 
	\begin{align}
	y(\xx_1,\xx_2,\uu_1) & =  D(L,S,\tilde{I}) + \alpha L \, S ~, \\
	& \text{where  } L  = 10 + 10\times \xx_1 ~,~ S = 1 + \xx_2 ~,~ \uu_1 = \tilde{I} ~,~\\
	& \text{and } (\xx_1,\xx_2) \in [0,1]^2 ~.
	\label{eq:bendingVol}
	\end{align}
{Here $\alpha$ is the weight balancing the two effects in the objective function. 
It is chosen as $\alpha=60$ so that $y$ has several local minima and only one global minimum. 
This global} 
solution is $(\xx_1^\star,\xx_2^\star,\uu_1^\star) = (0; 0.43; \tilde{I}(3))$ with output $y^\star = 1.287385\times 10^3$.
	
	\subsection{Experiments setup and metrics} 
	\label{sec:doesetup}   
	The optimization of each pair of algorithm and test case are repeated 50 times from different initial DoEs.
	The DoEs are generated by minimax Latin Hypercube Sampling. 
	The size of the DoEs is $\ndoe = 4 \times \ncont \times \ndis \times \text{max}(\nlev{i})$ and a budget of $\ndoe$ + 50 evaluations of the true objective function. 
	Remember that the true objective function is supposed to be computationally intensive although it is not in these experiments so that runs can be repeated.
	The evolution strategies are stopped after $\ndoe + 50$ evaluations of the true function, like the other algorithms.

	The internal local optimizer, COBYLA, is restarted 5 times during the likelihood maximization and 10 times during the maximization of the acquisition criterion.
	The focus-search algorithm has a sample size of $1000$ with $5$ boundary reduction iterations and $3$ multi-starts, for a total of 3000 calls to the acquisition criterion. 
	
	A summary of the dimensions involved in the different examples is given in Table~\ref{tab:dimTestCases}.
	\begin{table}[!ht]
		\centering
		\begin{tabular}{|c|c|c|c|c|}
			\hline
			~ & $\ncont$ & $\ndis$ & $\nlev{i}$  & $\ndoe$ \\
			\hline
			Branin-Hoo & 1 & 1 & 4 & 16 \\
			\hline
			Goldstein & 2 & 1 & 5 & 40 \\
			\hline
			Hartmann & 4 & 2 & \{5,4\} & 160 \\
			\hline
			Beam Bending & 2 & 1 & 12 & 96 \\
			\hline
		\end{tabular}
		\caption{Dimensions and DoE size of the test cases.\label{tab:dimTestCases}}
	\end{table}
	
	\section{Results and discussion}
	\label{sec:Resul} 
	The results are provided with 4 main metrics.
	The performance of an algorithm is classically described by the median objective function over the 50 repeated runs, calculated at each iteration. 
	The associated measure of dispersion of the performance is the interquartile over the repetitions as a function of the iteration.
	To discriminate between methods that are rapid but provide rough solutions from the ones that take more time but yield better solutions, the two other metrics are based on the definition of targets.
	For each test case, a target is a given quantile of all the objectives functions found by all the algorithms throughout all the repetitions. A 10\% target is difficult, while a 50\% target is the median performance.  
	The third metric is the iteration number at which the median objective function of a given algorithm reaches a given target. 
	The fourth metric is the success rate (given a target), which is the percentage of the runs that do better than the target.
	The metrics associated to the quantile targets have the advantage that they are normalized with respect to the test cases: 
	thanks to the quantiles, the definitions of an easy, a median or a hard target stands accross the different functions to optimize. 
	The target-based metrics will later be averaged over the different test cases. 
	
	Let us now review the performances of the algorithms on each test case.
	
	\subsection{Analytical test functions}
	\label{sec:Rind}
	
\begin{figure}[H] 
		\begin{subfigure}[b]{0.5\linewidth}
			\centering
			\includegraphics[width=\textwidth]{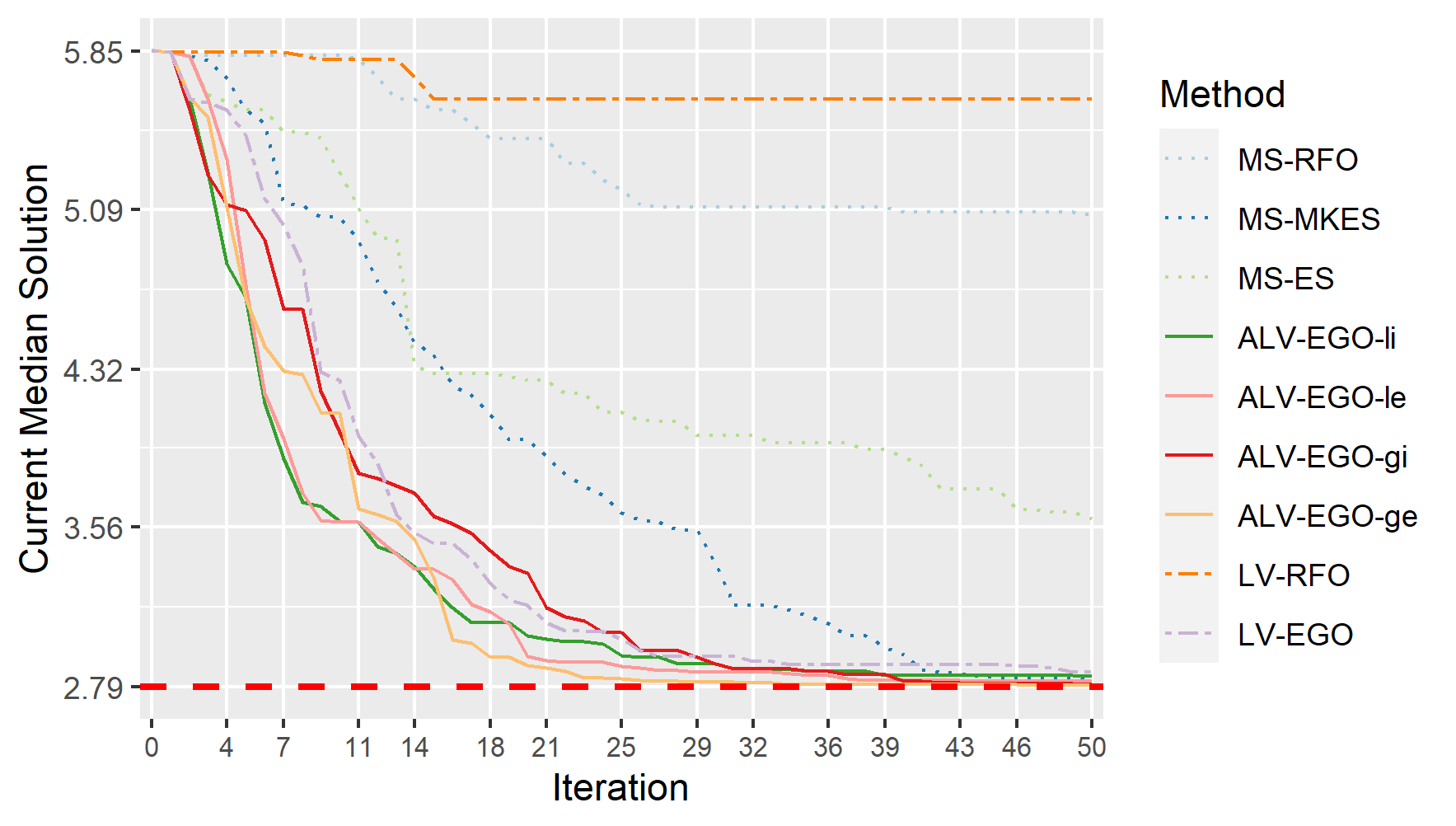}
			\caption{Median solution} 	
			\label{figUz0:ua} 
			\vspace{4ex}
		\end{subfigure}
		\begin{subfigure}[b]{0.5\linewidth}
			\centering
			\includegraphics[width=\textwidth]{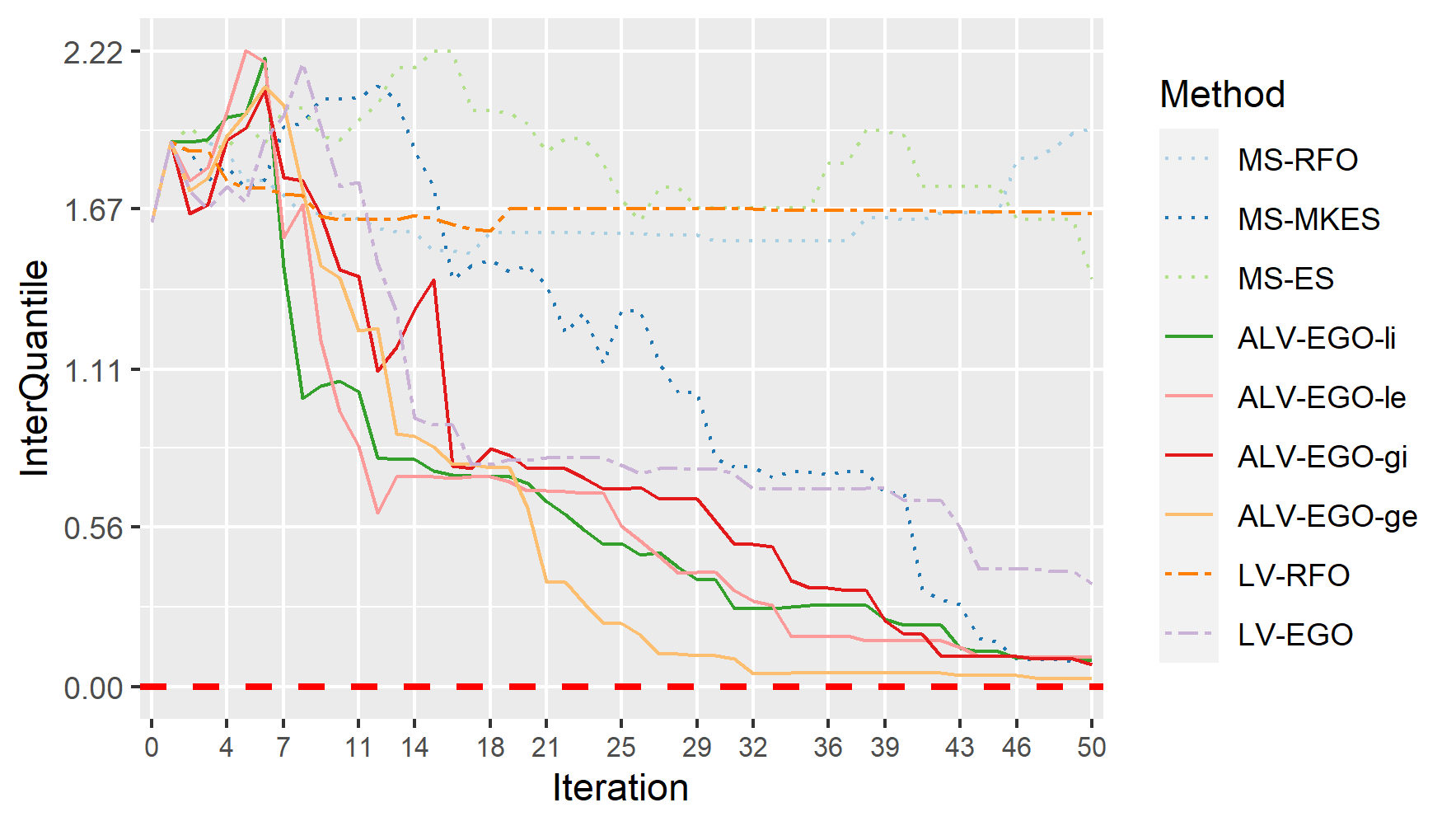}
			\caption{Interquartile range} 
			\label{figUz0:ub} 
			\vspace{4ex}
		\end{subfigure} 
		\begin{subfigure}[b]{0.5\linewidth}
			\centering
			\includegraphics[width=\textwidth]{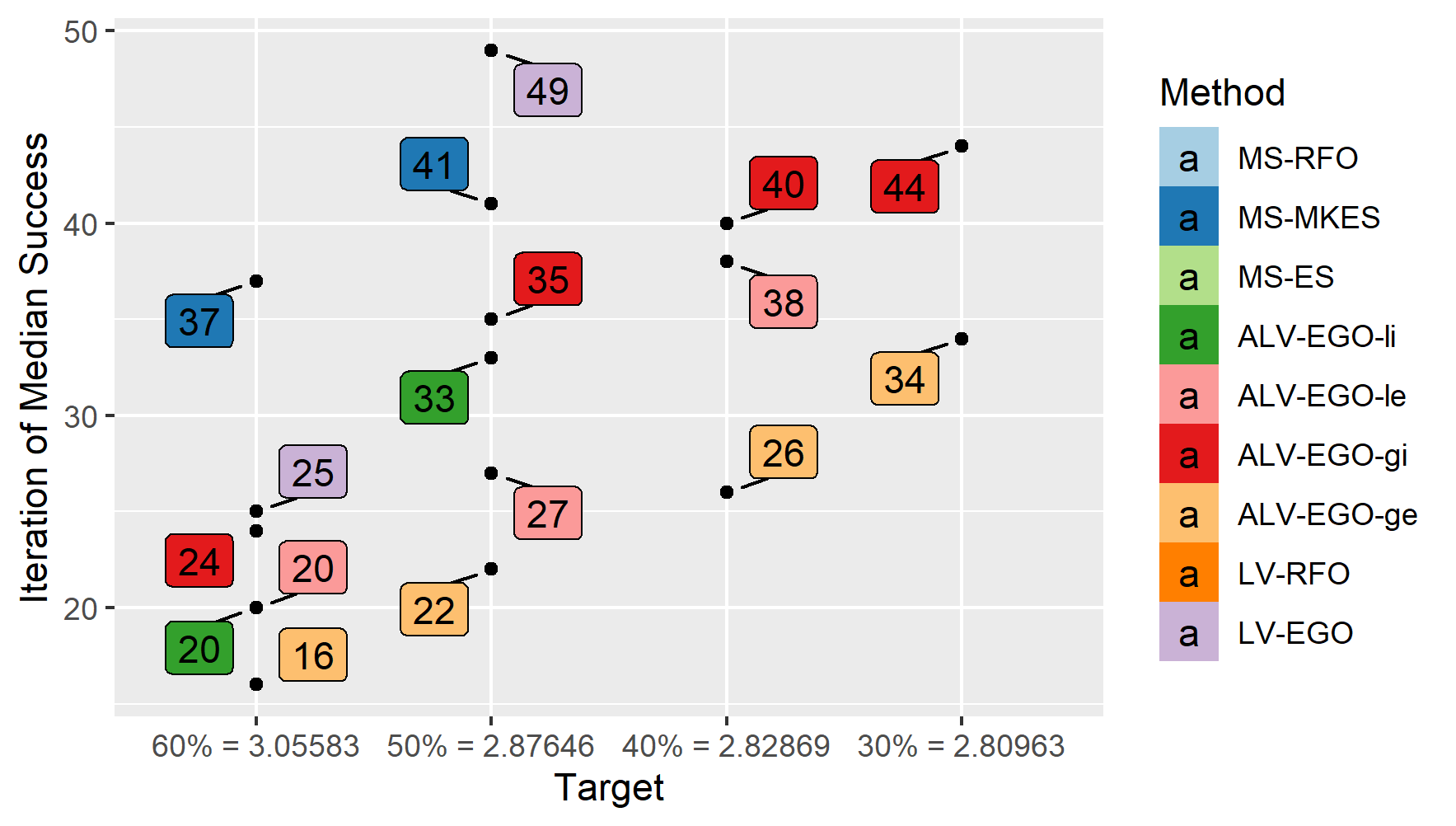}
			\caption{Iteration to median success} 			
			\label{figUz0:uc} 
		\end{subfigure}
		\begin{subfigure}[b]{0.5\linewidth}
			\centering
			\includegraphics[width=\textwidth]{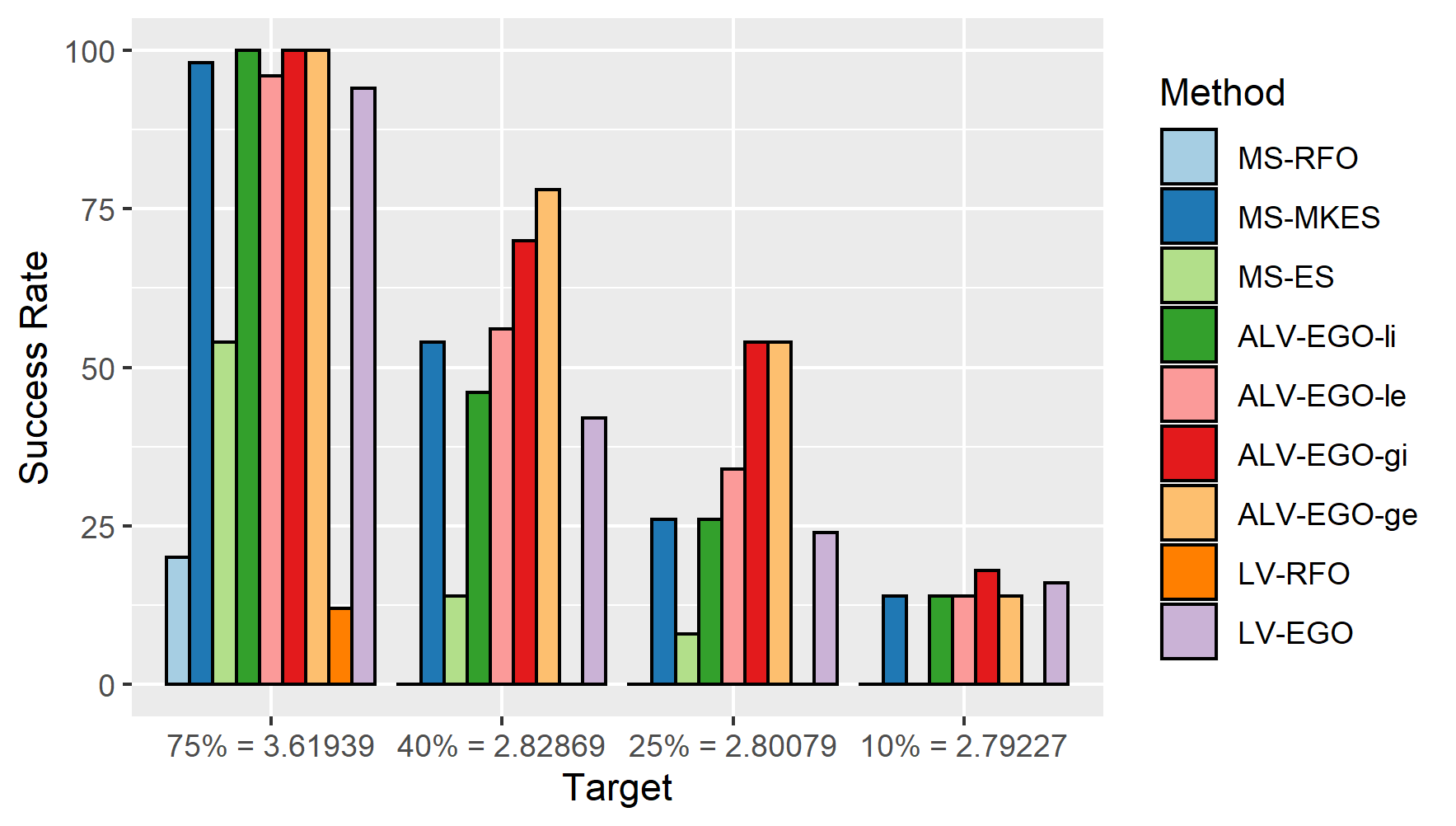}
			\caption{Success rate} 
			\label{figUz0:ud} 
		\end{subfigure} 
		\caption{Comparison of all 9 algorithms on the Branin function. $y^{\star} = 2.79118$. }
		\label{figUz0} 
	\end{figure}

	\paragraph{Branin function.}
	Figure \ref{figUz0} presents the results for the Branin function with the four metrics. 
	On the top left plot, showing the median value for the objective function, it is clear that the two methods that rely on the random forest metamodel (MS-RFO and LV-RFO) are overtaken by all other methods. 
	This indicates that, whether in the mixed or in the latent-augmented space, random forests do not represent sufficiently well the Branin function in comparison to Gaussian processes. 
	Looking at Figure~\ref{figUz0:ub}, it is observed that the fast methods typically have the lowest spread in performance and vice versa. This is expected as non converging runs may yield a wide range of performances.
	All methods involving the discrete constraint (i.e., the augmented Lagrangians) managed to improve over the LV-EGO performance; and including a mixed metamodel increased significantly the success rate  and the median solution for the evolutionary strategy.
	
	Regarding the success rate on Figure \ref{figUz0:ud}, the methods MS-MKES, LV-EGO, ALV-EGO-li, -le, -ge and -gi were the most prominent, the latter being capable to reach success rates of about $20\%$ for a $10\%$ target. 
	Notice that all these methods contain Gaussian processes. Indeed, the Branin function is easy to represent by a GP whether continuous or mixed.
	In the same vein, MS-MKES which differs from MS-ES by the use of a GP, clearly benefits from that metamodel. 
	\\
	All ALV- methods, which account for the discrete constraint, obtained the best median performances.
	ALV-EGO-ge in particular found all targets, in the median sense, earlier than the other algorithms as can be seen from Figure \ref{figUz0:uc}.
	
	A last comment is necessary regarding the bottom of Figure \ref{figUz0}: the plot on the left describes the median performance (in terms of targets reached) while the right plot counts the success rate at reaching a target over all runs. Therefore, some targets are reached on the right by some of the runs of a given algorithm, while they are never atteined on the left by the median of the same algorithm. This comment stands accross all test cases.
	
\begin{figure}[H] 
		\begin{subfigure}[b]{0.5\linewidth}
			\centering
			\includegraphics[width=\textwidth]{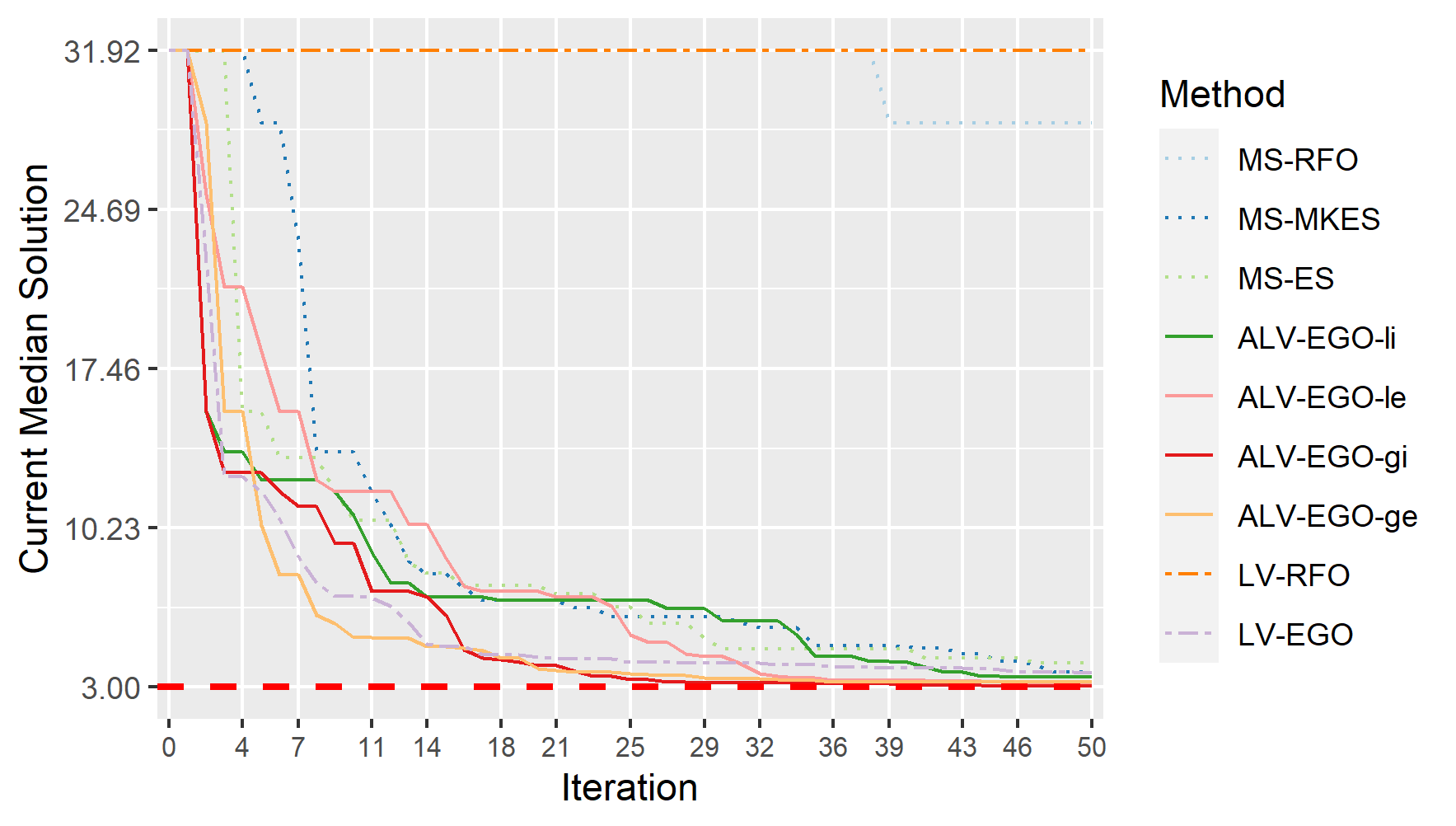}
			\caption{Median solution} 
			\label{figUg0:ua} 
			\vspace{4ex}
		\end{subfigure}
		\begin{subfigure}[b]{0.5\linewidth}
			\centering
			\includegraphics[width=\textwidth]{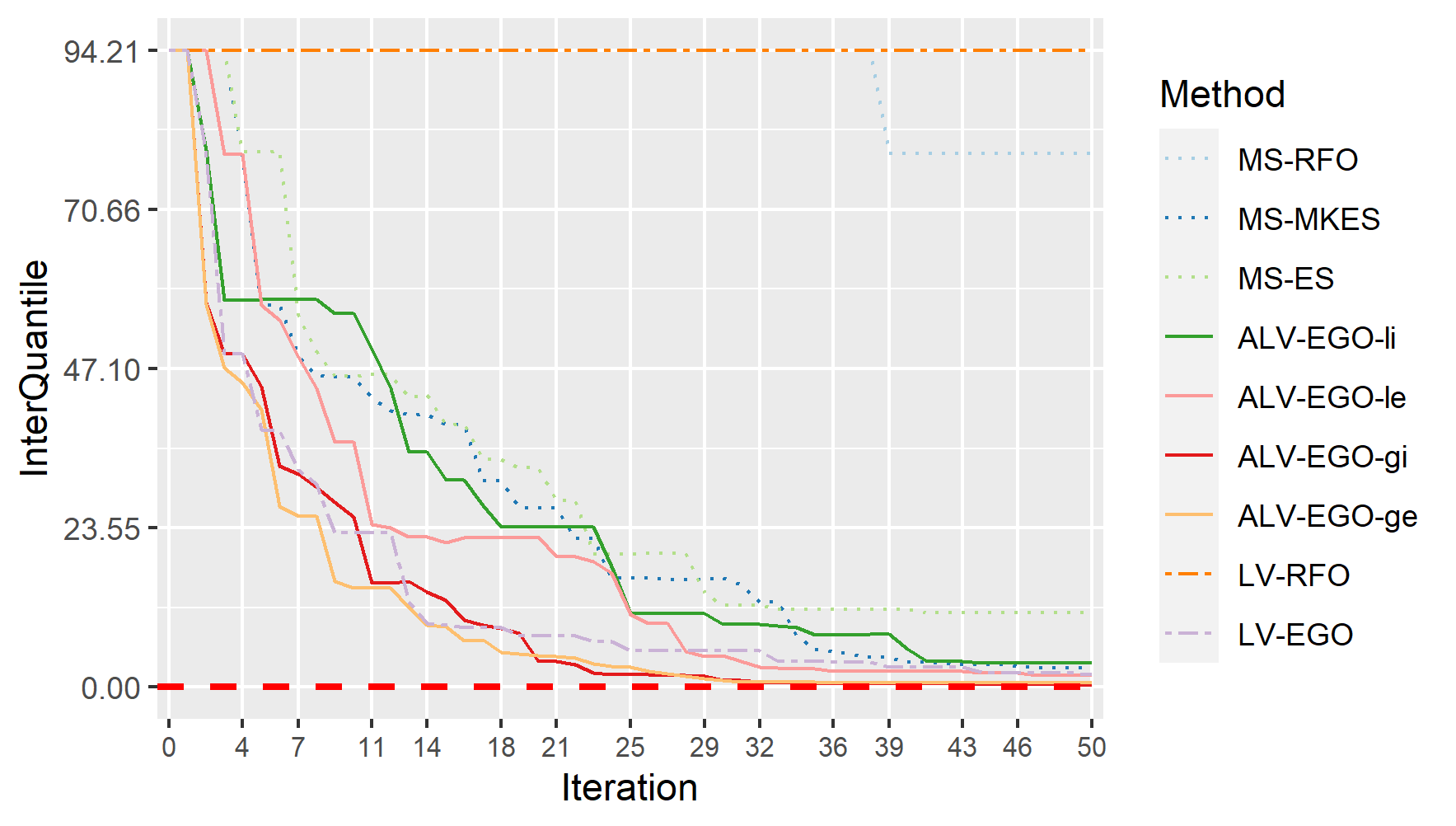}
			\caption{Interquartile range} 
			\label{figUg0:ub} 
			\vspace{4ex}
		\end{subfigure} 
		\begin{subfigure}[b]{0.5\linewidth}
			\centering
			\includegraphics[width=\textwidth]{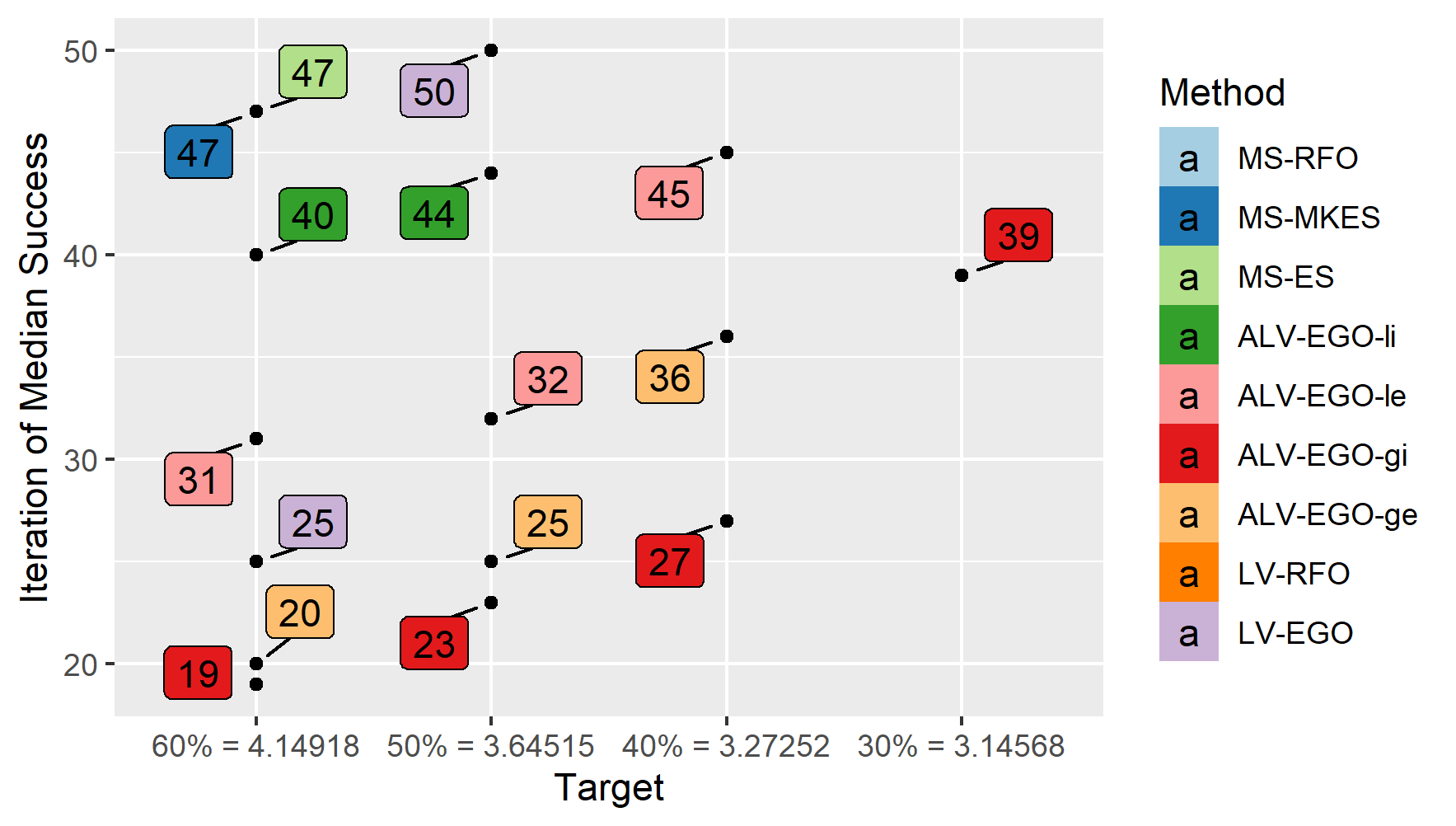}
			\caption{Iteration to median success} 
			\label{figUg0:uc} 
		\end{subfigure}
		\begin{subfigure}[b]{0.5\linewidth}
			\centering
			\includegraphics[width=\textwidth]{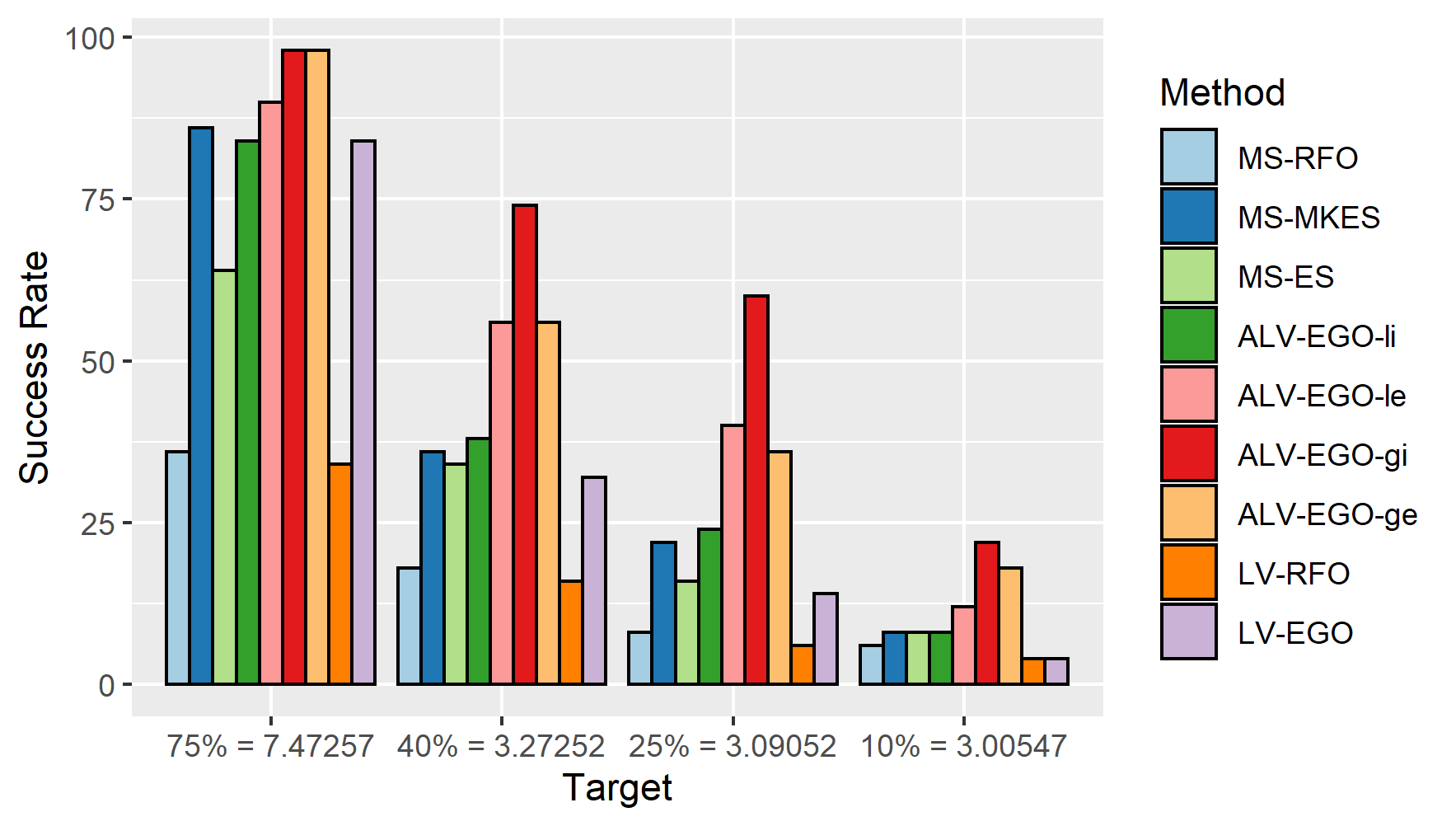}
			\caption{Success rate} 
			\label{figUg0:ud} 
		\end{subfigure} 
		\caption{Comparison of the 9 algorithms on the Golstein function. $y^{\star} = 3$.}
		\label{figUg0} 
	\end{figure}

	\paragraph{Goldstein function.}
	The experiments done with the Goldstein test function are summed up in Figure \ref{figUg0}. 
	Like with the Branin function, algorithms relying on random forests (LV-RFO and MS-RFO) showed both poor performance (top left plot). The associated high constant interquartile (top right) is that of the best points in the initial designs, which remains unchanged since no better point is found by these algorithms. 
	
	Considering the success rates for all targets (bottom plots), it is seen that accounting for the discreteness through a constraint (which is the distinctive feature of ALV- methods) is useful with the Goldstein function: like with Branin, ALV-EGO-gi is the best performer, but the other ALV- follow and outperform LV-EGO.
	All ALV- strategies almost reach the absolute target of percentile $25\%$ with a rate of $25\%$ or higher.
	The comparison of the plots~\ref{figUg0:uc} and \ref{figUg0:ud} also shows that, behind the ALV- methods, LV-EGO has a good median performance (cf. Figure~\ref{figUg0:uc}) but more of the MS-MKES searches manage to find difficult targets (the 25\% and 10\% quantiles).
	
\begin{figure}[H] 
		\begin{subfigure}[b]{0.5\linewidth}
			\centering
			\includegraphics[width=\textwidth]{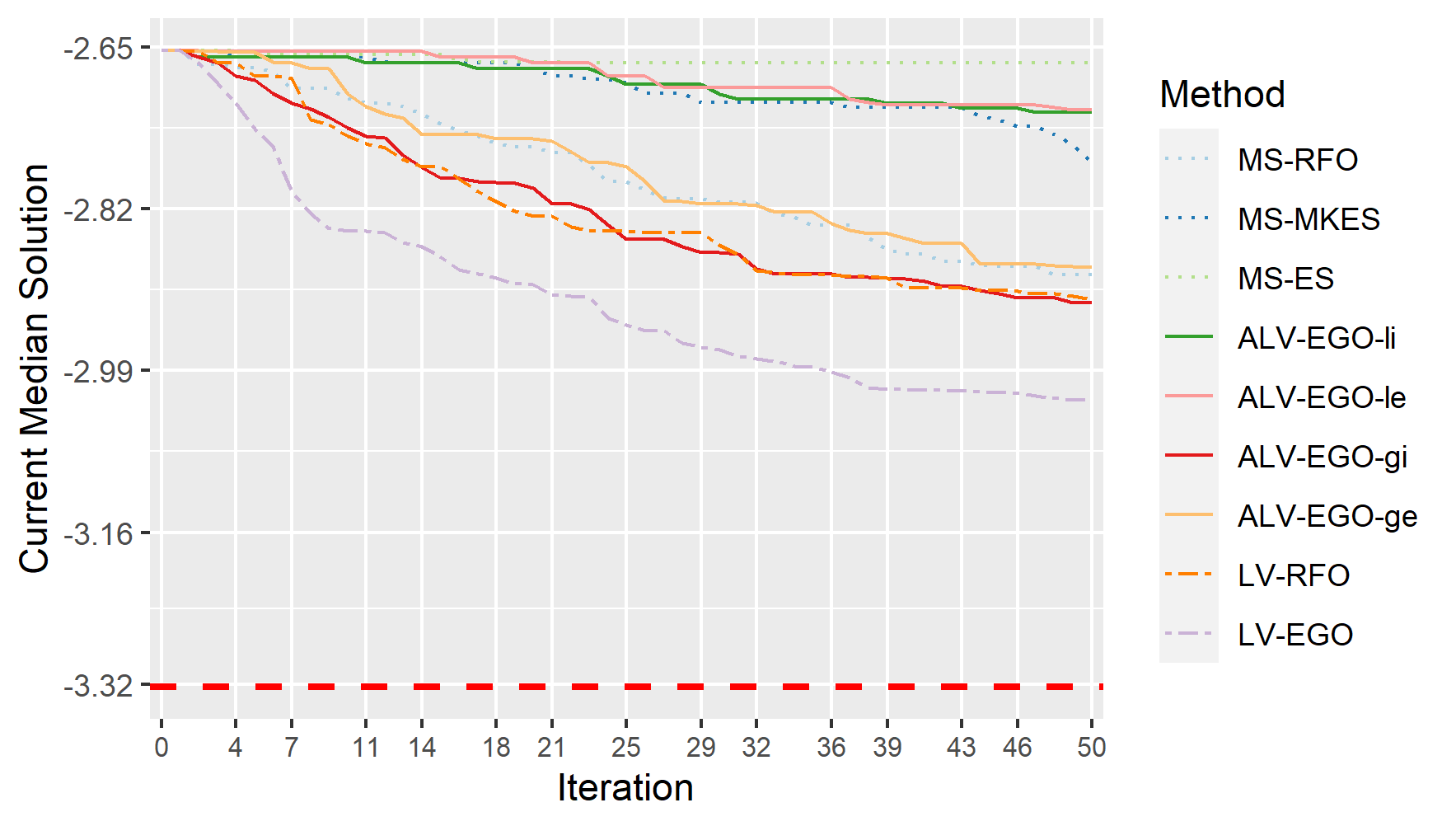}
			\caption{Median solution} 
			\label{figUh0:ua} 
			\vspace{4ex}
		\end{subfigure}
		\begin{subfigure}[b]{0.5\linewidth}
			\centering
			\includegraphics[width=\textwidth]{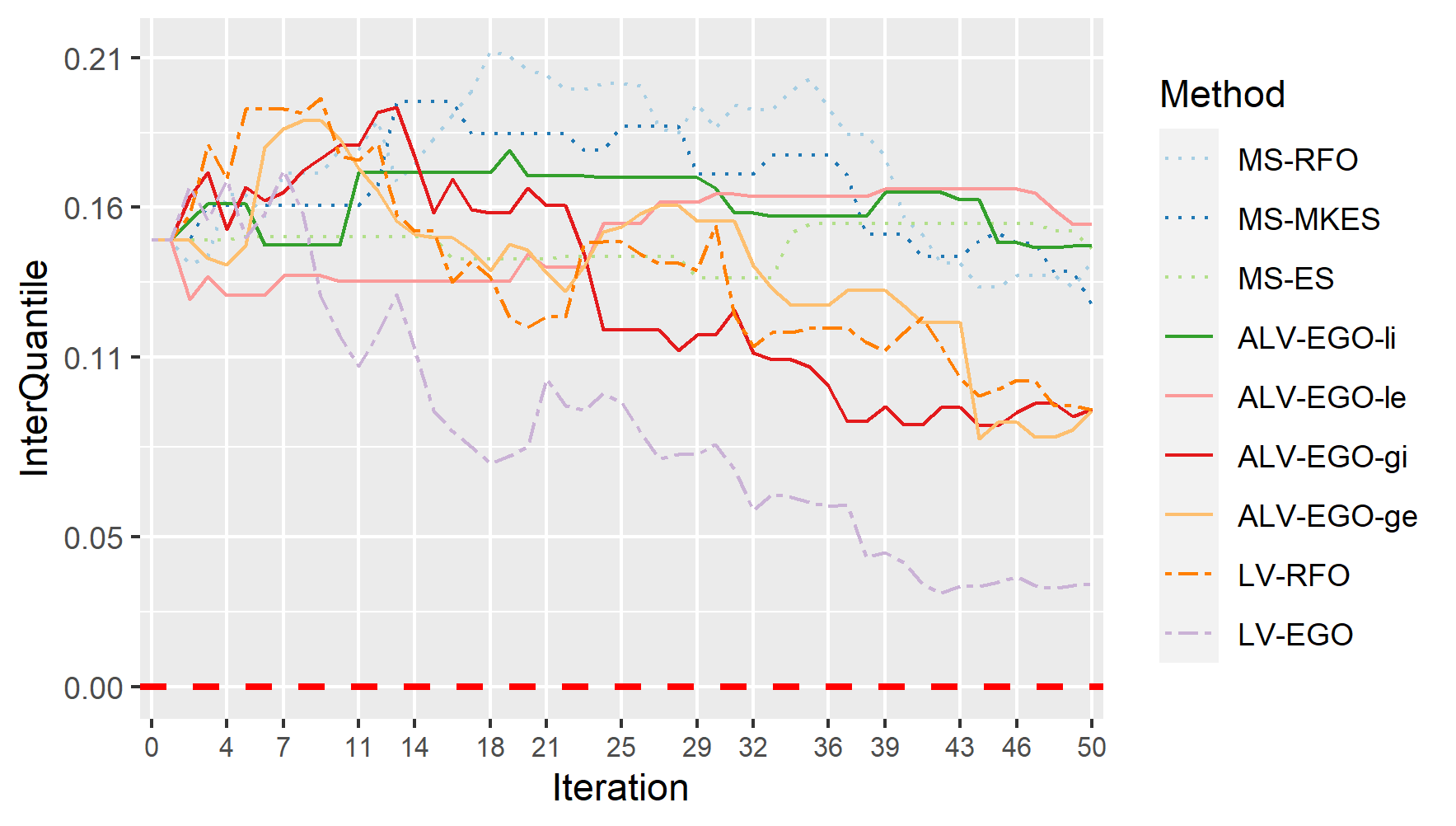}
			\caption{Interquartile range} 
			\label{figUh0:ub} 
			\vspace{4ex}
		\end{subfigure} 
		\begin{subfigure}[b]{0.5\linewidth}
			\centering
			\includegraphics[width=\textwidth]{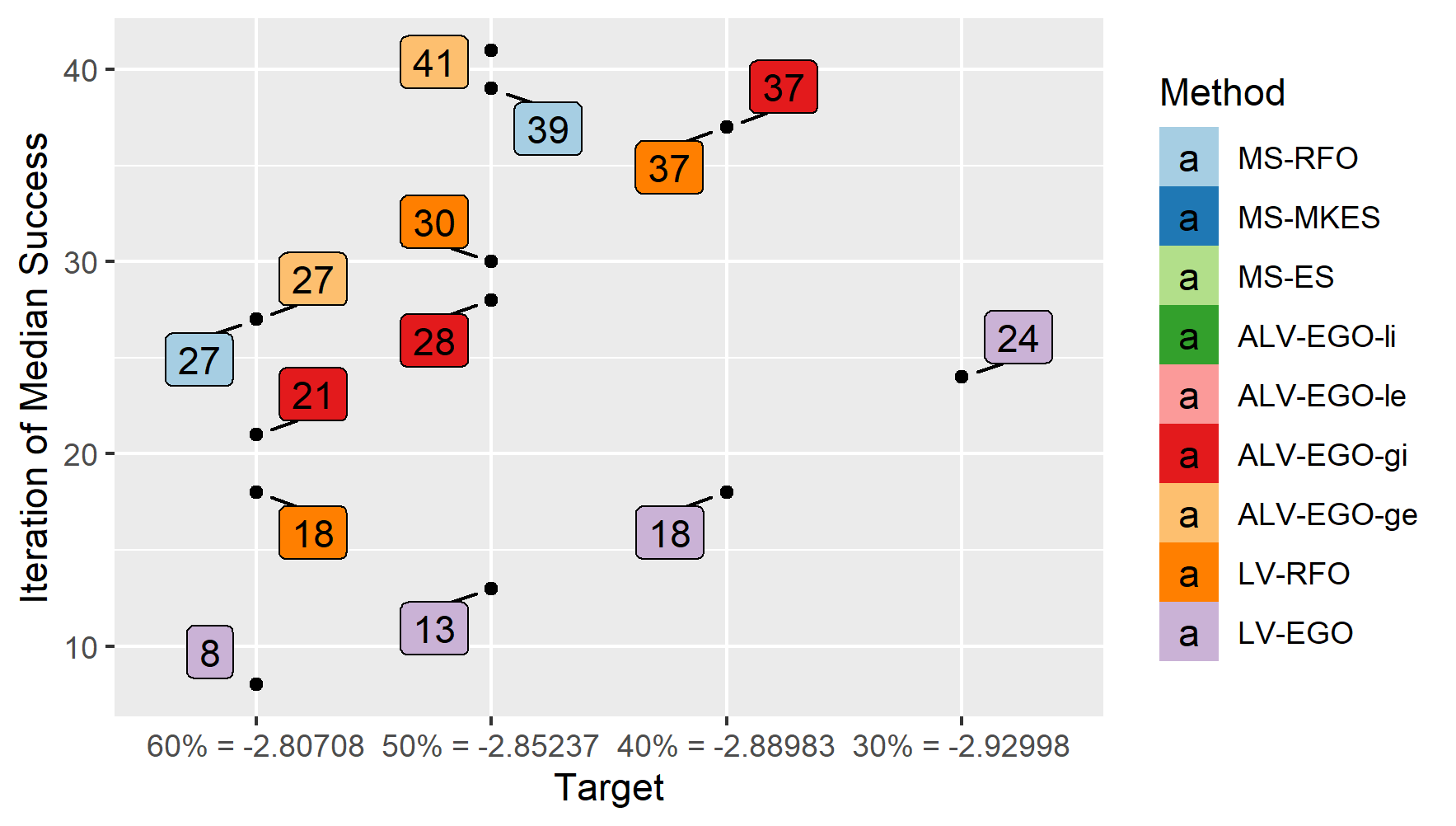}
			\caption{Iteration to median success} 
			\label{figUh0:uc} 
		\end{subfigure}
		\begin{subfigure}[b]{0.5\linewidth}
			\centering
			\includegraphics[width=\textwidth]{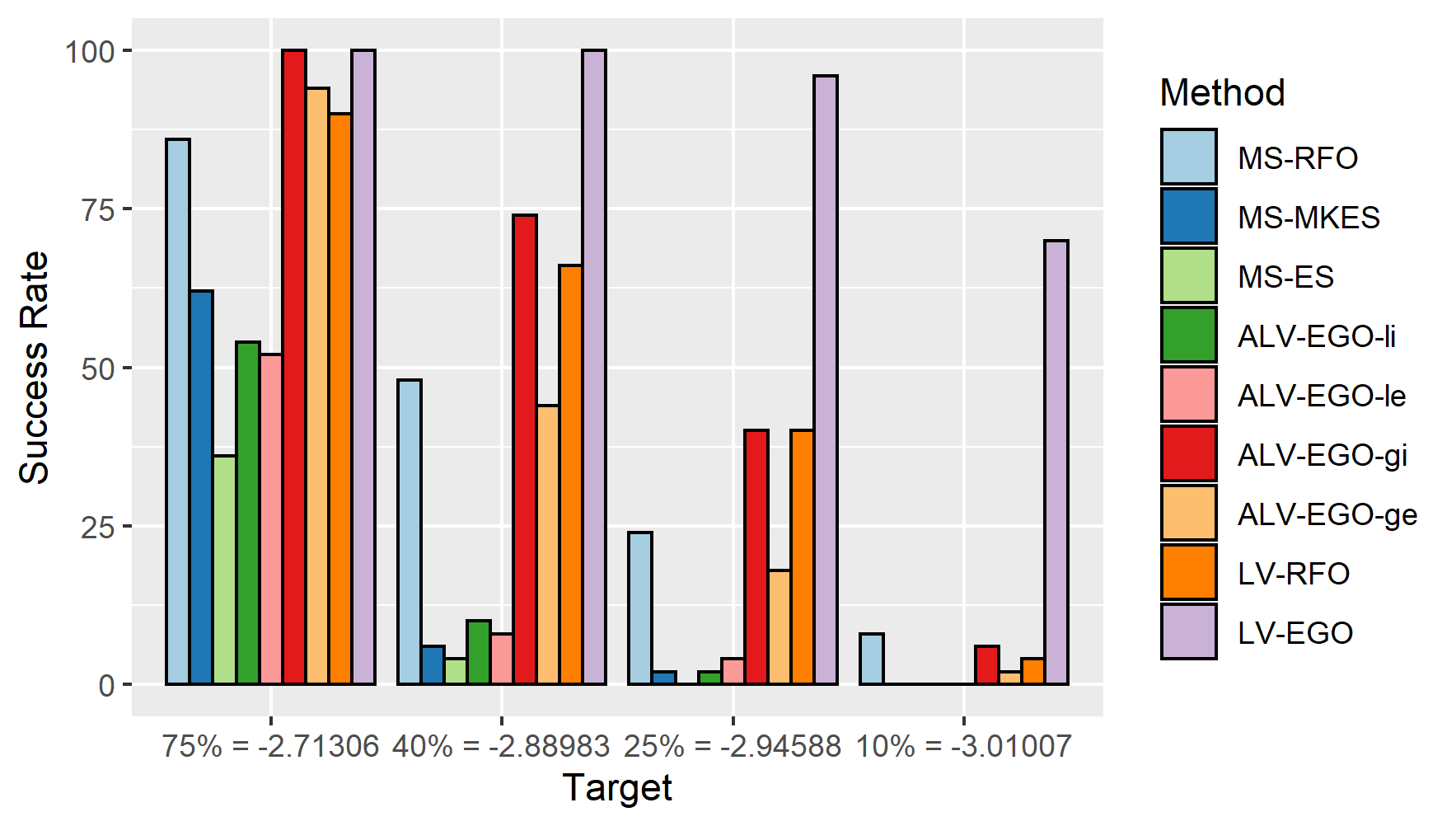}
			\caption{Success rate} 
			\label{figUh0:ud} 
		\end{subfigure} 
		\caption{Comparison of the 9 algorithms on the Hartmann function (for which $y^{\star} = -3.32237$). 
			\label{figUh0}}
	\end{figure}

	\paragraph{Hartmann function.}
	Results on the Hartmann function which has 4 continuous and 2 discrete variables, with a total of 9 discrete levels, will be impacted by the sensitivity of the algorithms to an increase in dimension. 
	These results are reported in Figure~\ref{figUh0}.
	
	LV-EGO stands out as the best method with respect to all criteria for Hartmann. 
	The next two best methods are LV-RFO and ALV-EGO-gi, followed by MS-RFO and ALV-EGO-ge.
	This time, LV-RFO and MS-RFO, which both rely on random forests, belong to the efficient methods: random forests gain in relative performance with respect to the GPs when the dimension and the size of the initial DoE increase. 
	For Hartmann, LV-EGO consistently outperforms the ALV- implementations.
	The importance of keeping the coupling between discrete and latent variables during the optimization seems less crucial, and even somewhat detrimental, in the Hartmann case. We think that this is due to the very tight budget (50 iterations after the initial DoE) which does not allow the convergence of the optimizers, as can be seen in the Plot~\ref{figUh0:ua} where the global optimum is not reached. 
	Because the optimum is not really found, constraints on discreteness are superfluous and their handling through the pre-image problem is sufficient.
	As in the other test cases, MS-ES was slower than the other methods.
		
	\subsection{Beam bending application}
	\label{sec:latdisc}
	
	\paragraph{Optimization results.}
	
\begin{figure}[H] 
		\begin{subfigure}[b]{0.5\linewidth}
			\centering
			\includegraphics[width=\textwidth]{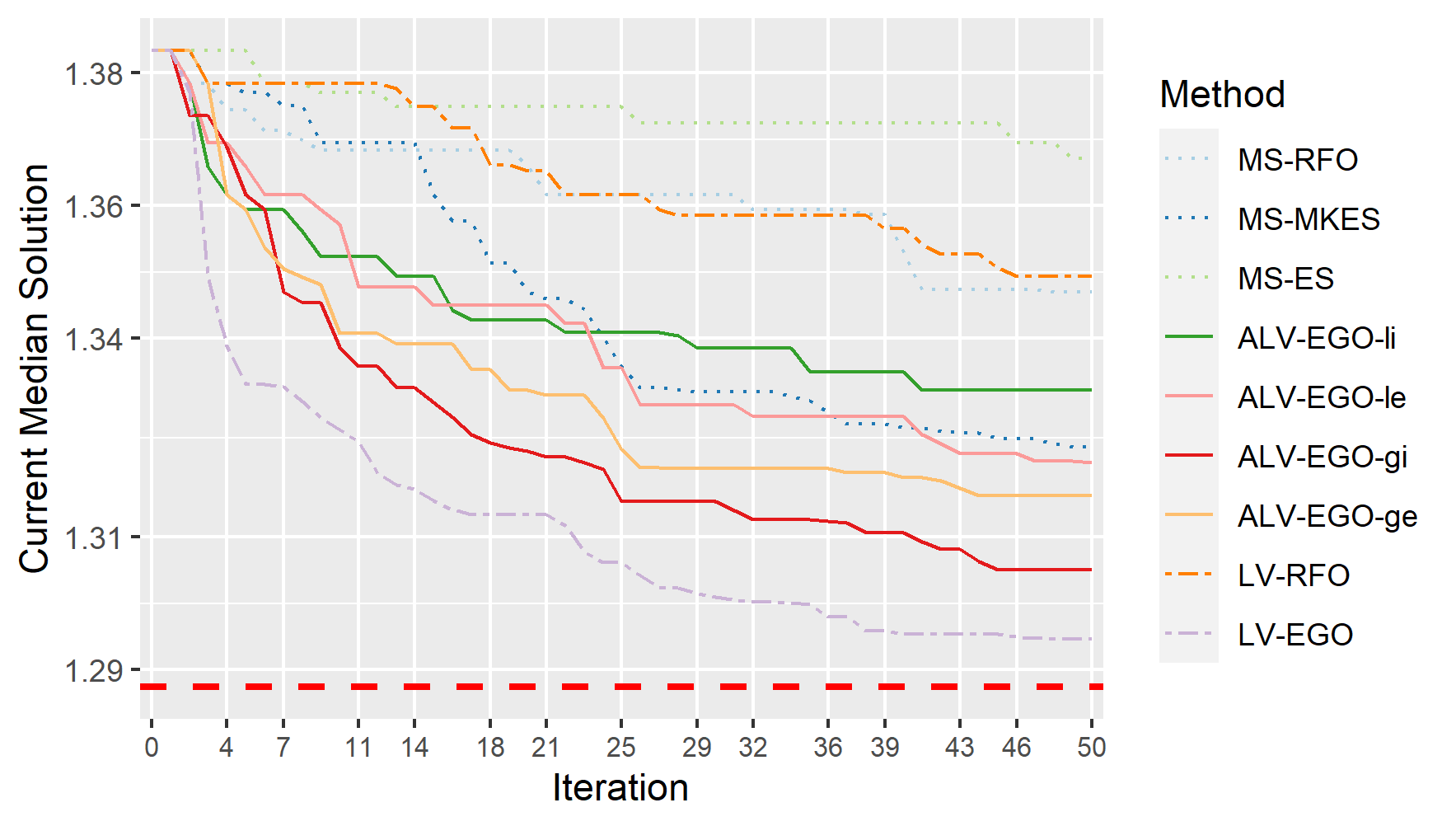}
			\caption{Median solution} 
			\label{figUb0:ua} 
			\vspace{4ex}
		\end{subfigure}
		\begin{subfigure}[b]{0.5\linewidth}
			\centering
			\includegraphics[width=\textwidth]{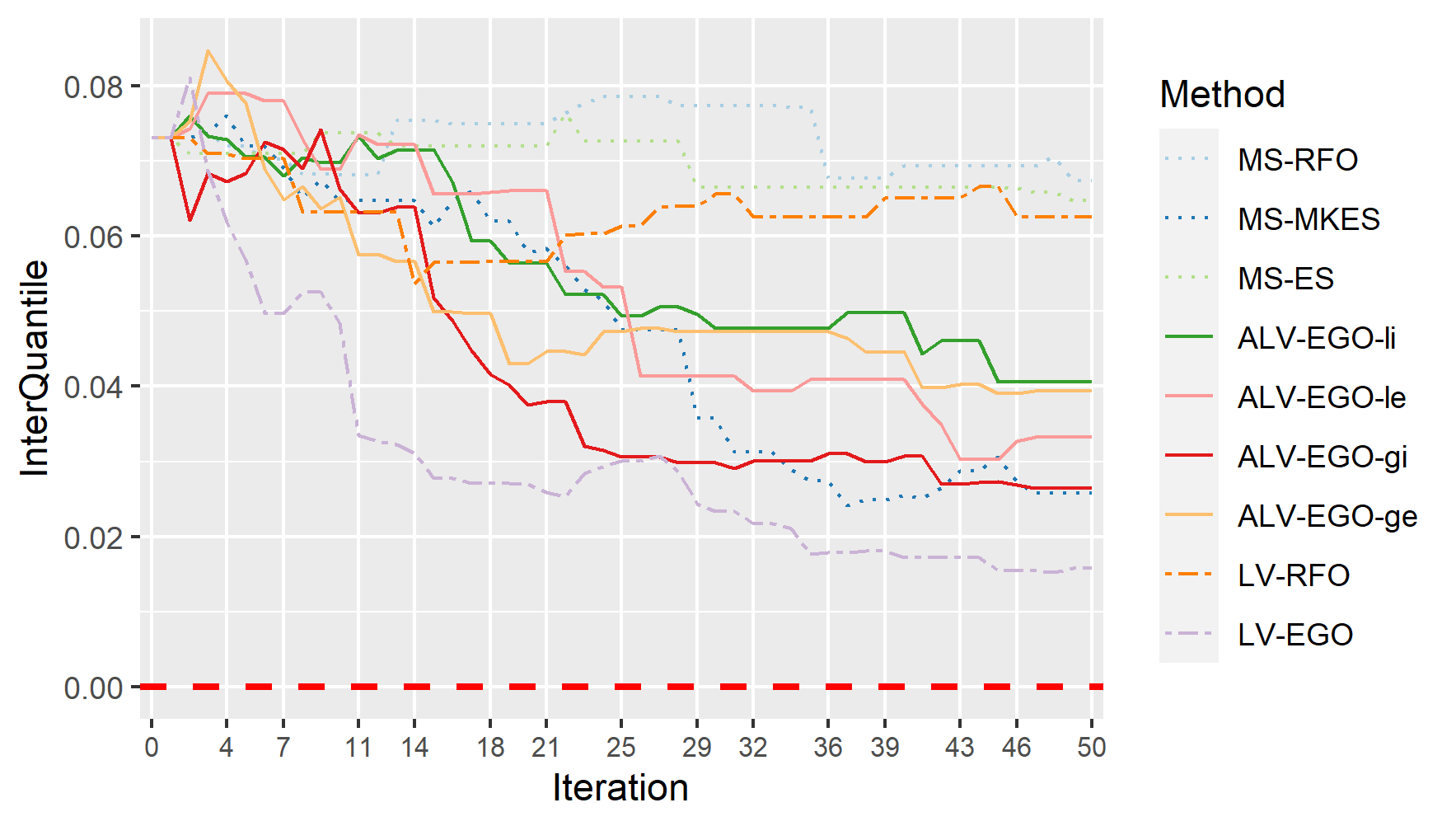}
			\caption{Interquartile range} 
			\label{figUb0:ub} 
			\vspace{4ex}
		\end{subfigure} 
		\begin{subfigure}[b]{0.5\linewidth}
			\centering
			\includegraphics[width=\textwidth]{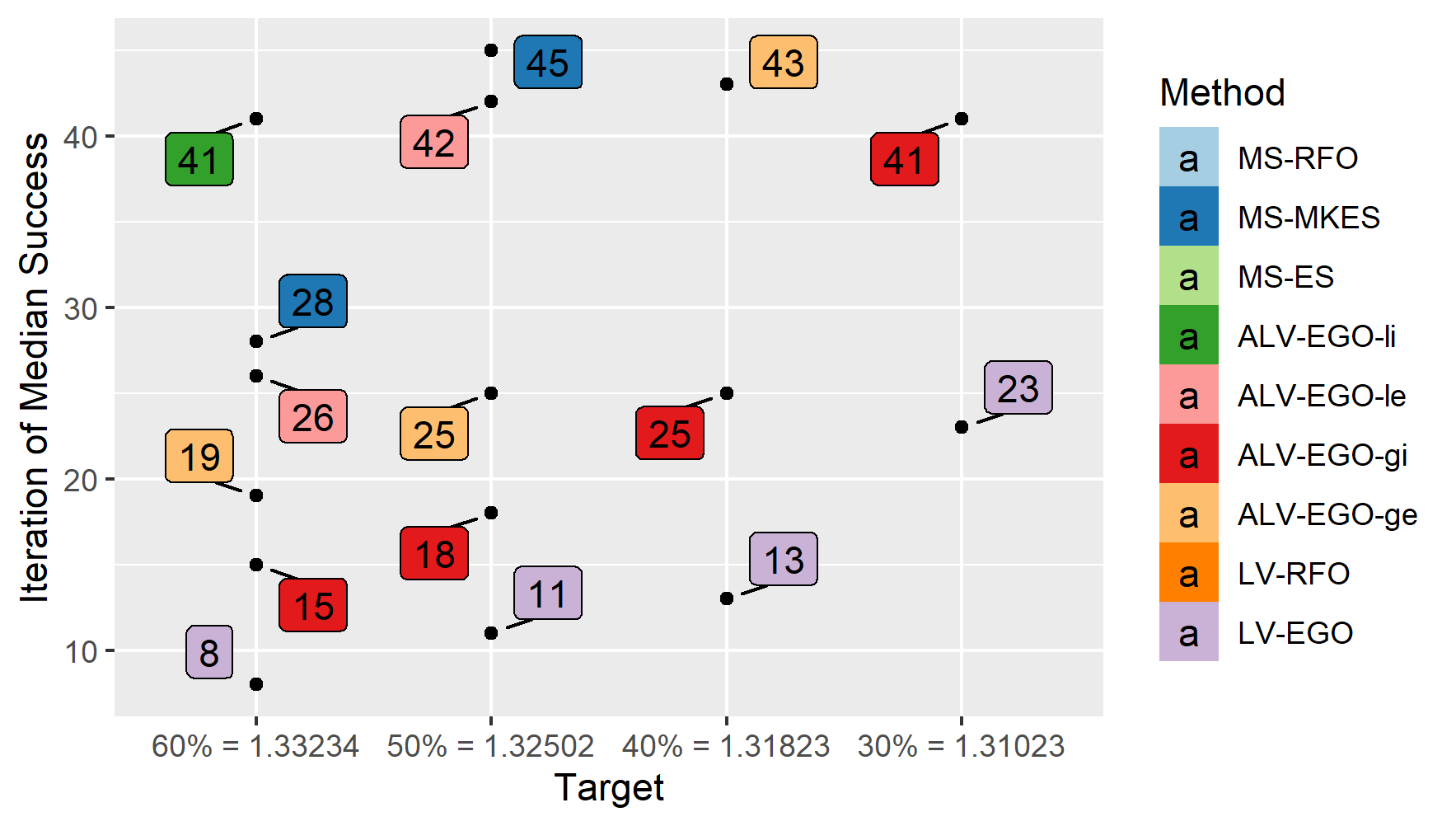}
			\caption{Iteration to median success} 
			\label{figUb0:uc} 
		\end{subfigure}
		\begin{subfigure}[b]{0.5\linewidth}
			\centering
			\includegraphics[width=\textwidth]{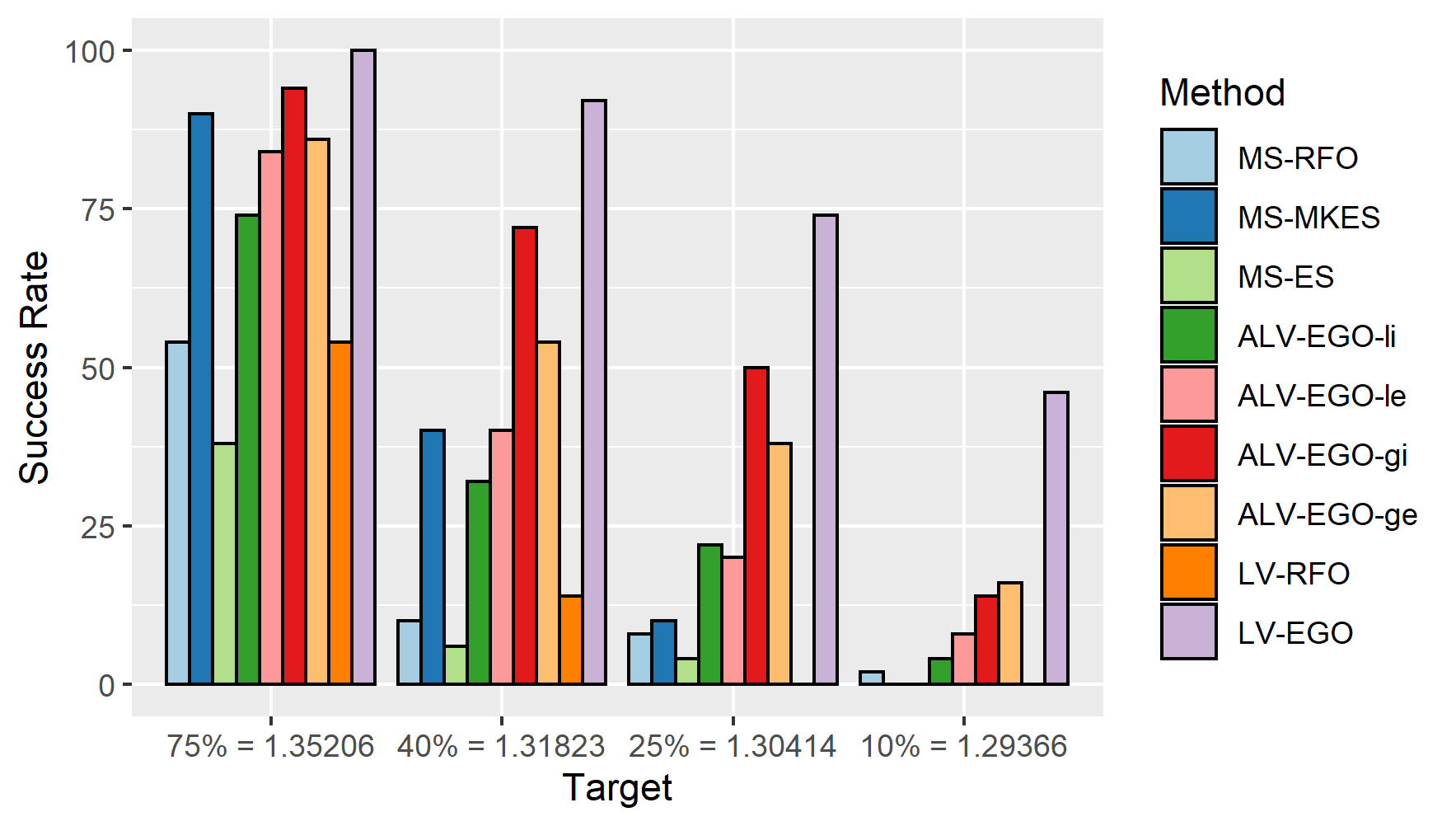}
			\caption{Success rate} 
			\label{figUb0:ud} 
		\end{subfigure} 
		\caption{Comparison of all 9 algorithms on the beam design test case ($y^{\star} = 1.28738$). }
		
		\label{figUb0} 
	\end{figure}

	Figure~\ref{figUb0} summarizes the 4 comparison metrics of all 9 algorithms in the bended beam test case. 
	The ranking of the algorithms is similar to that obtained with the Branin and Goldstein functions.
	LV-EGO has the best convergence both in terms of median speed (cf. plots of the left column) and accuracy (bottom right plot).  
	ALV-EGO-gi is the second most efficient method followed by ALV-EGO-ge. 
	Again, the algorithms that resort to random forests, LV-RFO and MS-RFO, are the slowest and most inaccurate. 
	They share this counter-performance with MS-ES.
	
\begin{figure}[h!] 
		\stackunder[5pt]{\includegraphics[width=1.cm,height=1.cm]{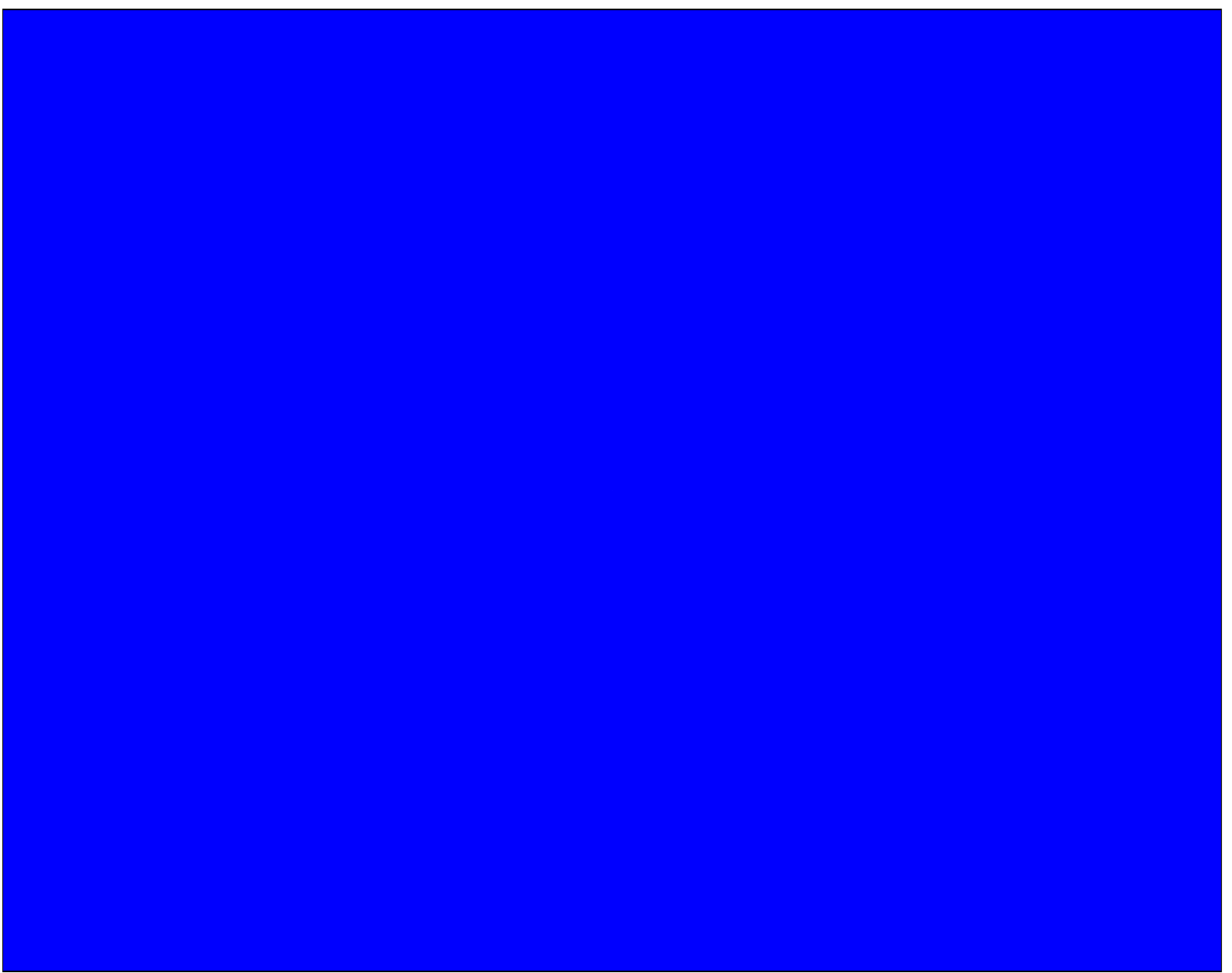}}{$\tilde{I}_1$}
		\hspace{0.2cm}
		\stackunder[5pt]{\includegraphics[width=1.cm,height=1.cm]{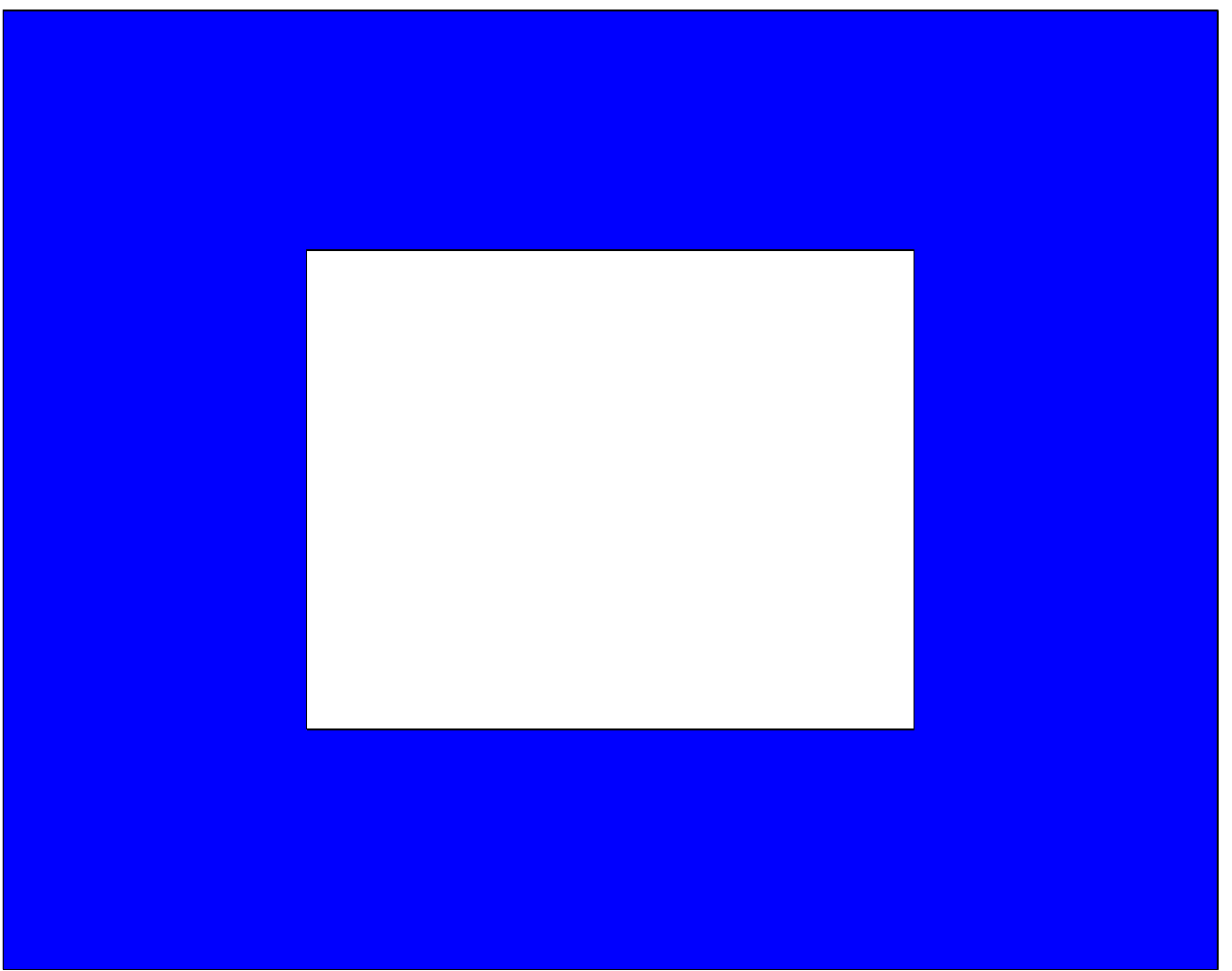}}{$\tilde{I}_2$}
		\hspace{0.2cm}
		\stackunder[5pt]{\includegraphics[width=1.cm,height=1.cm]{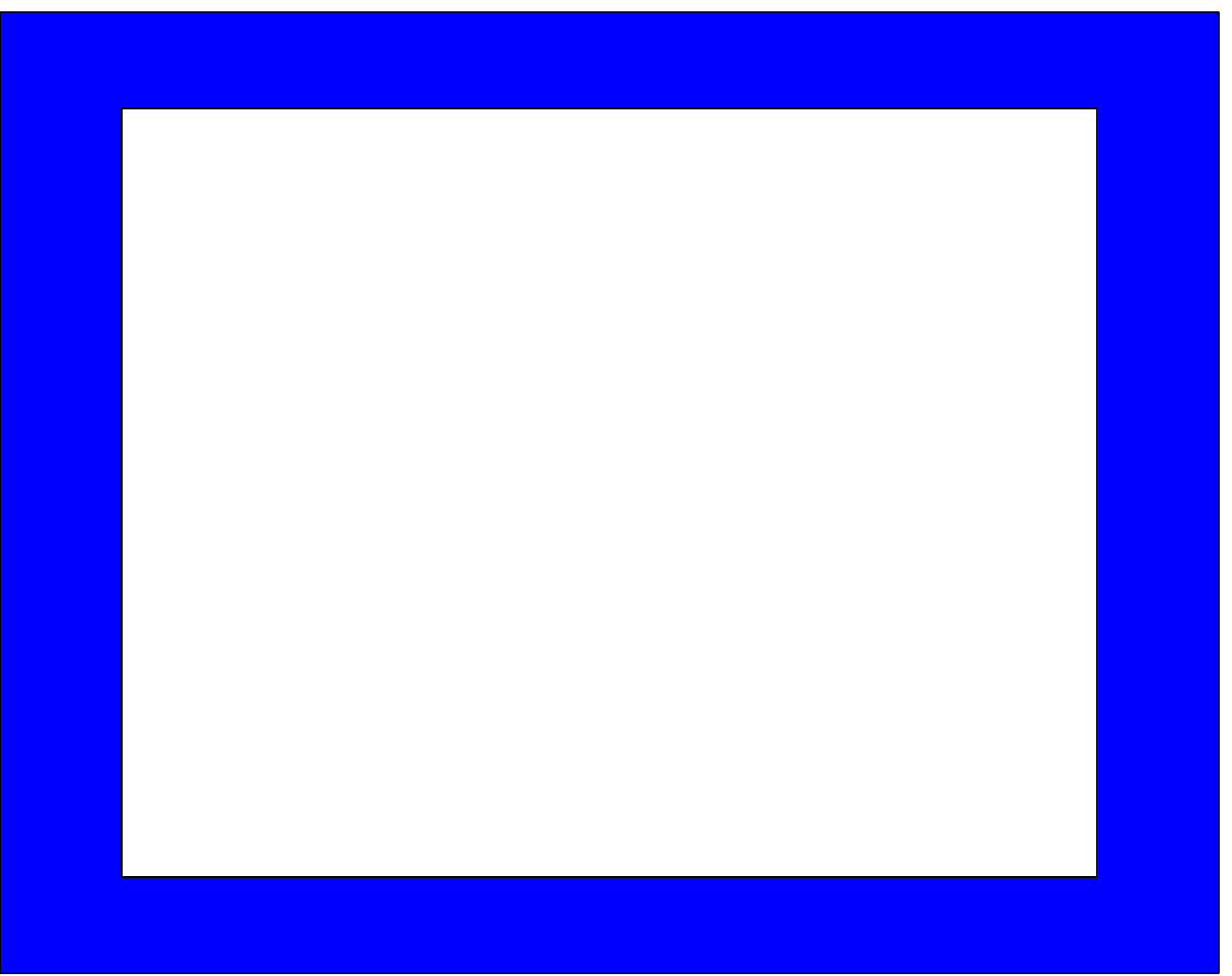}}{$\tilde{I}_3$}
		\hspace{0.2cm}
		\stackunder[5pt]{\includegraphics[width=1.cm,height=1.cm]{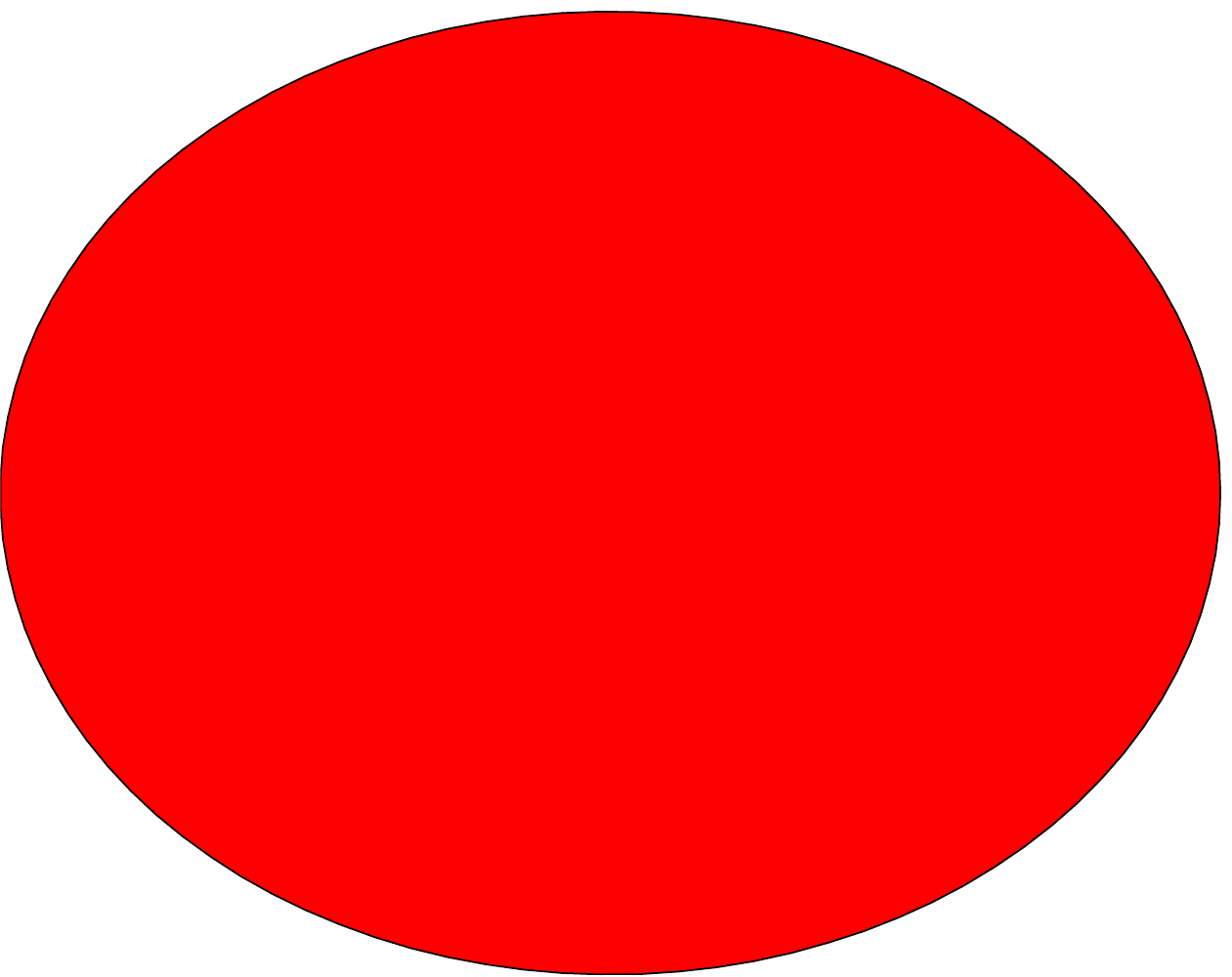}}{$\tilde{I}_4$}
		\hspace{0.2cm}
		\stackunder[5pt]{\includegraphics[width=1.cm,height=1.cm]{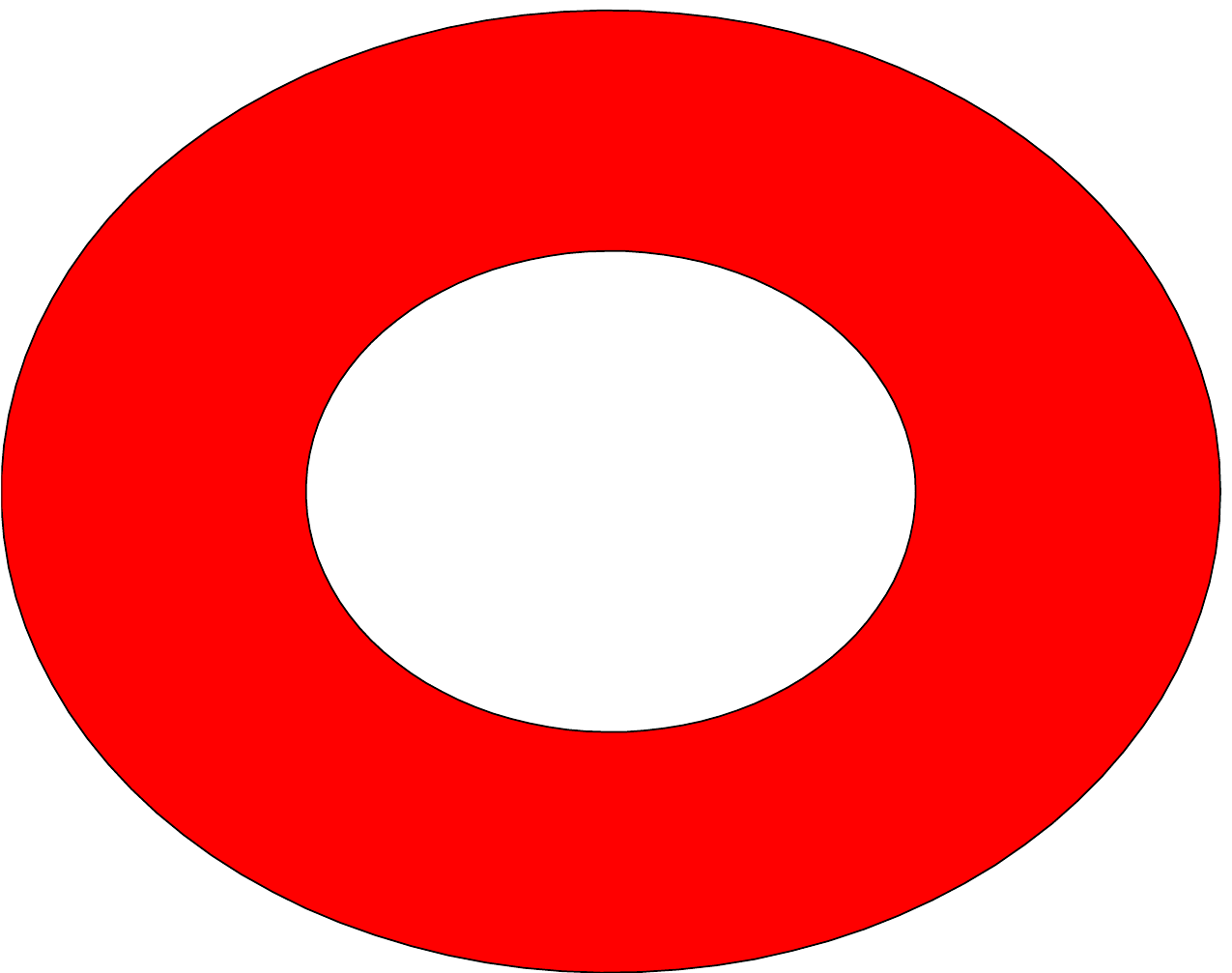}}{$\tilde{I}_5$}
		\hspace{0.2cm}
		\stackunder[5pt]{\includegraphics[width=1.cm,height=1.cm]{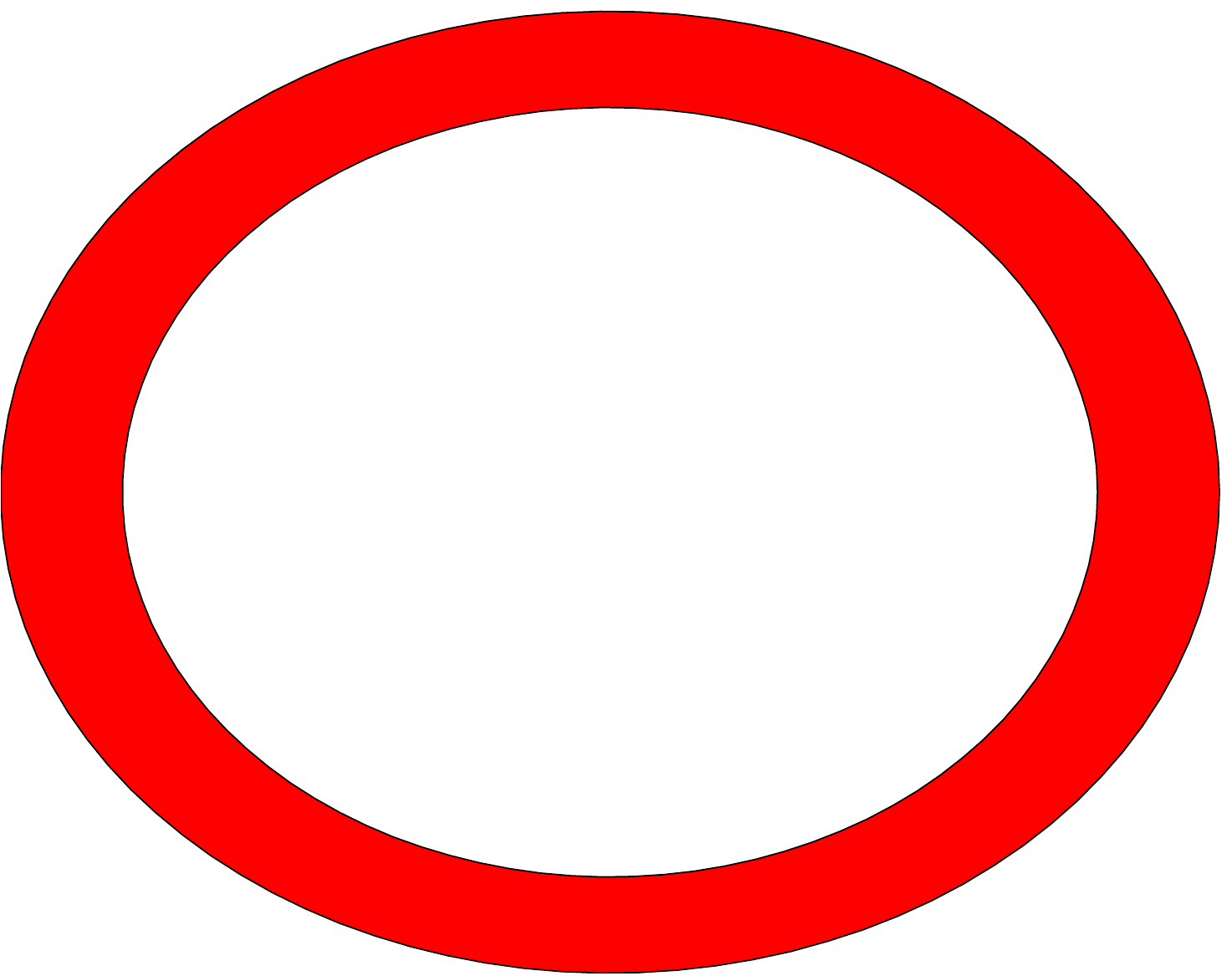}}{$\tilde{I}_6$}
		\hspace{0.2cm}
		\stackunder[5pt]{\includegraphics[width=1.cm,height=1.cm]{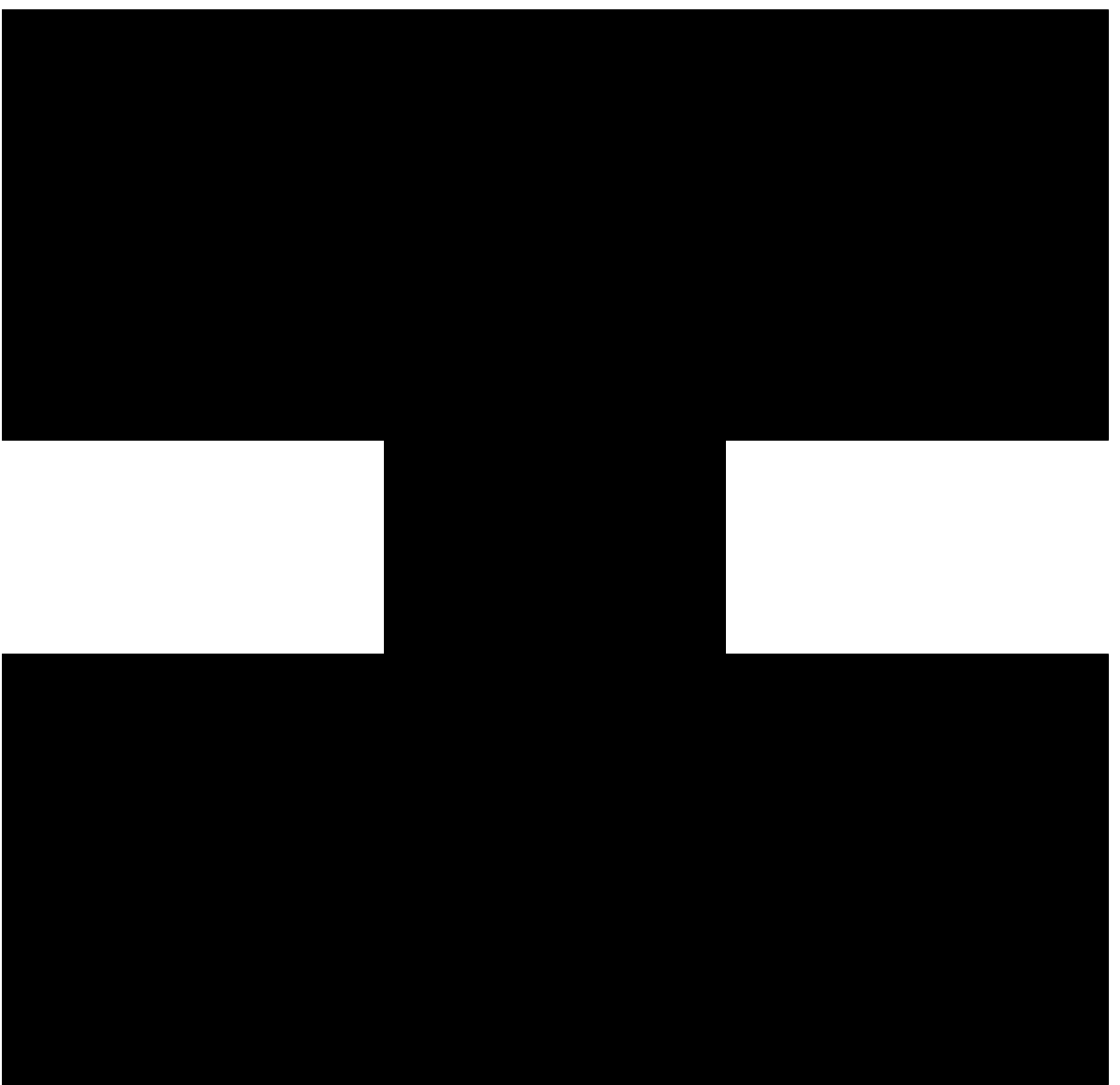}}{$\tilde{I}_7$}
		\hspace{0.2cm}
		\stackunder[5pt]{\includegraphics[width=1.cm,height=1.cm]{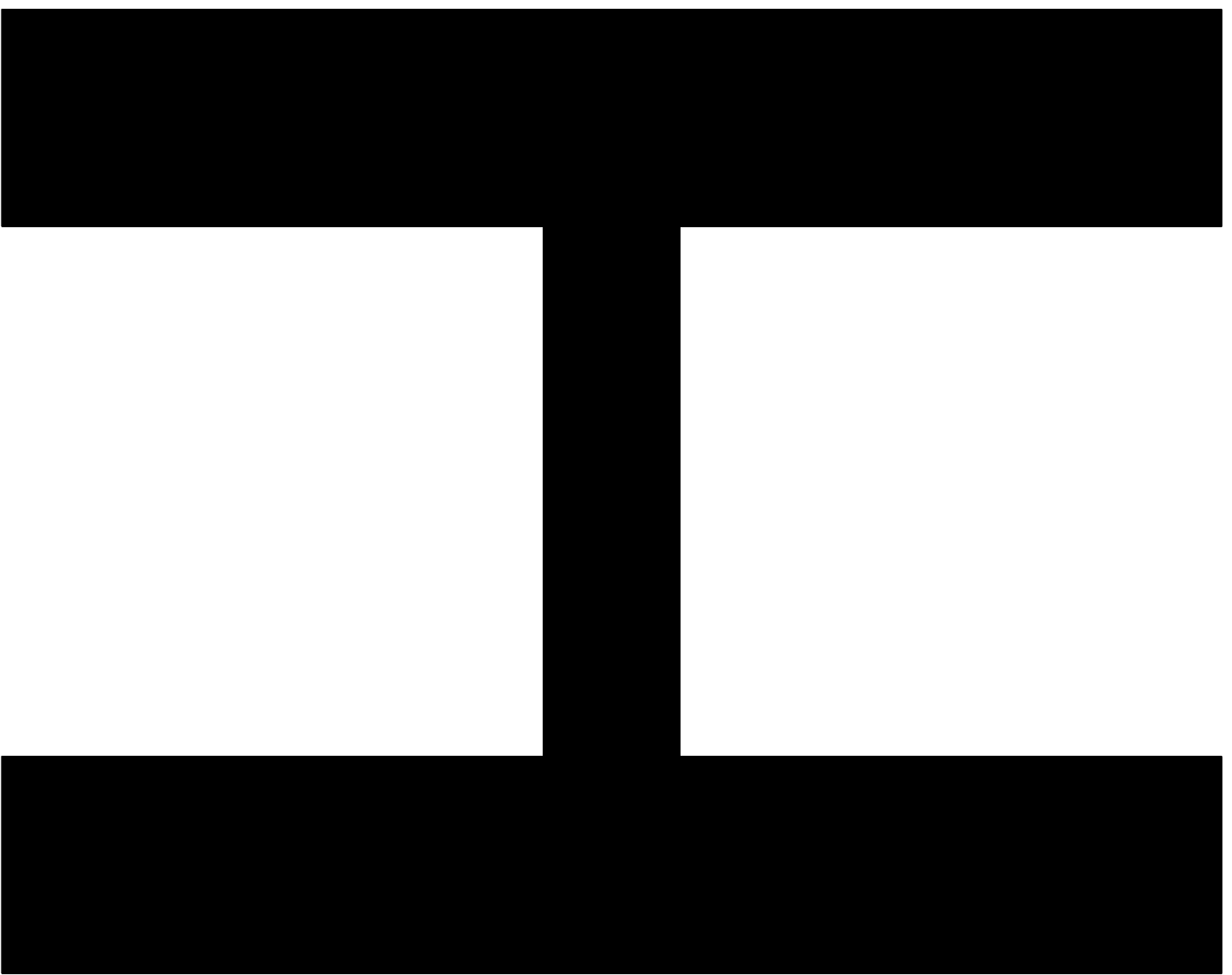}}{$\tilde{I}_8$}
		\hspace{0.2cm}
		\stackunder[5pt]{\includegraphics[width=1.cm,height=1.cm]{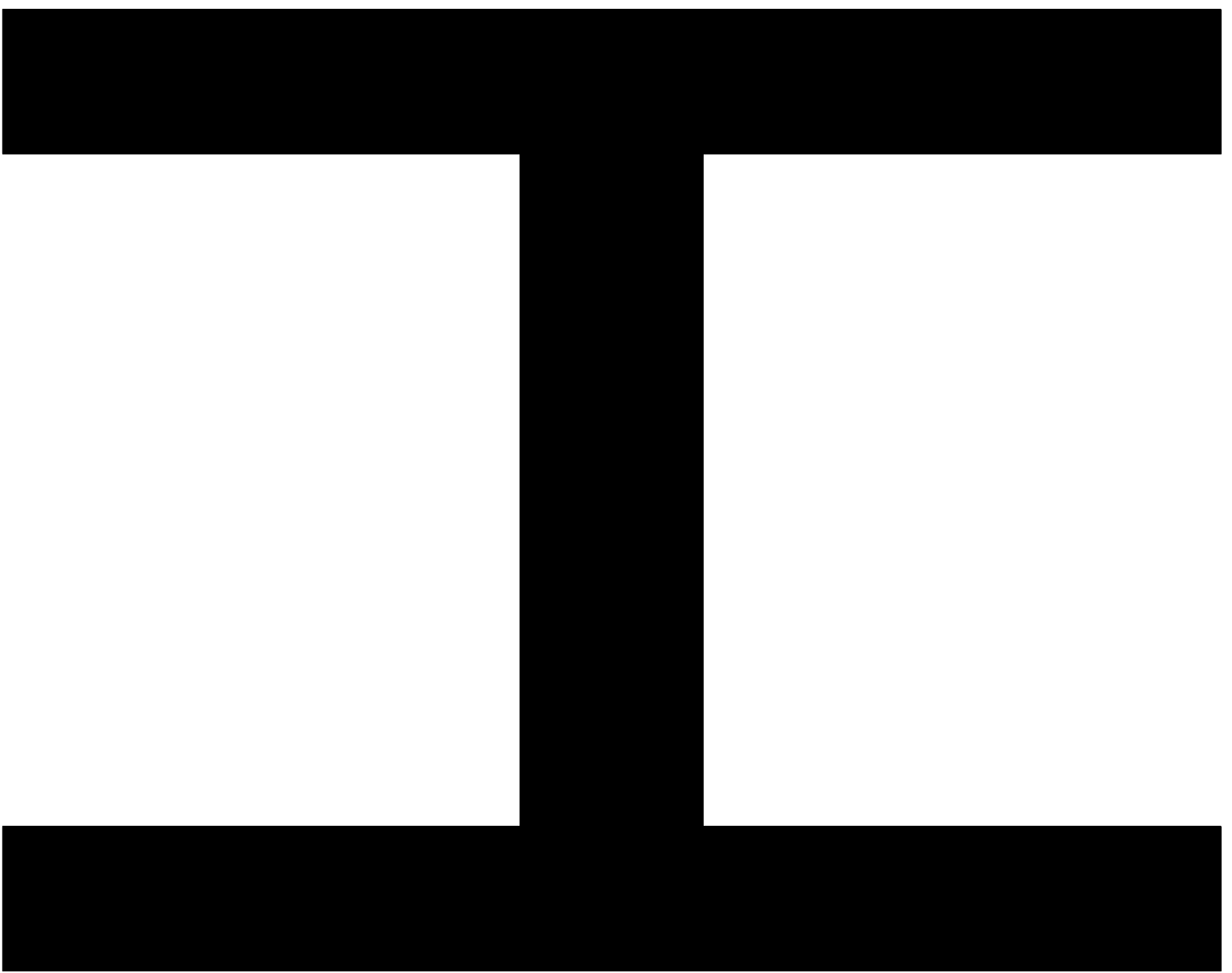}}{$\tilde{I}_9$}
		\hspace{0.2cm}
		\stackunder[5pt]{\includegraphics[width=1.cm,height=1.cm]{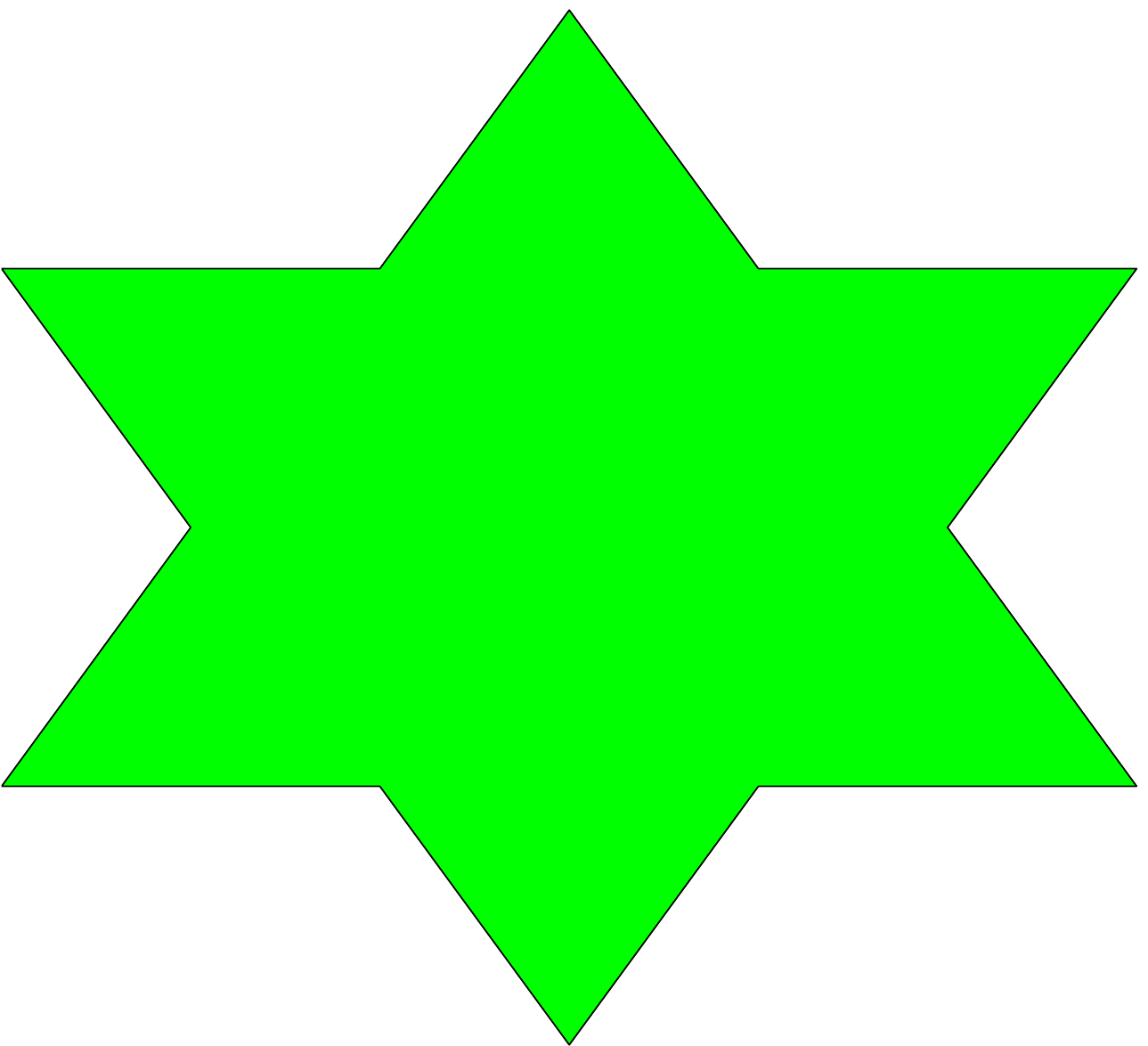}}{$\tilde{I}_{10}$}
		\hspace{0.2cm}
		\stackunder[5pt]{\includegraphics[width=1.cm,height=1.cm]{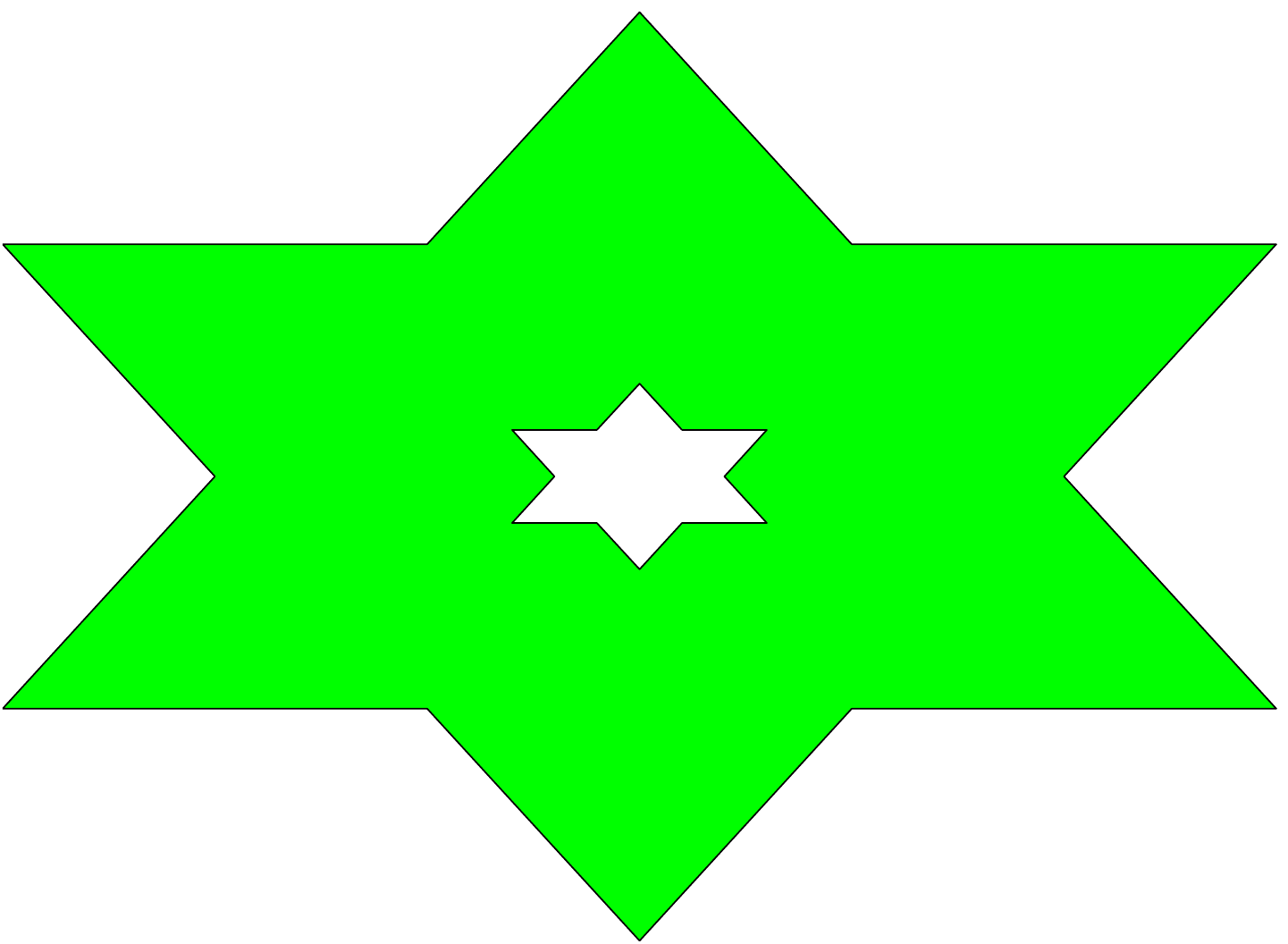}}{$\tilde{I}_{11}$}
		\hspace{0.2cm}
		\stackunder[5pt]{\includegraphics[width=1.cm,height=1.cm]{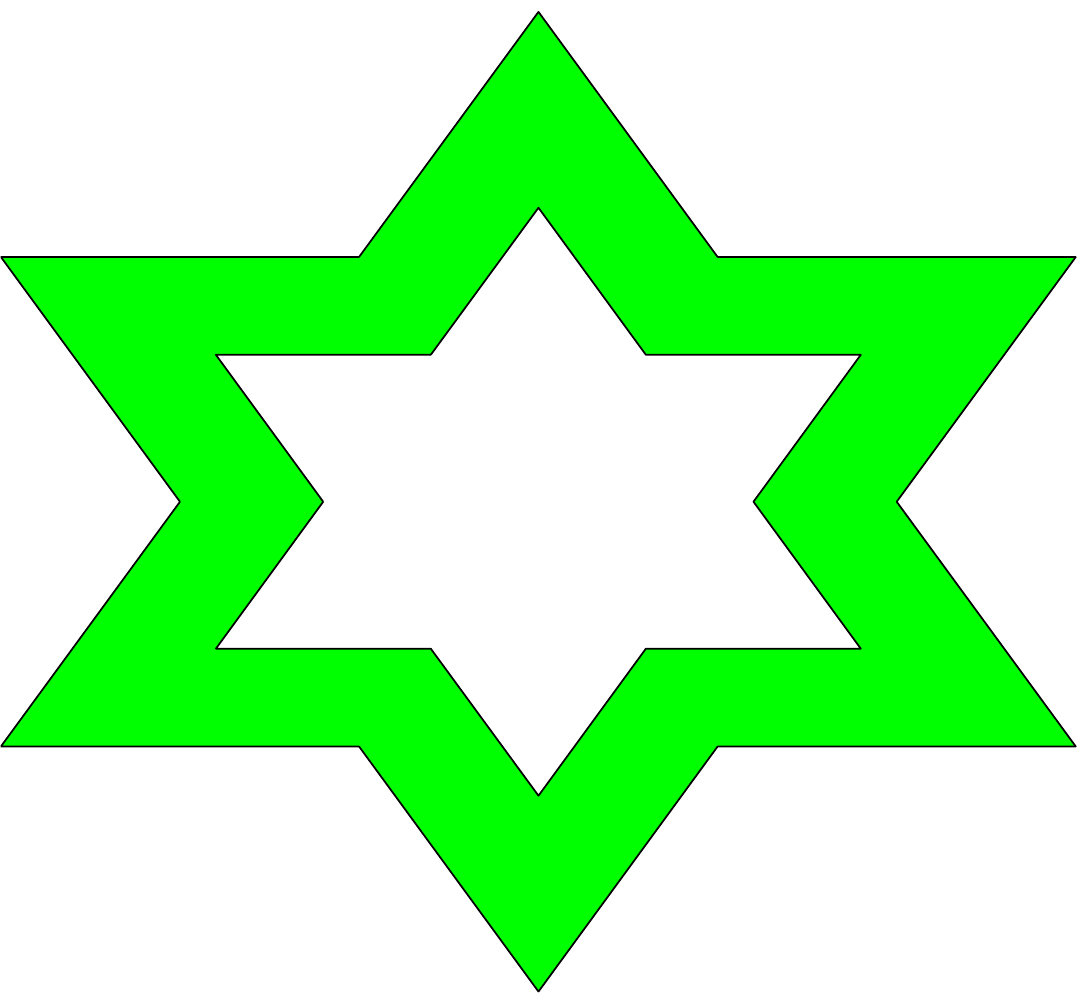}}{$\tilde{I}_{12}$}
		\caption{Shapes of the considered beam profiles. The scale differs from one picture to another, as the areas are supposed to be the same for each cross-section. From \cite{roustant_group_2020}.
			\label{fig:beamSections}
		}
	\end{figure}

\begin{figure}[H] 
		\begin{subfigure}[b]{0.5\linewidth}
			\centering
			\includegraphics[width=\textwidth]{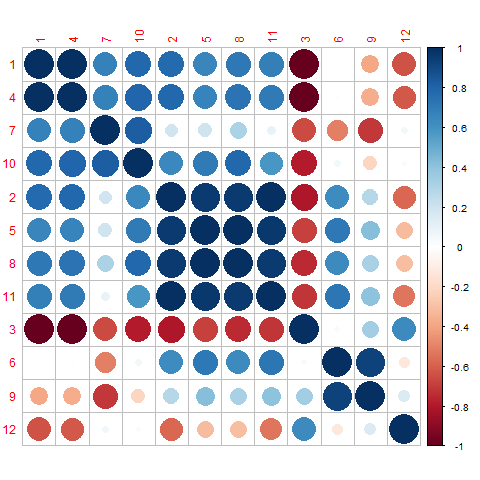}			
			
			\caption{Correlation of the latent variables at iteration $\# 1$} 
			\label{figlv0:ua} 
			\vspace{4ex}
		\end{subfigure}
		\begin{subfigure}[b]{0.5\linewidth}
			\centering
			\includegraphics[width=\textwidth]{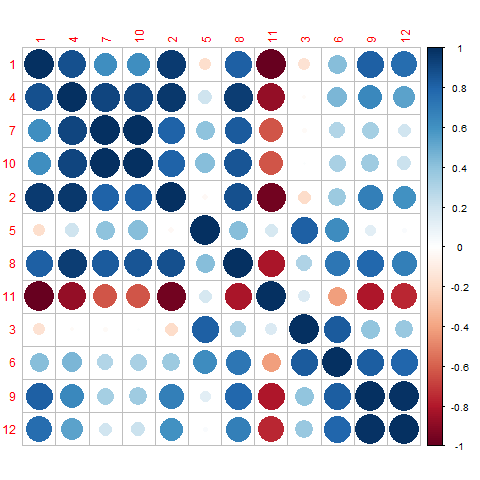}
			\caption{Correlation of the latent variables at iteration $\# 26$} 
			\label{figlv0:ub} 
			\vspace{4ex}
		\end{subfigure} 
		\begin{subfigure}[b]{0.5\linewidth}
			\centering
			\includegraphics[width=\textwidth]{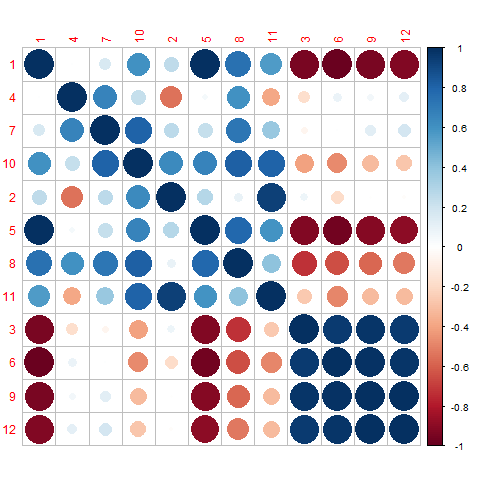}
			\caption{Correlation of the latent variables at iteration $\# 49$} 
			\label{figlv0:uc} 
		\end{subfigure}
		\begin{subfigure}[b]{0.5\linewidth}
			\centering
			\includegraphics[width=\textwidth]{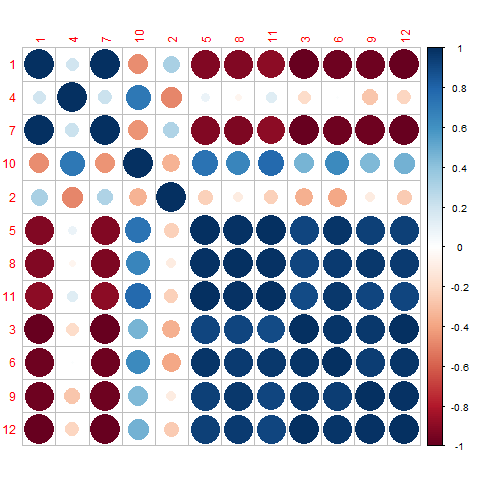}
			\caption{Correlation of the latent variables at iteration $\# 50$} 
			\label{figlv0:ud} 
		\end{subfigure} 
		\caption{
Representation of the correlation between the latent variables at various iterations $t$. 
{The size and the color of the circles correspond to the absolute and signed magnitudes of the correlations.}
The correlations are those of the categorical kernel of Equation (\ref{eq:kerCatLV}). 
The levels were grouped according to $\tilde{I}$:  $\{1, 4, 7, 10\}$, $\{2, 5, 8, 11\}$, $\{3, 6, 9, 12\}$.
			\label{figlv0} 
		}
	\end{figure}

	\paragraph{Latent variables in the beam application.}
	The beam subject to a bending load is a test case that allows to interprete the latent variables.
	Indeed, the normalized moment of inertia, $\tilde I$, is a candidate latent variable once it is allowed to take continuous values as it determines, 
	with the continuous cross-section $S$ and the length $L$, the output (the penalized beam deflection) $y$ in Equation~(\ref{eq:bendingVol}).
	The levels of $\tilde I$ (given in Equation~(\ref{eq:Iset})) correspond to 3 increasingly hollow profiles of 4 shapes, as illustrated in Figure~\ref{fig:beamSections}. 
	Because a relaxed $\tilde I$ is a possible latent variable, it is expected that the latent variables $\Latvec{t}$ learned from the data will be grouped in the same way as $\tilde I$. 
	Looking at $\tilde I$ values and at Figure \ref{fig:beamSections}, we thus expect, in the image space defined by latent variables, three groups of levels: those corresponding to solid forms (levels $\{1, 4, 7, 10\}$), medium-hollow forms (levels $\{2, 5, 8, 11\}$) and hollow forms (levels $\{3, 6, 9, 12\}$).
	For the sake of interpretation, we select $1$ run that found the global optimum with the Vanilla LV-EGO algorithm. In Figure \ref{figlv0}, we represent in a color scale the estimated correlation matrix corresponding to the categorical kernel of Equation (\ref{eq:kerCatLV}), at iterations $[1;26;49;50]$.
	At the beginning of the optimization, at iteration 1, we can see a block-structure which corresponds quite well to the three groups of forms described above. This structure becomes less clear for the next iterations of the LV-EGO algorithm. This may be explained by the fact that the algorithm creates an unbalanced design, with more points in the promising areas according to the optimizers, so that all levels are no longer properly represented. 

	\subsection{Summary and discussion}
	\label{sec:Rsumm}
\begin{figure}[H] 
		\begin{subfigure}[b]{0.5\linewidth}
			\centering
			\includegraphics[width=\textwidth]{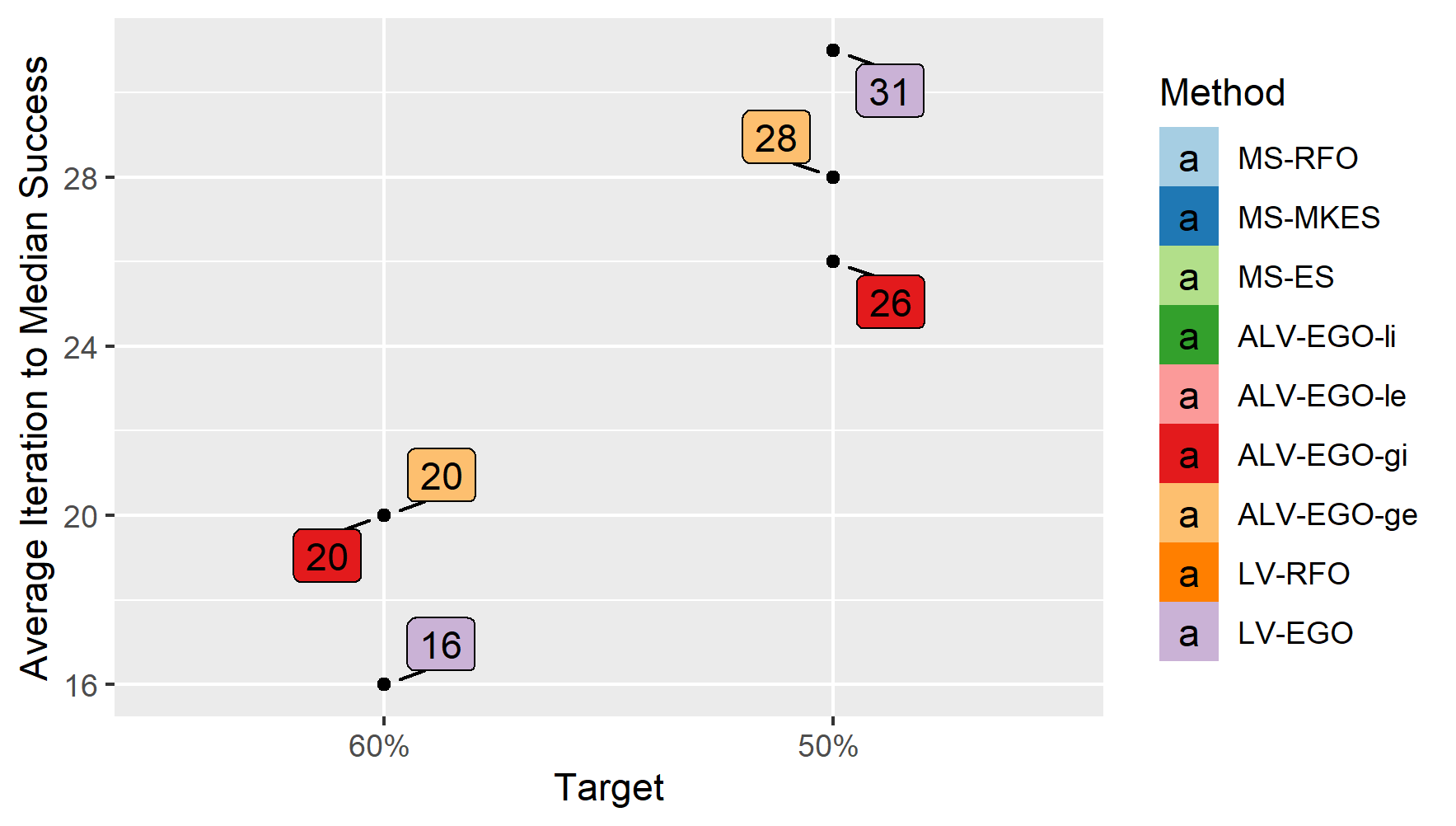}
			\caption{Average iteration to median success} 
			\label{figC:ua} 
			\vspace{4ex}
		\end{subfigure}
		\begin{subfigure}[b]{0.5\linewidth}
			\centering
			\includegraphics[width=\textwidth]{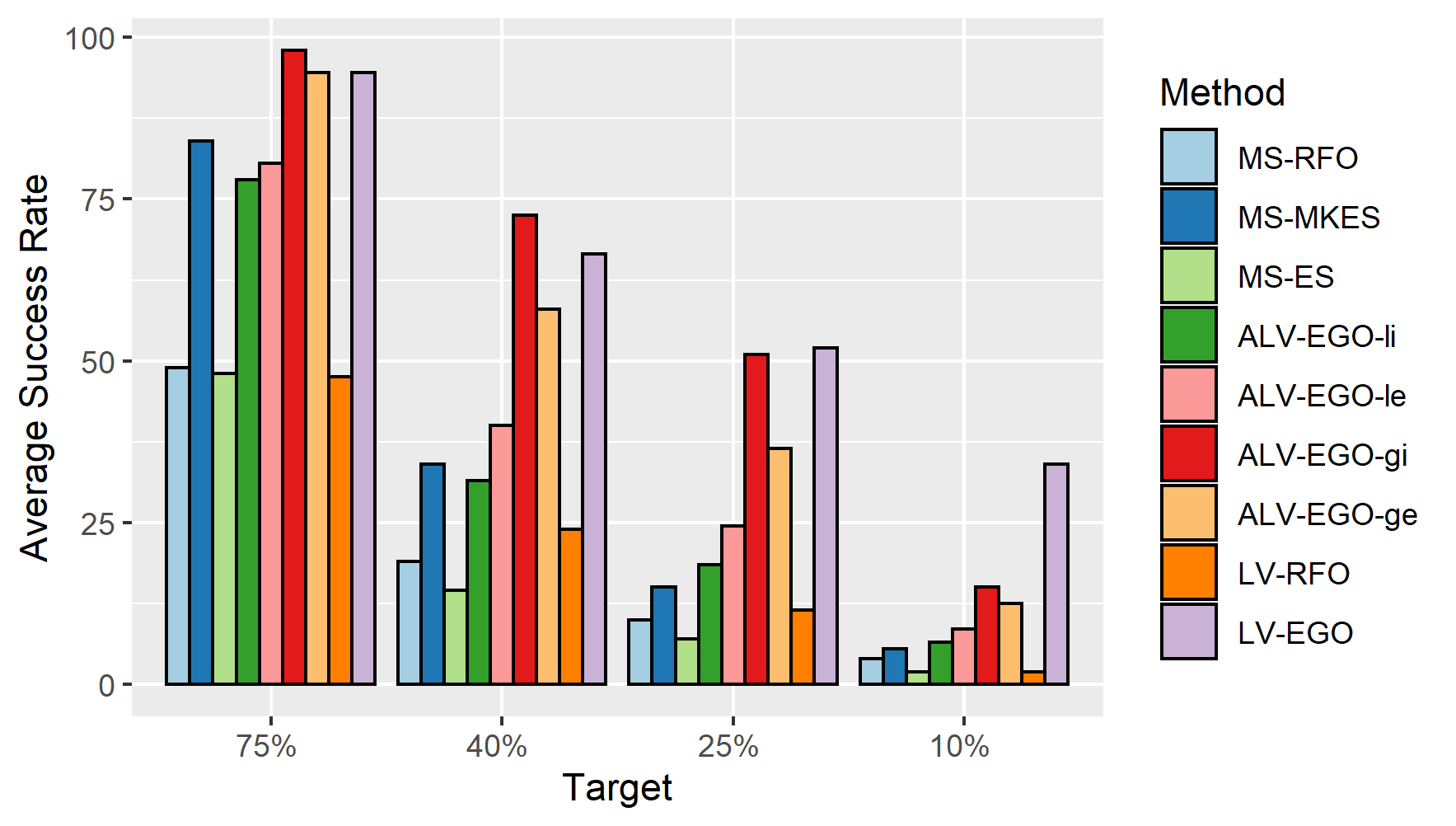}
			\caption{Average success rates} 
			\label{figC:ub} 
			\vspace{4ex}
		\end{subfigure} 
		\caption{Comparison of the 9 algorithms tested with results averaged over all test cases.}
		\label{figC} 
	\end{figure}

	The results of all the previous test cases which are measured through targets can be averaged.
	For example, the success rate of an algorithm at 25\% difficulty is the average of the rates for the 25\% quantiles of all test cases.
	The average results are presented in Figure~\ref{figC}.
	
	The three leading algorithms out of the 9 tested are ALV-EGO-gi, -ge and LV-EGO. 
	Among them, LV-EGO is slightly better at locating difficult targets (10\% quantile) while ALV-EGO-gi (closely followed by ALV-EGO-ge) is more robust at locating 50\% targets as can be seen from the median success plot in Figure~\ref{figC:ua}.
	All three algorithms have in common to use latent variables.
	In particular, these algorithms outperformed MS-MKES which benefits from a Gaussian process but works only in the mixed space, i.e., MS-MKES does not imply latent variables.
	This shows that latent variables are useful to speed up a Bayesian search for mixed problems.
	
	No clear advantage, on the average, was found for accounting for the discrete nature of the variables through constraints: LV-EGO, which ignores the link between latent variables and the discrete variables until the pre-image problem, is competitive with the best of the augmented Lagrangian ALV-EGO algorithms. 
	We hypothesize that the constraint on latent variables, by creating disconnected feasibility islands around $\Latvec{t}(\uu),~\uu \in \Uset$, makes the optimization of the acquisition criterion almost as difficult to solve as it originally was in the mixed space, therefore not allowing to fully benefit from the continuity of the $\Xset \times \Lset$ space.
	
	In our tests, the global updating of the Lagrange multipliers was always preferable to the local counterparts, ALV-EGO-gi and -ge eclipsing ALV-EGO-li and -le. 
	The ALV-EGO-gi approach, where the discrete constraint is relaxed and turned into an inequality (Equation~(\ref{eq:OP})), works better on the average than ALV-EGO-ge where the constraint is an equality. This illustrates the positive effect of the relaxation $\epsilon$, that softens the phenomenon we mentionned above where the feasible domain is broken into disconnected regions. 
	
	MS-ES is consistently less efficient than the other algorithms. It was expected, because there is no metamodel to save calls to the function. Furthermore, the sampling is done in the mixed space.
	The optimizers based on random forests have also rather poor average performances, to the exception of the 6 dimensional Hartmann function. We believe the random forests need a sufficiently large initial DoE (which happened with a higher dimension) to fruitfully guide the search.
	
\begin{figure}[H] 
		\begin{subfigure}[b]{0.5\linewidth}
			\centering
			\includegraphics[width=\textwidth]{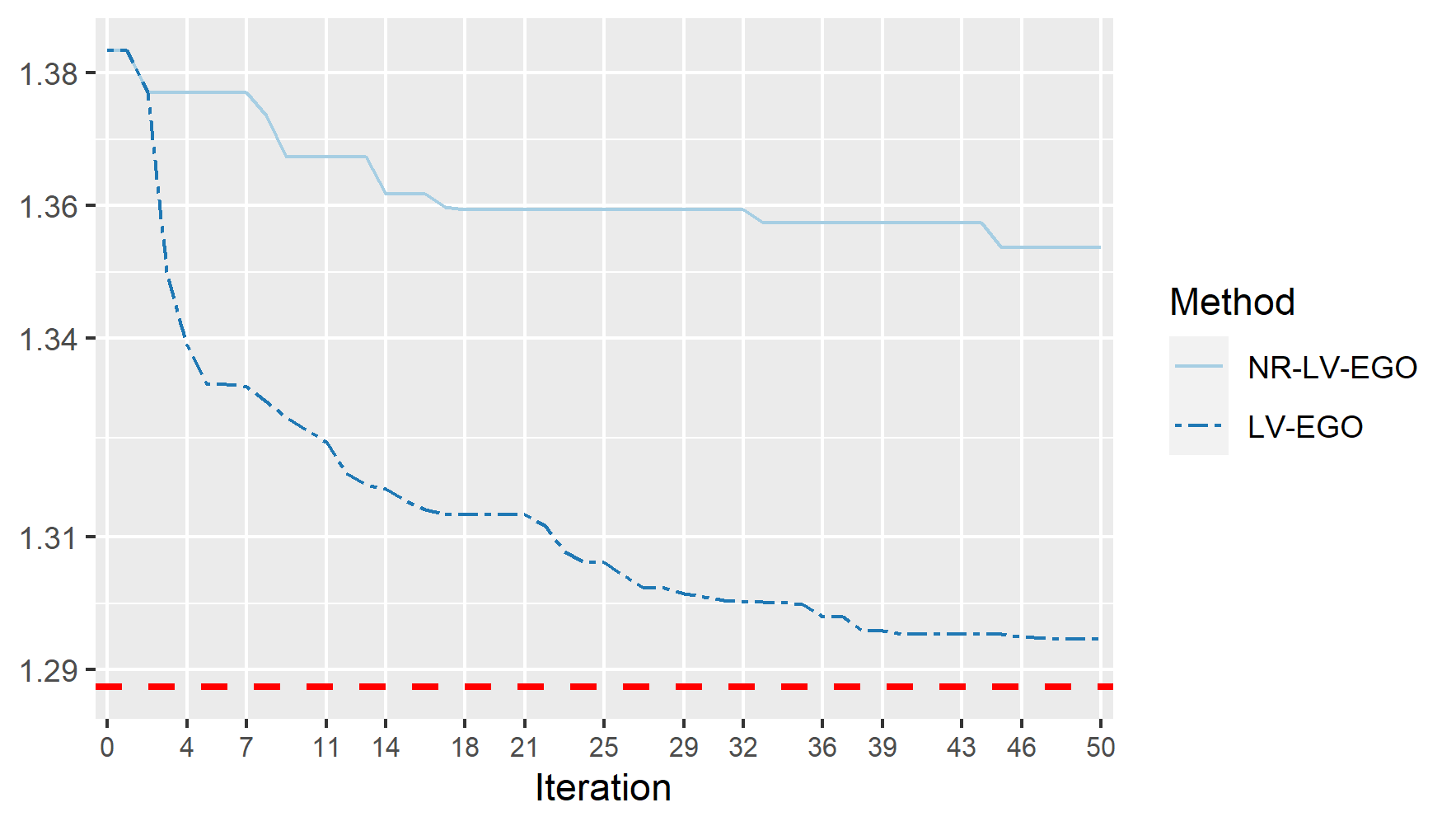}
			\caption{Median objective function} 
			\label{figlvc1:ua} 
			\vspace{4ex}
		\end{subfigure}
		\begin{subfigure}[b]{0.5\linewidth}
			\centering
			\includegraphics[width=\textwidth]{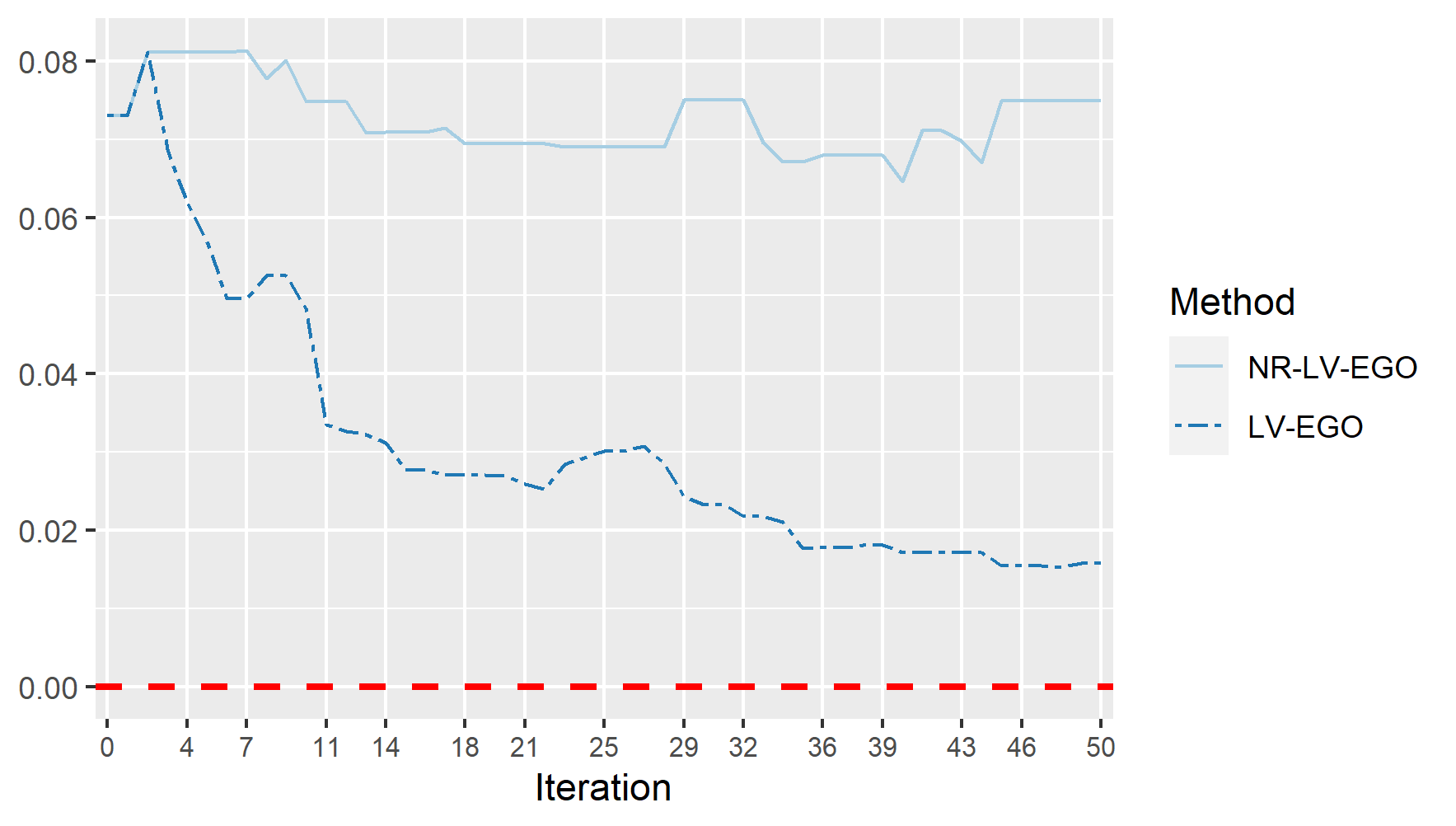}
			\caption{Interquartile of the objective functions} 
			\label{figlvc1:ub} 
			\vspace{4ex}
		\end{subfigure} 
		\begin{subfigure}[b]{0.5\linewidth}
			\centering
			\includegraphics[width=\textwidth]{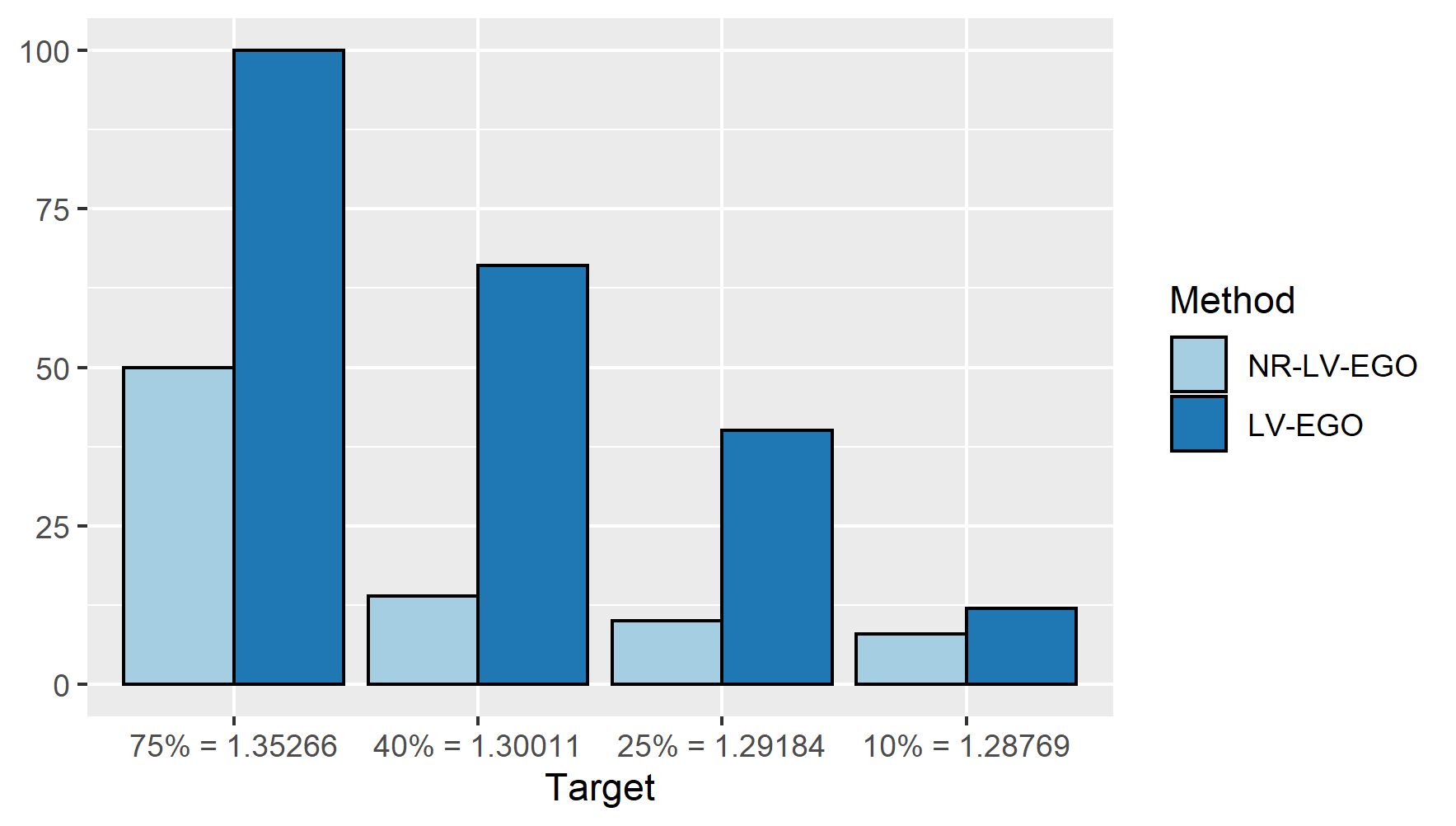}
			\caption{Success rate} 
			\label{figlvc1:uc} 
		\end{subfigure}
		\begin{subfigure}[b]{0.5\linewidth}
			\centering
			\includegraphics[width=\textwidth]{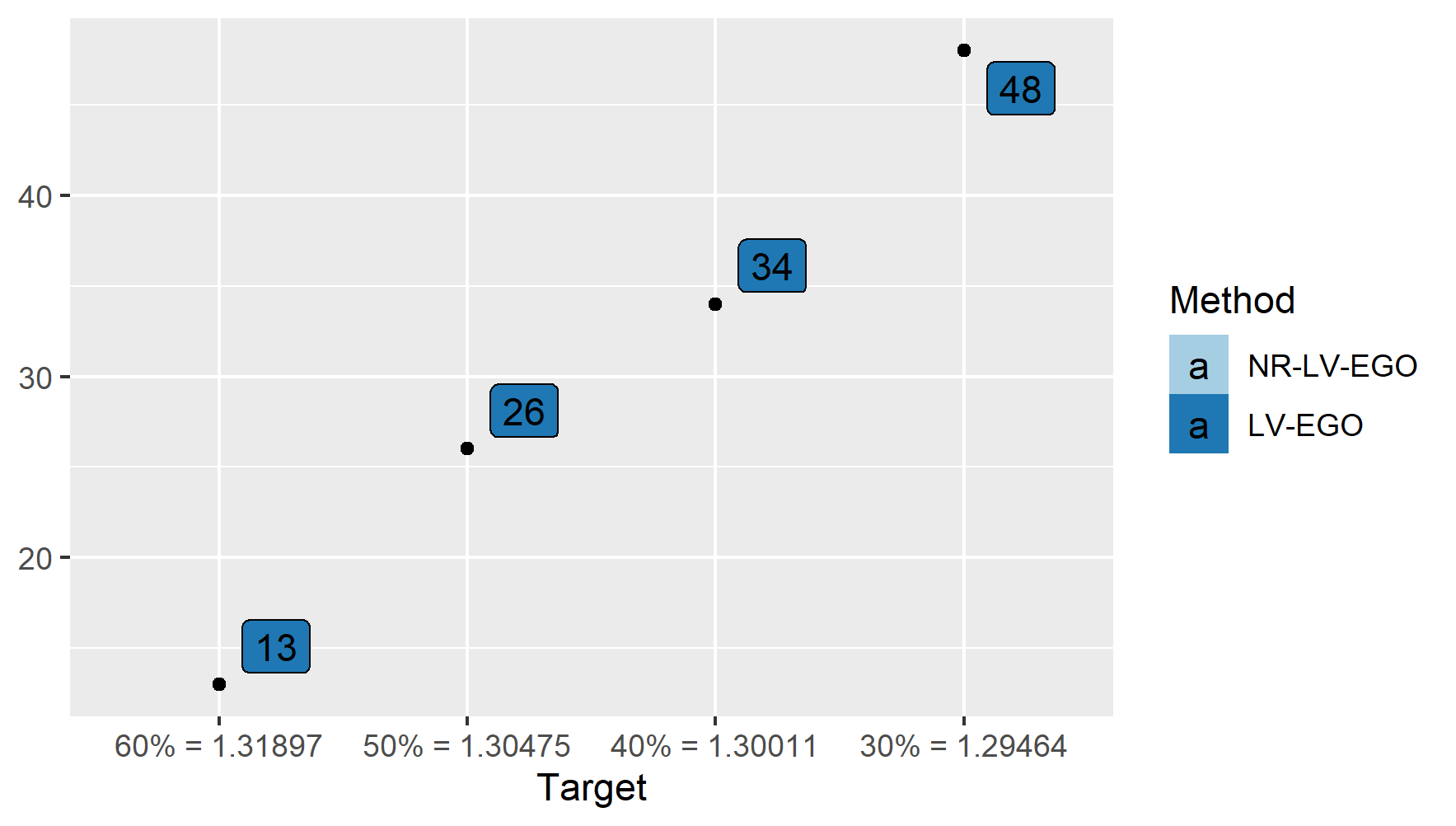}
			\caption{Iteration to median success} 
			\label{figlvc1:ud} 
		\end{subfigure} 
		\caption{Comparison of LV-EGO with and without (NR-LV-EGO) a repeated estimation of the latent variables. Results for the beam design application. }
		\label{figlvc1} 
	\end{figure}

	As a final comment, we discuss the necessity of re-estimating the latent variables at each iteration.
	The estimation of the latent variables has an important numerical 
	cost of about $\nlat t^3 \sum_{i=1}^\ndis \nlev{i}$ operations at each iteration $t$ (cf. Table \ref{tab:numComplexities}).
	It was repeated at each iteration in the algorithms with latent variables considered so far.
	In the experiment reported in Figure \ref{figlvc1}, a version of the LV-EGO algorithm is considered where the latent variables 
	are estimated once only, with the initial DoE, yielding the NR-LV-EGO algorithm (for Non Repeated estimation of $\Latvec{}()$).

	As can be seen in Figure \ref{figlvc1} when comparing LV-EGO with NR-LV-EGO, 
	the re-estimation of the latent variables at each iteration, as implemented in the LV-EGO algorithm and its ALV-EGO variants, {considerably improves} its performance. 
	An accompanying result is the visualization of the correlation matrix of the discrete variable provided in Figure \ref{figlv0}, where one notices that the correlation (hence the latent variables) evolves in time. 
	Our experiments indicate that this evolution is beneficial to the optimization efficiency.
	
	\section{Conclusions and perspectives}
	\label{sec:Concl}
	This work has investigated five Bayesian optimization approaches to small and medium size mixed problems that hinged on latent variables.
	They differed in the way the coupling between the discrete variables and their relaxed pendants, the latent variables, is implemented.
	
	Algorithms involving latent variables were compared to other algorithms directly working in the mixed space and 
	were found to consistently outperform them. 
	LV-EGO and ALV-EGO-gi were more efficient (in terms of calls to the true objective function) than MS-MKES which also benefits from the Gaussian process. 
	These first results show that latent variables provide a flexible way to handle mixed problems where the total number of levels and of variables is less or equal to about 10 variables and 10 levels in total.
	
	Accounting for the discrete nature of some variables through a constraint during the relaxed optimization with augmented Lagrangians was not clearly found to further increase the performance of the search as LV-EGO competed equally and even sometimes outperformed the ALV versions of the algorithms.
	It was also observed that expressing the discreteness as an inequality constraint by adding a tolerance was a better option than expressing it as an equality.
	The global updating strategy of the Lagrange multipliers, which to the best of our knowledge is original, improved over the more common local updating schemes.
	Finally, the random forests metamodels did not do as well as the Gaussian processes, whether in their continuous or mixed forms, within the Bayesian optimization algorithm.
	
	Our study needs to be completed in three ways. 
	To fully leverage on the continuous latent space, the gradient of the acquisition function should be analytically calculated and used to guide its maximization.
	The implementation we proposed creates more latent variables than there are discrete levels, which limits its application to about 10 levels. This limitation can be overcome with under-parameterized kernels {based on groups \cite{roustant_group_2020} or warping techniques \cite{qian2008gaussian,kergp}}.
	Mixed Bayesian optimization through latent variables would also gain in credibility if the convergence results of EGO were generalized to it.
	
	\subsection*{Data availability}
	The source code of this work will be made available upon request to the corresponding author.
	
	\bibliographystyle{plainnat}
	\bibliography{lv_ego_BO_feb_2022}

\begin{thebibliography}{33}
\providecommand{\natexlab}[1]{#1}
\providecommand{\url}[1]{\texttt{#1}}
\expandafter\ifx\csname urlstyle\endcsname\relax
  \providecommand{\doi}[1]{doi: #1}\else
  \providecommand{\doi}{doi: \begingroup \urlstyle{rm}\Url}\fi

\bibitem[Audet and Dennis~Jr(2001)]{audet2001pattern}
Charles Audet and John~E Dennis~Jr.
\newblock Pattern search algorithms for mixed variable programming.
\newblock \emph{SIAM Journal on Optimization}, 11\penalty0 (3):\penalty0
  573--594, 2001.

\bibitem[Bartz-Beielstein et~al.(2019)Bartz-Beielstein, Filipi{\v{c}},
  Koro{\v{s}}ec, and Talbi]{Bartz-Beielstein2019}
T.~Bartz-Beielstein, B.~Filipi{\v{c}}, P.~Koro{\v{s}}ec, and E.G. Talbi.
\newblock \emph{High-Performance Simulation-Based Optimization}.
\newblock Studies in Computational Intelligence. Springer International
  Publishing, 2019.
\newblock ISBN 9783030187644.
\newblock URL \url{https://books.google.fr/books?id=8yGbDwAAQBAJ}.

\bibitem[Bartz-Beielstein and Zaefferer(2017)]{bartz2017model}
Thomas Bartz-Beielstein and Martin Zaefferer.
\newblock Model-based methods for continuous and discrete global optimization.
\newblock \emph{Applied Soft Computing}, 55:\penalty0 154--167, 2017.

\bibitem[Belotti et~al.(2013)Belotti, Kirches, Leyffer, Linderoth, Luedtke, and
  Mahajan]{belotti2013mixed}
Pietro Belotti, Christian Kirches, Sven Leyffer, Jeff Linderoth, James Luedtke,
  and Ashutosh Mahajan.
\newblock Mixed-integer nonlinear optimization.
\newblock \emph{Acta Numerica}, 22:\penalty0 1--131, 2013.

\bibitem[Bischl et~al.(2018)Bischl, Richter, Bossek, Horn, Thomas, and
  Lang]{mlrMBO}
Bernd Bischl, Jakob Richter, Jakob Bossek, Daniel Horn, Janek Thomas, and
  Michel Lang.
\newblock {mlrMBO}: A modular framework for model-based optimization of
  expensive black-box functions, 2018.

\bibitem[Cao et~al.(2000)Cao, Jiang, and Wu]{cao2000evolutionary}
YJ~Cao, L~Jiang, and QH~Wu.
\newblock An evolutionary programming approach to mixed-variable optimization
  problems.
\newblock \emph{Applied Mathematical Modelling}, 24\penalty0 (12):\penalty0
  931--942, 2000.

\bibitem[Deville et~al.(2017--2021)Deville, Ginsbourger, Roustant, and
  Durrande]{kergp}
Yves Deville, David Ginsbourger, Olivier Roustant, and Nicolas Durrande.
\newblock kergp.
\newblock \url{https://cran.r-project.org/package=kergp}, 2017--2021.

\bibitem[Emmerich et~al.(2008)Emmerich, Zhang, Li, Flesch, and
  Lucas]{emmerich2008mixed}
Michael Emmerich, A~Zhang, R~Li, I~Flesch, and Peter~J. Lucas.
\newblock Mixed-integer bayesian optimization utilizing a-priori knowledge on
  parameter dependences.
\newblock \emph{Journal of Physical Chemistry A - J PHYS CHEM A}, pages 65--72,
  01 2008.

\bibitem[Frazier(2018)]{frazier_tutorial_2018}
Peter~I. Frazier.
\newblock A {Tutorial} on {Bayesian} {Optimization}.
\newblock \emph{arXiv e-prints}, page arXiv:1807.02811, July 2018.

\bibitem[Hestenes(1969)]{hestenes1969multiplier}
Magnus~R Hestenes.
\newblock Multiplier and gradient methods.
\newblock \emph{Journal of optimization theory and applications}, 4\penalty0
  (5):\penalty0 303--320, 1969.

\bibitem[Hutter et~al.(2011)Hutter, Hoos, and
  Leyton-Brown]{hutter2011sequential}
Frank Hutter, Holger~H Hoos, and Kevin Leyton-Brown.
\newblock Sequential model-based optimization for general algorithm
  configuration.
\newblock In \emph{International Conference on Learning and Intelligent
  Optimization}, pages 507--523. Springer, 2011.

\bibitem[Jones et~al.(1998)Jones, Schonlau, and Welch]{Jones1998}
Donald~R. Jones, Matthias Schonlau, and William~J. Welch.
\newblock Efficient global optimization of expensive black-box functions.
\newblock \emph{Journal of Global Optimization}, 13\penalty0 (4):\penalty0
  455--492, Dec 1998.
\newblock ISSN 1573-2916.
\newblock \doi{10.1023/A:1008306431147}.
\newblock URL \url{https://doi.org/10.1023/A:1008306431147}.

\bibitem[Le~Riche and Guyon(2002)]{leriche_dual_evo_2002}
Rodolphe Le~Riche and Fr\'ed\'eric Guyon.
\newblock Dual evolutionary optimization.
\newblock \emph{Lecture Notes in Computer Science}, \penalty0 (2310):\penalty0
  281--294, 2002.
\newblock selected papers of the 5th Int. Conf. Evolution Artificielle.

\bibitem[Le~Riche and Picheny(2021)]{leriche2021revisiting}
Rodolphe Le~Riche and Victor Picheny.
\newblock Revisiting bayesian optimization in the light of the coco benchmark.
\newblock \emph{Structural and MultiDisciplinary Optimization}, 2021.
\newblock to appear.

\bibitem[Li et~al.(2013)Li, Emmerich, Eggermont, B{\"a}ck, Sch{\"u}tz,
  Dijkstra, and Reiber]{li2013mixed}
Rui Li, Michael~TM Emmerich, Jeroen Eggermont, Thomas B{\"a}ck, Martin
  Sch{\"u}tz, Jouke Dijkstra, and Johan~HC Reiber.
\newblock Mixed integer evolution strategies for parameter optimization.
\newblock \emph{Evolutionary computation}, 21\penalty0 (1):\penalty0 29--64,
  2013.

\bibitem[Lin et~al.(2018)Lin, Liu, Chen, and Zhang]{lin_hybrid_2018}
Ying Lin, Yu~Liu, Wei-Neng Chen, and Jun Zhang.
\newblock A hybrid differential evolution algorithm for mixed-variable
  optimization problems.
\newblock \emph{Information Sciences}, 466:\penalty0 170--188, 2018.
\newblock ISSN 00200255.
\newblock \doi{10.1016/j.ins.2018.07.035}.
\newblock URL
  \url{https://linkinghub.elsevier.com/retrieve/pii/S0020025516318163}.

\bibitem[Minoux(1986)]{minoux1986mathematical}
M.~Minoux.
\newblock \emph{Mathematical Programming: Theory and Algorithms}.
\newblock A Wiley-Interscience publication. Wiley, 1986.
\newblock ISBN 9780471901709.
\newblock URL \url{https://books.google.fr/books?id=5kDvAAAAMAAJ}.
\newblock translated by Vajda, S.

\bibitem[Mockus et~al.(1978)Mockus, Tiesis, and
  Zilinskas]{mockus1978application}
Jonas Mockus, Vytautas Tiesis, and Antanas Zilinskas.
\newblock The application of bayesian methods for seeking the extremum.
\newblock \emph{Towards global optimization}, 2\penalty0 (117-129):\penalty0 2,
  1978.

\bibitem[Nocedal and Wright(2006)]{nocedal_numerical_2006}
Jorge Nocedal and Stephen~J. Wright.
\newblock \emph{Numerical optimization}.
\newblock Springer series in operations research. Springer, New York, 2nd ed
  edition, 2006.
\newblock ISBN 978-0-387-30303-1.
\newblock OCLC: ocm68629100.

\bibitem[Ocenasek and Schwarz(2002)]{ocenasek2002estimation}
Jiff Ocenasek and Josef Schwarz.
\newblock Estimation of distribution algorithm for mixed continuous-discrete
  optimization problems.
\newblock In \emph{2nd Euro-International Symposium on Computational
  Intelligence}, pages 227--232. IOS Press Kosice, Slovakia, 2002.

\bibitem[Pelamatti et~al.(2019)Pelamatti, Brevault, Balesdent, Talbi, and
  Guerin]{pelamatti_efficient_2019}
Julien Pelamatti, Loic Brevault, Mathieu Balesdent, El-Ghazali Talbi, and
  Yannick Guerin.
\newblock Efficient global optimization of constrained mixed variable problems.
\newblock \emph{Journal of Global Optimization}, 73\penalty0 (3):\penalty0
  583--613, 2019.
\newblock ISSN 0925-5001, 1573-2916.
\newblock \doi{10.1007/s10898-018-0715-1}.
\newblock URL \url{http://link.springer.com/10.1007/s10898-018-0715-1}.

\bibitem[Picheny et~al.(2016)Picheny, Gramacy, Wild, and
  Le~Digabel]{picheny_bayesian_2016}
Victor Picheny, Robert~B Gramacy, Stefan Wild, and Sebastien Le~Digabel.
\newblock Bayesian optimization under mixed constraints with a slack-variable
  augmented {L}agrangian.
\newblock In D.~Lee, M.~Sugiyama, U.~Luxburg, I.~Guyon, and R.~Garnett,
  editors, \emph{Advances in Neural Information Processing Systems}, volume~29.
  Curran Associates, Inc., 2016.
\newblock URL
  \url{https://proceedings.neurips.cc/paper/2016/file/31839b036f63806cba3f47b93af8ccb5-Paper.pdf}.

\bibitem[Powell(1994)]{Powell1994}
M.~J.~D. Powell.
\newblock \emph{A Direct Search Optimization Method That Models the Objective
  and Constraint Functions by Linear Interpolation}, pages 51--67.
\newblock Springer Netherlands, Dordrecht, 1994.
\newblock ISBN 978-94-015-8330-5.
\newblock \doi{10.1007/978-94-015-8330-5_4}.
\newblock URL \url{https://doi.org/10.1007/978-94-015-8330-5_4}.

\bibitem[Rockafellar(1993)]{rockafellar_lagrange_1993}
R.~Tyrrell Rockafellar.
\newblock Lagrange {Multipliers} and {Optimality}.
\newblock \emph{SIAM Review}, 35\penalty0 (2):\penalty0 183--238, 1993.
\newblock URL \url{http://www.jstor.org/stable/2133143}.

\bibitem[Roustant et~al.(2020)Roustant, Padonou, Deville, Clément, Perrin,
  Giorla, and Wynn]{roustant_group_2020}
Olivier Roustant, Espéran Padonou, Yves Deville, Aloïs Clément, Guillaume
  Perrin, Jean Giorla, and Henry Wynn.
\newblock Group kernels for gaussian process metamodels with categorical
  inputs.
\newblock \emph{SIAM/ASA Journal on Uncertainty Quantification}, 8\penalty0
  (2):\penalty0 775--806, 2020.
\newblock \doi{10.1137/18M1209386}.
\newblock URL \url{https://doi.org/10.1137/18M1209386}.

\bibitem[Thi et~al.(2019)Thi, Le, and Dinh]{Thi2019}
H.A.L. Thi, H.M. Le, and T.P. Dinh.
\newblock \emph{Optimization of Complex Systems: Theory, Models, Algorithms and
  Applications}.
\newblock Advances in Intelligent Systems and Computing. Springer International
  Publishing, 2019.
\newblock ISBN 9783030218034.
\newblock URL \url{https://books.google.fr/books?id=R46dDwAAQBAJ}.

\bibitem[Vazquez and Bect(2010)]{vazquez_convergence_2010}
Emmanuel Vazquez and Julien Bect.
\newblock Convergence properties of the expected improvement algorithm with
  fixed mean and covariance functions.
\newblock \emph{Journal of Statistical Planning and Inference}, 140\penalty0
  (11):\penalty0 3088--3095, 2010.
\newblock ISSN 03783758.
\newblock \doi{10.1016/j.jspi.2010.04.018}.
\newblock URL
  \url{https://linkinghub.elsevier.com/retrieve/pii/S0378375810001850}.

\bibitem[Wang et~al.(2016)Wang, Hutter, Zoghi, Matheson, and
  de~Feitas]{wang2016bayesian}
Ziyu Wang, Frank Hutter, Masrour Zoghi, David Matheson, and Nando de~Feitas.
\newblock Bayesian optimization in a billion dimensions via random embeddings.
\newblock \emph{Journal of Artificial Intelligence Research}, 55:\penalty0
  361--387, 2016.

\bibitem[Wilson et~al.(2018)Wilson, Hutter, and
  Deisenroth]{wilson_maximizing_2018}
James~T. Wilson, Frank Hutter, and Marc~Peter Deisenroth.
\newblock Maximizing acquisition functions for bayesian optimization.
\newblock In \emph{Proceedings of the 32nd International Conference on Neural
  Information Processing Systems}, NIPS’18, page 9906–9917, Red Hook, NY,
  USA, 2018. Curran Associates Inc.

\bibitem[Zaefferer(2014--2021)]{cego}
Martin Zaefferer.
\newblock {CEGO}.
\newblock \url{https://cran.r-project.org/package=CEGO}, 2014--2021.

\bibitem[Zhang et~al.(2019)Zhang, Tao, Chen, and Apley]{zhang_latent_2019}
Yichi Zhang, Siyu Tao, Wei Chen, and Daniel~W. Apley.
\newblock A {Latent} {Variable} {Approach} to {Gaussian} {Process} {Modeling}
  with {Qualitative} and {Quantitative} {Factors}.
\newblock \emph{Technometrics}, pages 1--12, 2019.
\newblock ISSN 0040-1706, 1537-2723.
\newblock \doi{10.1080/00401706.2019.1638834}.
\newblock URL
  \url{https://www.tandfonline.com/doi/full/10.1080/00401706.2019.1638834}.

\bibitem[Zhang et~al.(2020)Zhang, Apley, and Chen]{zhang_bayesian_2020}
Yichi Zhang, Daniel~W. Apley, and Wei Chen.
\newblock Bayesian {Optimization} for {Materials} {Design} with {Mixed}
  {Quantitative} and {Qualitative} {Variables}.
\newblock \emph{Scientific Reports}, 10\penalty0 (1), December 2020.
\newblock ISSN 2045-2322.
\newblock \doi{10.1038/s41598-020-60652-9}.
\newblock URL \url{http://www.nature.com/articles/s41598-020-60652-9}.

\bibitem[Zuniga and Sinoquet(2020)]{Zuniga_Sinoquet_mixedEGO}
Miguel~Munoz Zuniga and Delphine Sinoquet.
\newblock Global optimization for mixed categorical-continuous variables based
  on gaussian process models with a randomized categorical space exploration
  step.
\newblock \emph{INFOR: Information Systems and Operational Research},
  58\penalty0 (2):\penalty0 310--341, 2020.
\newblock \doi{10.1080/03155986.2020.1730677}.
\newblock URL \url{https://doi.org/10.1080/03155986.2020.1730677}.

\end{thebibliography}
	
	\subsection*{Acknowledgments}	
	This work was supported in part by the OQUAIDO research chair in applied mathematics and by the CIROQUO consortium.
	
	\subsection*{Authors' contributions}
	The kernels with latent variables were developed jointly by O. Roustant, G. Perrin and J. Cuesta-Ramirez.
	The Bayesian optimization formulation was developed jointly by R. Le Riche, J. Cuesta-Ramirez, O. Roustant, G. Perrin and C. Durantin.
	The augmented Lagrangians schemes were developped jointly by R. Le Riche and J. Cuesta-Ramirez.
	The test cases were proposed by G. Perrin, C. Durantin, A. Glière and J. Cuesta-Ramirez.
	J. Cuesta-Ramirez did the code implementation. 
	All authors reviewed the manuscript.
	
	\subsection*{Competing interests}
	The authors declare no competing interests.

	\appendix
	
	\section{Complements on the augmented Lagrangians}
	\label{sec-AL}
	\subsection*{Case of an equality constraint}
\begin{figure}
		\begin{center}
			\includegraphics[width=0.75\textwidth]{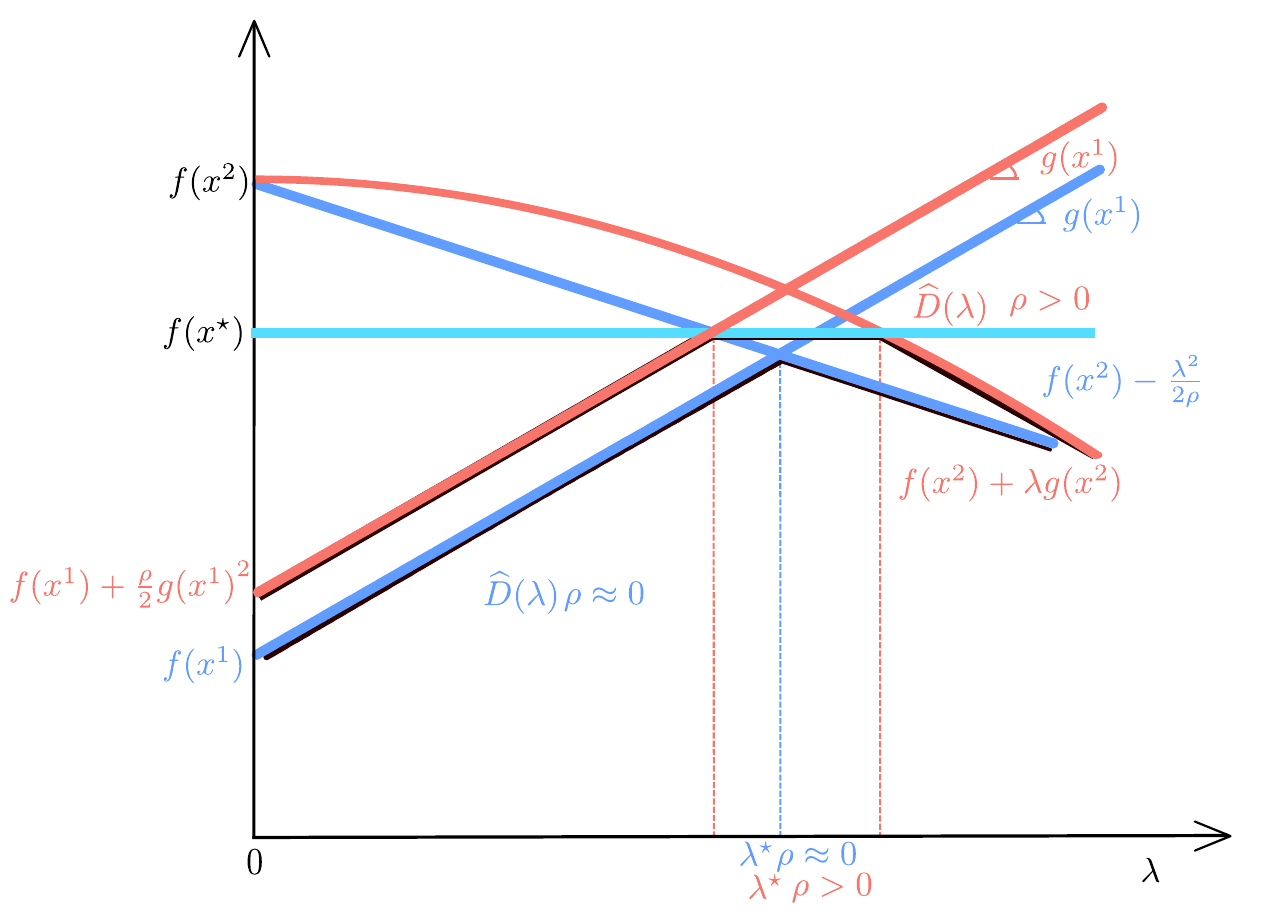}
		\end{center}
		\caption{Sketch of Rockafellar's augmented Lagrangian for $\rho\approx0$ in blue and $\rho>0$ in red. $x^1$ is infeasible, $x^2$ feasible (and $g(x^2) < -\lambda/\rho$) and $x^\star$ is an optimum with $g(x^\star)=0$. 
The black highlighted curves are the approximation to the dual function, $\widehat{\dualF}(\lambda)$ for $\Xdoe = \{x^1,x^2,x^\star\}$, for $\rho\approx0$ and $\rho>0$. There is no saddle point and a duality gap with the blue set of curves in that 
			$x^\star \notin \arg \min_x L_A(x;\lambda^\star,\rho\approx0)$ and 
			$\widehat{\dualF}(\lambda^\star) = min_x L_A(x;\lambda^\star,\rho\approx0) < L_A(x^\star;\lambda^\star,\rho \approx0)$, i.e., minimizing the augmented Lagrangian does not lead to the result of the problem. However, by increasing $\rho$, it is visible that the $y$-intercept of the infeasible points increase so that one always reaches a state where $x^\star = \arg \min_x L_A(x;\lambda^\star,\rho)$ as in the red set of curves. A similar illustration can be done with the augmented Lagrangian with equality constraint: $f(x)+\rho/2 h^2(x)$ is the $y$-intercept and $h(x)$ is the slope of the augmented Lagrangian associated to $x$. The main difference is that all points contribute linearly in terms of $\lambda$ to $L_A(x;\lambda,\rho)$.  
		}
		\label{fig:ALg}
	\end{figure}

	Let us first consider an optimization problem with an equality constraint, 
	\begin{equation}
	\left\{
	\begin{array}{l}
	\min_{x \in \mathcal X} f(x) \\
	\text{such that } h(x) = 0
	\end{array}
	\right.
	\label{eq:OPh}
	\end{equation}
	At this point, $f()$ and $h()$ are very general functions on a $d$-dimensional general set $\mathcal X$. 
	We only require that $\mathcal X$ is not empty, that $f()$ and $h()$ are bounded, and that there is at least one solution to (\ref{eq:OPh}), $x^\star \in \mathcal X$, which can be attained.
	$f()$ and $h()$ are not necessarily continuous, a fortiori not necessarily differentiable. 
	With respect to the main body of the article, the notations are simplified in this Section: $\mathcal X$ stands for the cartesian product of \Xset and \Lset, $f(x)$ generalizes $-\log(1+\EI[t]{\xx,\rlatvec})$ and $h(x)$ corresponds to $g^{(t)}(\rlatvec)$ when $\epsilon=0$.
	Note that $g^{(t)}()$, being made of the minimum distance to a discrete set of points (cf. Eq.~(\ref{eq:OP})), is not differentiable.
	$g^{(t)}()$ is the only constraint in the article. This appendix considers one constraint too, but all the results given readily generalize to many constraints by replacing the products by vector scalar products.\\
	Problem~(\ref{eq:OPh}) can be equivalently reformulated as 
	\begin{equation}
	\left\{
	\begin{array}{l}
	\min_{x \in \mathcal X} f(x) + \frac{1}{2}\rho h^2(x)\\
	\text{such that } h(x) = 0
	\end{array}
	\right.
	\label{eq:OPhrho}
	\end{equation}
	where $\rho \ge 0$ is a penalty parameter.
	The two above formulations have the same solution $x^\star$ and the same value of optimal objective function since $x^\star$ is feasible, $h(x^\star)=0$, therefore $f(x^\star) = f(x^\star) + \frac{1}{2}\rho h^2(x^\star)$.
	However, as proved in \cite{minoux1986mathematical} and sketched in Figure~\ref{fig:ALg},
	there is always a positive lower bound on the penalty parameters, $\rho\ge \rho^\star \ge 0$, such that Problem~(\ref{eq:OPhrho}) can be equivalently solved through the dual formulation, 
	\begin{equation}
	\begin{split}
	& \max_{\lambda \in \IR} \dualF(\lambda,\rho) \\
	& \text{where } \dualF(\lambda,\rho) = \min_{x \in \mathcal X} L_A(x;\lambda,\rho) \\
	& \text{and } L_A(x;\lambda,\rho) = f(x) + \lambda h(x) + \frac{1}{2} \rho h^2(x)
	\end{split}
	\label{eq:OPdual}
	\end{equation}
	In this way, the augmented Lagrangian of \cite{hestenes1969multiplier} is the classical Lagrangian of the penalized problem~(\ref{eq:OPhrho}). 
	We write $\lambda^\star,\rho^\star$ a solution to (\ref{eq:OPdual}).
	$\dualF(\lambda,\rho)$ is the lower front of all augmented Lagrangians for varying $x$ at a given $\lambda,\rho$. 
	The ``global dual'' update of $(\lambda,\rho)$ comes from the resolution of (\ref{eq:OPdual}) where the set $\mathcal X$ is approximated by the finite subset of samples $\Xdoe$. \\
	Let us denote
	\begin{equation}
	x(\lambda,\rho) = \arg \min_{x \in \mathcal X} L_A(x;\lambda,\rho)
	\label{eq:xlambda}
	\end{equation}
	a solution at given multiplier and penalty parameter. 
	The function $\dualF(\lambda,\rho)$ is concave in $\lambda$ and $\rho$ and $h(x(\lambda,\rho))$ is a subgradient with respect to $\lambda$ \cite{minoux1986mathematical}. 
	This is at the root of updating strategies that we called ``local dual'' earlier and which consist in a gradient step in the dual space, 
	\begin{equation}
	\lambda_{t+1} = \lambda_t + \alpha \partial_\lambda \dualF(\lambda_t,\rho_t) = \lambda_t + \alpha h(x(\lambda_t,\rho_t)) \quad ,
	\label{eq:lambdaUpdateGrad}
	\end{equation}
	where $\alpha > 0$ is a step size factor.

	More specific update strategies such as those given in \cite{nocedal_numerical_2006,picheny_bayesian_2016} stem from the Karush Kuhn and Tucker (KKT) optimality conditions and require the additional assumption that $\mathcal X \in \IR^d$ and $f()$ and $h()$ are differentiable. At $x^\star$, since $h(x^\star)=0$ and $\lambda^{KKT}$ being the KKT multiplier\footnote{
		The Lagrange multiplier that maximizes the dual function is equal to the KKT multiplier only when the functions are differentiable, the constraints qualification conditions apply, and there is a saddle point i.e., $min_x max_{\lambda} L_A(x;\lambda,\rho) = max_{\lambda} min_x L_A(x;\lambda,\rho)$.
	}, one has
	\begin{eqnarray}
	\nabla f(x^\star) + \rho h(x^\star) \nabla h(x^\star) + \lambda^{KKT} \nabla h(x^\star) & =& 0 \nonumber \\
	\Rightarrow \qquad \nabla f(x^\star) + \lambda^{KKT} \nabla h(x^\star) & = & 0 \label{eq:kkt}
	\end{eqnarray}
	At iteration $t$, the necessary conditions for $x^t = x(\lambda_t,\rho_t)$ to be the minimum of $L_A(;\lambda_t,\rho_t)$ are
	\begin{equation}
	\nabla f(x^t) + \rho_t h(x^t) \nabla h(x^t) + \lambda_t \nabla h(x^t) = 
	\nabla f(x^t) + (\rho_t h(x^t) + \lambda_t) \nabla h(x^t) = 0  
	\label{eq:xtcond}
	\end{equation}
	Comparing equations~(\ref{eq:kkt}) and (\ref{eq:xtcond}), $x^t$ can be driven to $x^\star$ if 
	\begin{equation}
	\lambda_{t+1} = \lambda_t + \rho_t h(x^t)
	\label{eq:lambdaUpdateH}
	\end{equation}
	The updates (\ref{eq:lambdaUpdateGrad}) and (\ref{eq:lambdaUpdateH}) have the same form,  (\ref{eq:lambdaUpdateH}) is more restrictive since the KKT conditions must apply but the step size is known.
	
	\vskip\baselineskip
	The equality constraint of the article (Equation~(\ref{eq:OP}) with $\epsilon=0$) is a minimum over distances. 
	It has the additional feature that it is always positive or null, $\forall x \in \mathcal X~,~h(x) \ge 0$. 
	Because of this, if $h$ is locally differentiable around $x^\star$, $\nabla h(x^\star) = 0$ since $h$ has a minimum at $x^\star$. The constraint qualification condition is not satisfied ($\nabla h(x^\star)$ does not span a non-empty set) and the KKT conditions do not apply.
	Another consequence is that the optimal Lagrange multiplier must be positive and the search for $\lambda$ can be written $\max_{\lambda \ge 0} \dualF(\lambda,\rho)$ in Problem~(\ref{eq:OPdual}), as in Problem~(\ref{eq:dualFormulation}).\\
	\textit{Proof:} 
	Assume $\rho$ is large enough for Problem~(\ref{eq:OPhrho}) to have a saddle point at its optimum, $f(x^\star) \le f(x) + \rho/2 h^2(x) + \lambda^\star h(x)~,~ \forall x$ where $\lambda^\star$ is the optimum Lagrange multiplier.
	Since the optimization problem has an active constraint, there is a point $x^I$ that is infeasible, $h(x^I)>0$, and has a better objective function than the feasible solution (otherwise the constraint is useless), $f(x^I)+\frac{\rho}{2}h^2(x^I) \le f(x^\star)$. If the optimum Lagrange multiplier is negative, $\lambda^\star < 0$, $f(x^I)+\frac{\rho}{2}h^2(x^I)+\lambda^\star h(x^I) < f(x^\star)$ which contradicts the fact that $x^\star$ is a solution to the dual problem.
	$\square$

	\subsection*{Inequality constraint}
	When $\epsilon>0$, Problem~(\ref{eq:OP}) has an inequality constraint which we rewrite here more simply, 
	\begin{equation}
	\left\{
	\begin{array}{l}
	\min_{x \in \mathcal X} f(x) \\
	\text{such that } g(x) \le 0
	\end{array}
	\right.
	\label{eq:OPg}
	\end{equation}
	The considerations on augmented Lagragian done above for equality constraints readily extend to inequality constraints by introducing a slack variable, 
	\begin{equation}
	\left\{
	\begin{array}{l}
	\min_{x,s \in \mathcal X \times \IR} f(x) \\
	\text{such that } g(x)+s^2 = 0
	\end{array}
	\right.
	\label{eq:OPgs}
	\end{equation}
	and the expression for the augmented Lagrangian (\ref{eq:OPdual}) becomes
	\begin{equation}
	L_A(x,s;\lambda,\rho = f(x) + \lambda (g(x)+s^2) + \frac{1}{2} \rho (g(x)+s^2)^2
	\label{eq:ALgs}
	\end{equation}
	The minimization of $L_A()$ on the slack variable $s$ can be done analytically:
	\begin{equation*}
	\frac{\partial L_A(x,s;\lambda,\rho)}{\partial s} = 0
	~\iff~
	s^2 = - \frac{\lambda}{\rho} - g(x)
	\end{equation*}
	Since $s^2$ needs to be positive, all cases are summed up in 
	\begin{equation}
	s^2  = \max\left(0,- \frac{\lambda}{\rho} - g(x)\right)
	\label{eq:s2}
	\end{equation}
	Reinjecting the expression of $s^2$ into the augmented Lagrangian yields
	\begin{equation}
	L_A(x;\lambda,\rho) = f(x) + \frac{1}{2\rho}\left[(\max(0,\lambda+\rho g(x)))^2 - \lambda^2 \right]
	\label{eq:ALineq}
	\end{equation}
	which is equivalent to the expression of Rockafellar with the 2 cases given in Equation~(\ref{eq-AL_Rockafellar}) (recall $-\log(1+EI)$ is $f(x)$).

	The update equations for $\lambda$ are the same as those for the equality case where the slack variable $s^2$ takes its optimal value. 
	On the one hand, it is possible to solve the approximated dual problem as in (\ref{eq-dualFHat}).
	On the other hand, a step along a subgradient in the dual space can be taken,  
	\begin{eqnarray}
	\lambda_{t+1} &=&\lambda_t + \alpha (g(x^t) + s_t^2)  \nonumber \\
	\Rightarrow \lambda_{t+1} &=&\lambda_t + \alpha \left( g(x^t) + \max(0,-\frac{\lambda}{\rho} - g(x^t)) \right)  \label{eq:lambdaUpdIneqDual}
	\end{eqnarray}
	where $\alpha$ is again a positive step factor. It has the same form as Equation~(\ref{eq-lambdaUpdateKKT}).
	The update (\ref{eq-lambdaUpdateKKT}) is fully recovered from the KKT conditions as above for equalities, (\ref{eq:lambdaUpdateH}), 
	\begin{eqnarray}
	\lambda_{t+1} &=&\lambda_t + \rho (g(x^t) + s_t^2)  \nonumber \\
	\Rightarrow \lambda_{t+1} &=&\lambda_t + \rho \left( g(x^t) + \max(0,-\frac{\lambda}{\rho} - g(x^t)) \right)  \label{eq:lambdaUpdIneqDualRho}
	\end{eqnarray}
	Equations~(\ref{eq:lambdaUpdIneqDual}) and (\ref{eq:lambdaUpdIneqDualRho}) are the same but in the latest the step factor $\alpha$ is known and 
	equal to $\rho$, which comes at the additional expense of the KKT validity conditions.

\end{document}